\documentclass[12pt,a4paper]{article}
\usepackage{fullpage}
\usepackage{amscd}
\usepackage{amssymb}
\def\so{\mathfrak{so}}
\def\SO{\mathbf{SO}}
\def\sp{\mathfrak{sp}}
\def\Sp{\mathbf{Sp}}
\def\u{\mathfrak{u}}
\def\su{\mathfrak{su}}
\def\SU{\mathbf{SU}}
\def\spin{\mathfrak{spin}}
\def\Spin{\mathbf{Spin}}

\def\G{\mathbf{G}}
\def\1{\mathbf{1}}
\def\#{\sharp}

\def\Aut{\mathrm{Aut}}
\def\b{\flat}
\def\C{\mathbb{C}}
\def\Cas{\mathrm{Cas}}
\def\D{\Delta}
\def\e{\varepsilon}
\def\End{\mathrm{End}\,}
\def\g{\mathfrak{g}}
\def\h{\mathfrak{h}}
\def\H{\mathbb{H}}
\def\Hom{\mathrm{Hom}\,}
\def\id{\mathrm{id}}
\def\im{\mathrm{im}\,}
\def\ins{\,\lrcorner\,\,}
\def\L{\Lambda}
\def\l{\lambda}
\def\o{\omega}
\def\pr{\mathrm{pr}}
\def\Q{\mathcal{Q}}
\def\R{\mathbb{R}}
\def\Ric{\mathrm{Ric}}
\def\S{\mathrm{Sym}}
\def\t{\mathfrak{t}}
\def\tr{\mathrm{tr}}
\def\U{\mathcal{U}\,}
\def\W{\mathfrak{W}}
\def\Zent{\mathrm{Zent}\,}
\def\<#1,#2>{\langle\,#1,\,#2\,\rangle}

\def\rectangle(#1,#2)[#3,#4]#5{
 \multiput(#1,#2)(#3,0)2{\line(0,1){#4}}\multiput(#1,#2)(0,#4)2{\line(1,0){#3}}
 \put(#1,#2){\vbox to #4pt{\hbox to #3pt{\hfill}\vfill}}}
\def\recttext(#1,#2)[#3,#4]#5{\put(#1,#2)
 {\vbox to #4pt{\vfill\hbox to #3pt{\hss#5\hss}\vfill}}}
\def\tfrac#1#2{{\textstyle\frac{#1}{#2}}}
\def\proof{\noindent\textbf{Proof:}\quad}
\def\pfill{\par\vskip3mm plus1mm minus1mm\noindent}
\def\qed{\ensuremath{\hfill\Box}}
\newtheorem{Lemma}{Lemma}[section]
\newtheorem{Proposition}[Lemma]{Proposition}
\newtheorem{Theorem}[Lemma]{Theorem}
\newtheorem{Corollary}[Lemma]{Corollary}
\newtheorem{Definition}[Lemma]{Definition}
\newtheorem{Remark}[Lemma]{Remark}
\newtheorem{Example}[Lemma]{Example}
\begin{document}
\title{The Weitzenb{\"o}ck Machine}
\author{Uwe Semmelmann \& Gregor Weingart}
\maketitle
 \begin{abstract}
   Weitzenb\"ock formulas are an important tool in relating
   local differential geometry to global topological properties by means
   of the so--called Bochner method. In this article we give a unified
   treatment of the construction of all possible Weitzenb\"ock formulas
   for all irreducible, non--symmetric holonomy groups. The resulting
   classification is two--fold, we construct explicitly a basis of the
   space of Weitzenb\"ock formulas on the one hand and characterize
   Weitzenb\"ock formulas as eigenvectors for an explicitly known matrix
   on the other. Both classifications allow us to find tailor--suit
   Weitzenb\"ock formulas for applications like eigenvalue estimates
   or Betti number estimates.
 \end{abstract}
\tableofcontents
\section{Introduction}
 Weitzenb\"ock formulas are an important tool for linking differential
 geometry and topology of compact Riemannian manifolds. They feature
 prominently in the Bochner method, where they are used to prove the
 vanishing of Betti numbers under suitable curvature assumptions or the
 non--existence of metrics of positive scalar curvature on spin manifolds
 with non--vanishing $\hat A$--genus. Moreover they are used to proof
 eigenvalue estimates for Laplace and Dirac type operators. In these
 applications one tries to find a (positive) linear combination of
 hermitean squares $D^*D$ of first order differential operators $D$,
 which sums to an expression in the curvature only. In this approach
 one need only consider special first order differential operators $D$
 known as generalized gradients or Stein--Weiss operators, which are
 defined as projections of a covariant derivative $\nabla$. Examples
 for generalized gradients include the exterior derivative $d$ and
 its adjoint $d^*$ and the Dirac and twistor operator in spin geometry.

 In this article we present two different classifications of all possible
 linear combinations of hermitean squares $D^*D$ of generalized gradients
 $D$, which sum to pure curvature expressions, if the underlying connection
 is the Levi--Civita connection $\nabla$ of a Riemannian manifold $M$ of
 reduced holonomy $\SO(n)$, $\mathbf{U}(n)$, $\SU(n)$ or the exceptional
 holonomies $\G_2$ and $\Spin(7)$. Both classifications are interesting
 in their own right, the first describes a recursive procedure to calculate
 a generating set of Weitzenb\"ock formulas, the second classification
 provides a simple means to decide, whether a given linear combination
 of hermitean squares of Stein--Weiss operators is actually a pure
 curvature expression.

 \pfill
 In order to describe the setup of the article in more detail we recall
 that every representation $G\longrightarrow\Aut\,V$ of the holonomy group
 $G$ of a Riemannian manifold $(M,g)$ on a complex vector space $V$ defines
 a complex vector bundle $VM$ on $M$ with a covariant derivative induced
 from the Levi--Civita connection, in particular the complexified holonomy
 representation $T$ of $G$ defines the complexified tangent bundle $TM$.
 The generalized gradients on $VM$ are the parallel first order differential
 operators $T_\e$ defined as the projection of $\nabla:\,VM\longrightarrow
 TM\otimes VM$ to the parallel subbundles $V_\e M\subset TM\otimes VM$ arising
 from a decomposition $T\otimes V\,=\,\oplus_\e V_\e$ into irreducible
 subspaces. It will be convenient in this article to call every (finite)
 linear combination $\sum_\e c_\e T^*_\e T^{\phantom{*}}_\e$ of hermitean
 squares of generalized gradients a {\em Weitzenb\"ock formula}.

 Our first important observation is that the space $\W(V)$ of all Weitzenb\"ock
 formulas on a vector bundle $VM$ can be identified with the vector space
 $\End_\g(T\otimes V)$ and thus is an algebra, which is commutative for irreducible
 representations $V$. Moreover it is easy to see that the algebra $\W(V)$
 has a canonical involution, the twist $\tau:\,
 \W(V)\longrightarrow\W(V)$, such that a Weitzenb\"ock formula reduces
 to a pure curvature expression if and only if it is an eigenvector of
 $\tau$ of eigenvalue $-1$. Of course there are interesting Weitzenb\"ock
 formulas, which are eigenvectors of $\tau$ for the eigenvalue $+1$,
 perhaps the most prominent example is the connection Laplacian $\nabla^*
 \nabla$. The classical examples of Weitzenb\"ock formulas like
 $$
  \Delta\;\;=\;\;d\,d^*\;+\;d^*\,d\;\;=\;\;\nabla^*\nabla\;+\;q(R)
 $$
 the original Weitzenb\"ock formula or the Lichnerowicz--Weitzenb\"ock formula
 $$
  D^2\;\;=\;\;\nabla^*\nabla\;+\;\frac\kappa4
 $$
 reduce in this setting to the statements that $\Delta-\nabla^*\nabla$ and
 $D^2-\nabla^*\nabla$ respectively are eigenvectors of $\tau$ for the
 eigenvalue $-1$ and thus pure curvature expressions.

 Starting with the connection Laplacian $\nabla^*\nabla$, corresponding to
 $\1 \in \End_\g(T \otimes V)$, we will describe
 a recursion procedure to construct a basis of the space $\W(V_\l)$ of
 Weitzenb\"ock formulas on an irreducible vector bundle $V_\l M$ on $M$
 such that the base vectors are eigenvectors of $\tau$ with alternating
 eigenvalues $\pm1$. Interestingly this recursive procedure makes essential
 use of a second fundamental Weitzenb\"ock formula $B\in\W(V)$, the so--called
 {\em conformal weight operator}, which was considered for the first time
 in the work of Paul Gauduchon on conformal geometry \cite{pg1}.
 The details
 of this recursion procedure and the first initial elements are discussed
 for the holonomies $\SO(n)$, $\G_2$ and $\Spin(7)$ only, because the
 discussion of Weitzenb\"ock formulas in the K\"ahler holonomies $\mathbf{U}
 (n)$ and $\SU(n)$ is better done differently, whereas the hyperk\"ahler
 holonomies $\Sp(1)\,\Sp(n)$ and $\Sp(n)$ will be discussed in more detail
 in a forthcoming paper.

 Eventually we obtain a sequence of $B$--polynomials $p_i(B)$ such that
 $p_{2i}(B)$ is in the $(+1)$-- and $p_{2i-1}(B)$ is in the $(-1)$--eigenspace
 of $\tau$. If $b_\e$ are the $B$-eigenvalues on $V_\e\subset T \otimes V$,
 then the coefficient of $T^*_\e T^{\phantom{*}}_\e$ in the Weitzenb\"ock
 formula corresponding to $p_i(B)$ is given by $p_i(b_\e)$. An interesting
 feature appears for holonomy $\G_2$ and $\Spin(7)$. Here we have the
 decomposition $\Hom_\g(\L^2T,\End V) \cong \Hom_g(\g, \End V) \oplus
 \Hom_\g(\g^\perp, \End V)$ and because of the holonomy reduction any
 Weitzenb\"ock formula in the second summand has a zero curvature term.

 Finally we would like to mention that the problem of finding all
 possible Weitzenb\"ock formulas is also considered in the work
 of Y.~Homma (e.g. in \cite{homma2}). He gives a solution in
 the case of Riemannian, K\"ahlerian and HyperK\"ahlerian manifolds.
 Even if there are some similarities in the results, it seems fair
 to say that our method is completely different. In particular we
 describe an recursive procedure for obtaining the coefficients
 of Weitzenb\"ock formulas. The main difference is of course that
 we give a unified approach including the case of exceptional
 holonomies.
\section{The Holonomy Representation}\label{tables}
 For the rest of this article we will essentially restrict to irreducible
 non--symmetric holonomy algebras $\g$. Most of the statements easily
 generalize to holonomy algebras $\g$ with no symmetric irreducible factor
 in their local de Rham--decomposition, which could be called properly
 non--symmetric holonomy algebras. Some of the concepts introduced are
 certainly interesting for symmetric holonomy algebras as well, in particular
 the central idea used to find the matrix of the twist through the Recursion
 Formula \ref{rec}. Turning to irreducible non--symmetric holonomy algebras
 $\g$ leaves us with seven different cases
 \begin{equation}
  \label{ht}
  \hbox{\begin{tabular}{|l|c|l|c|}
   \hline & & & \\[-3mm]
    & algebra $\g_\R$ & holonomy representation $T_\R$
    & $T$ \\[1mm]
   \hline & & & \\[-3mm]
   general Riemannian
    & $\so_n$ & defining representation $\R^n$
    & $T$ \\[1mm]
   K\"ahler
    & $\u_n\cong i\R\oplus\su_n$ & defining representation $\C^n$
    & $\bar E\oplus E$ \\[1mm]
   Calabi--Yau
    & $\su_n$ & defining representation $\C^n$
    & $\bar E\oplus E$ \\[1mm]
   quaternionic K\"ahler
    & $\sp(1)\oplus\sp(n)$ & representation $\H^1\otimes_\H\H^n$
    & $H\otimes E$ \\[1mm]
   hyper--K\"ahler
    & $\sp(n)$ & defining representation $\H^n$
    & $\C^2\otimes E$ \\[1mm]
   exceptional $\G_2$
    & $\g_2$ & standard representation $\R^7$
    & $[7]$ \\[1mm]
   exceptional $\Spin(7)$
    & $\spin_7$ & spinor representation $\R^8$
    & $[8]$ \\[1mm]
   \hline
  \end{tabular}}
 \end{equation}
 according to a theorem of Berger, where $T$ denotes the complexified holonomy
 representation $T\,:=\,T_\R\otimes_\R\C$ endowed with the $\C$--bilinear
 extension $\<,>$ of the scalar product. For simplicity we will work with
 the complexified holonomy representation $T$ and the complexified holonomy
 algebra $\g\,:=\,\g_\R\otimes_\R\C$ throughout as well as with irreducible
 complex representations $V_\l$ of $\g$ of highest weight $\l$. Notations
 like $E$ or $H$ in the table above fix nomenclature for particularly
 important representations in special holonomy, say $E$ and $\bar E$ refer
 to the spaces of $(1,0)$-- and $(0,1)$--vectors in $T$ in the K\"ahler and
 Calabi--Yau case, while $[7]$ and $[8]$ are the standard $7$--dimensional
 representation of $\G_2$ and $8$--dimensional spinor representation
 of $\Spin(7)$ respectively. In passing we note that the complexified
 holonomy representation $T$ is not isotypical in the K\"ahler and the
 Calabi--Yau case and this is precisely the reason why these two cases
 differ significantly from the rest.

 In order to understand Weitzenb\"ock formulas or parallel second order
 differential operators it is a good idea to start with parallel first
 order differential operators usually called generalized gradients or
 Stein--Weiss operators. Their representation theoretic background is
 the decomposition of tensor products $T\otimes V$ of the holonomy
 representation $T$ with an arbitrary complex representation $V$. The
 general case immediately reduces to studying irreducible representations
 $V\,=\,V_\l$ of highest weight $\l$. In this section we will see that the
 isotypical components of $T\otimes V_\l$ are always irreducible for a
 properly non--symmetric holonomy algebra $\g$ and isomorphic to irreducible
 representations $V_{\l+\e}$ of highest weight $\l+\e$ for some weight
 $\e$ of the holonomy representation $T$. Thus the decomposition of
 $T\otimes V_\l$ is completely described by the subset of relevant
 weights $\e$:

 \begin{Definition}[Relevant Weights]
 \label{weights}\hfill\break
  A weight $\e$ of the holonomy representation $T$ is called relevant
  for an irreducible representation $V_\l$ of highest weight $\l$  if the
  irreducible representation $V_{\l+\e}$ of highest weight $\l+\e$ occurs
  in the tensor product $T\otimes V_\l$. We will write $\e\subset\l$ for a
  relevant weight $\e$ for a given irreducible representation $V_\l$.
 \end{Definition}

 \begin{Lemma}[Characterization of Relevant Weights]
 \label{relw}\hfill\break
  Consider the holonomy representation $T$ of an irreducible non--symmetric
  holonomy algebra $\g$ and an irreducible representation $V_\l$ of highest
  weight $\l$. The decomposition of the tensor product $T\otimes V_\l$ is
  multiplicity free in the sense that all irreducible subspaces are pairwise
  non--isomorphic. The complete decomposition of $T\otimes V_\l$ is thus the
  sum
  $$
   T\,\otimes\,V_\l\;\;\cong\;\;\bigoplus_{\e\subset\l}\,V_{\l+\e}
  $$
  over all relevant weights $\e$. A weight $\e\,\neq\,0$ is relevant if and
  only if $\l+\e$ is dominant. The zero weight $\e\,=\,0$ only occurs for the
  holonomy algebras $\so_n$ with $n$ odd and $\g_2$, it is relevant if $\l-
  \l_\Sigma$ or $\l-\l_T$ respectively is still dominant, where $\l_\Sigma$
  and $\l_T$ are the highest weights of the spinor representation of $\so_n$
  and the standard representation of $\g_2$.
 \end{Lemma}

 \proof
 The proof is essentially an exercise in Weyl's character formula
 \qed

 \pfill
 A particular consequence of Lemma \ref{relw} is that for sufficiently
 complicated representations $V_\l$ all weights $\e$ of the holonomy
 representation $T$ are relevant. With this motivation we will call a
 highest weight $\l$ {\em generic} if $\l + \e$ is dominant for all
 weights $\e$ of the holonomy representation $T$. A simple consideration
 shows that $\l$ is generic if and only if $\l-\rho$ is dominant, where
 $\rho$ is the half sum of positive roots or equivalently the sum of
 fundamental weights, unless we consider odd--dimensional generic holonomy
 $\g=\so_{2r+1}$ or $\g=\g_2$. In the latter holonomies the generic weights
 $\l$ must have $\l-\rho-\l_\Sigma$ or $\l-\rho-\l_T$ dominant respectively.
 In any case the number of relevant weights for the representation $V_\l$
 $$
  N(G, \l)
  \;\;:=\;\;
  \#\,\{\;\;\e\;|\;\;\e\textrm{\ is relevant for\ }\l\;\;\}
  \;\;\leq\;\;
  \dim\,T
 $$
 is bounded above by $\dim\,T$ with equality if and only if $\l$ is
 generic. In particular there are at most $\dim\,T$ summands in the
 decomposition of $T\otimes V_\l$ into irreducibles, exactly one copy
 of $V_{\l+\e}$ for every relevant weight $\e$.

 On the other hand the number $N(G,\l)$ of irreducible summands in the
 decomposition of $T\otimes V_\l$ agrees with the dimension of the algebra
 $\End_\g(\,T\otimes V_\l\,)$ of $\g$--invariant endomorphisms of $T\otimes
 V_\l$, because all isotypical components are irreducible by Lemma \ref{relw}.
 In the next section we will study the identification $\End_\g(\,T\otimes
 V_\l\,)\,=\,\Hom_\g(\,T\otimes T,\,\End\,V_\l\,)$ extensively, which allows
 us to break up $\End_\g(\,T\otimes V_\l\,)$ into interesting subspaces called
 Weitzenb\"ock classes, whose dimension can be calculated in the following
 way:

 \begin{Lemma}[Dimension of Weitzenb\"ock Classes]
 \hfill\label{generic}\break
  Let us call the space $W^\t\,:=\,\Hom_\t(\R,W)\,\subset\,W$ of elements
  of a $G$--representation $W$ invariant under a fixed Cartan subalgebra
  $\,\t\subset\g$ the zero weight space of $W$. The dimension of the zero
  weight space provides an upper bound
  $$
   \dim\,\Hom_\g(\;W,\,\End\,V_\l\;)
   \;\;\le\;\;
   \dim\,W^\t
  $$
  for the dimension of the space $\Hom_\g(\,W,\,\End\,V_\l\,)$ for an
  irreducible representation $V_\l$. For sufficiently dominant highest
  weight $\l$ in dependence on $W$ this upper bound is sharp.
 \end{Lemma}

 \noindent
 The lemma follows again from the Weyl character formula, but it is also an
 elementary consequence of Kostant's theorem \ref{kostant} formulated below.
 We will mainly use Lemma \ref{generic} for the subspaces $W_\alpha$ occuring
 in the decomposition $\,T \otimes T\,=\,\oplus\,W_\alpha$ into irreducibles.
 In the case of the holonomy algebras $\so_n,\, \g_2$ and $\spin_7$
 we have the decomposition $T\otimes T\,=\,\C\oplus\S^2_0T\oplus\g\oplus
 \g^\perp$ and the following dimensions of the zero weight spaces:
 \begin{equation}\label{wt}
  \hbox{\begin{tabular}{|l|c|c|c|c|c|}
   \hline & & & & & \\[-3mm]
    & $\dim\,T$
    & $\dim\,[\,\C\,]^\t$ & $\dim\,[\,\S^2_0T\,]^\t$
    & $\dim\,[\,\g\,]^\t$ & $\dim\,[\,\g^\perp\,]^\t$ \\[1mm]
   \hline & & & & & \\[-3mm]
    $\so_n$ & $n$ &
    $1$ & $\lfloor\frac{n-1}2\rfloor$ & $\lfloor\frac n2\rfloor$ & ---\\[1mm]
    $\g_2$   & $7$  & $1$ & $2$   & $3$   & $1$ \\[1mm]
    $\spin_7$ & $8$ & $1$ & $3$   & $3$   & $1$ \\[1mm]
   \hline
  \end{tabular}}
 \end{equation}
 Note in particular that the dimensions of the zero weight spaces sum up to
 $\dim\,T$.

 \pfill
 Although complete the decision criterion given in Lemma \ref{relw} is not
 particularly straightforward in general. At the end of this section we
 want to give a graphic interpretation of this decision criterion for all
 irreducible holonomy groups in order to simplify the task of finding the
 relevant weights. For a fixed holonomy algebra $\g$ the information
 necessary in this graphic algorithm is encoded in a single diagram
 featuring the weights of the holonomy representation $T$ and labeled
 boxes. A weight $\e$ is relevant for an irreducible representation
 $V_\l$ if and only if the highest weight $\l=\l_1\o_1+\ldots+\l_r\o_r$
 of $V_\lambda$, written as a linear combination of fundamental weights
 $\o_1,\ldots,\o_r$, satisfies all inequalities labeling the boxes containing
 $\e$. The notation introduced for the weights of the holonomy representation
 $T$ and the fundamental weights will be used throughout this article.

 To begin with let us consider even dimensional Riemannian geometry
 with generic holo\-nomy $\g\,=\,\so_{2r},\,r\geq 2$. In this case the
 holonomy representation $T$ is the defining representation, whose weights
 $\pm\e_1,\,\pm\e_2,\,\ldots,\,\pm\e_r$ form an orthonormal basis for a
 suitable scalar product $\<,>$ on the dual $\t^*$ of the maximal torus.
 The ordering of weights can be chosen in such a way that the fundamental
 weights $\omega_1,\,\ldots,\,\omega_r$ are given by:
 $$
  \begin{array}{lclclcl}
   \o_1       &=& \e_1 &\qquad&
   \pm\e_1    &=& \pm\,\o_1 \\
   \o_2       &=& \e_1+\e_2 &&
   \pm\e_2    &=& \pm(\o_2-\o_1) \\
   \;\;\vdots & & \quad\vdots && \;\;\vdots && \quad\vdots \\
   \o_{r-2}   &=& \e_1+\ldots+\e_{r-2} &&
   \pm\e_{r-2}&=& \pm(\o_{r-2}-\o_{r-3}) \\
   \o_{r-1}   &=& \frac12(\e_1+\ldots+\e_{r-1}+\e_r) &&
   \pm\e_{r-1}&=& \pm(\o_{r-1}+\o_r-\o_{r-2}) \\
   \o_r       &=& \frac12(\e_1+\ldots+\e_{r-1}-\e_r) &&
   \pm\e_r    &=& \pm(\o_{r-1}-\o_r)
  \end{array}
 $$
 Every dominant integral weight of $\so_{2r}$ can be written
 $\l\,=\,\l_1\o_1+\ldots+\l_r\o_r$ with natural numbers $\l_1,
 \,\ldots,\,\l_r\geq 0$ and the criterion of Lemma \ref{relw} becomes:
 \begin{center}\begin{picture}(430,130)(0,0)
  \linethickness{.8pt}
  \rectangle(  0, 90)[ 50, 20]{Gray}
  \recttext (  0, 90)[ 50, 20]{$+\e_1$}
  \linethickness{.2pt}
  \rectangle( 50, 70)[ 50, 40]{White}
  \recttext ( 50,110)[ 50, 20]{$\l_1\geq 1$}
  \recttext ( 50, 90)[ 50, 20]{$-\e_1$}
  \recttext ( 50, 70)[ 50, 20]{$+\e_2$}
  \rectangle(100, 50)[ 50, 40]{White}
  \recttext (100, 90)[ 50, 20]{$\l_2\geq 1$}
  \recttext (100, 70)[ 50, 20]{$-\e_2$}
  \recttext (100, 50)[ 50, 20]{$+\e_3$}
  \multiput (162, 60)(4,-1){14}{$\cdot$}
  \rectangle(230, 20)[ 50, 40]{White}
  \recttext (230, 60)[ 50, 20]{$\l_{r-2}\geq 1$}
  \recttext (230, 40)[ 50, 20]{$-\e_{r-2}$}
  \recttext (230, 20)[ 50, 20]{$+\e_{r-1}$}
  \rectangle(280,  0)[ 50, 40]{White}
  \recttext (283, 40)[ 50, 20]{$\l_{r}\geq 1$}
  \recttext (280, 20)[ 50, 20]{$-\e_{r-1}$}
  \recttext (280,  0)[ 50, 20]{$+\e_r$}
  \rectangle(282, 20)[ 98, 19]{White}
  \recttext (382, 20)[ 50, 20]{$\l_{r-1}\geq 1$}
  \recttext (330, 20)[ 50, 20]{$-\e_r$}
 \end{picture}\end{center}
 A weight $\e$ of the holonomy representation $T$ of $\so_{2r}$ is relevant
 for the irreducible representation $V_\l$ if and only if $\l$ satisfies
 all the conditions labeling the boxes containing $\e$. Say the weights
 $-\e_1$ and $+\e_2$ are relevant for all irreducible representations $V_\l$
 with $\l_1\geq1$, whereas $-\e_{r-1}$ is relevant for $V_\l$ if and only
 if $\l_{r-1}\geq1$ {\em and} $\l_r\geq1$.

 \pfill
 Odd dimensional Riemannian geometry $\g\,=\,\so_{2r+1},\,r\geq 1,$ with
 generic holonomy is of course closely related to $\g\,=\,\so_{2r}$. The
 weights $\pm\e_1,\,\pm\e_2,\,\ldots,\,\pm\e_r$ of the defining holonomy
 representation $T$ besides the zero weight form an orthonormal basis for
 a suitable scalar product $\<,>$ on the dual $\t^*$ of the maximal torus.
 With a suitable choice of ordering of weights the fundamental weights
 $\omega_1,\,\ldots,\,\omega_r$ and the weights of $T$ relate via:
 $$
  \begin{array}{lclclcl}
   \o_1       &=& \e_1 &\qquad&
   \pm\e_1    &=& \pm\,\o_1 \\
   \o_2       &=& \e_1+\e_2 &&
   \pm\e_2    &=& \pm(\o_2-\o_1) \\
   \;\;\vdots & & \quad\vdots && \;\;\vdots && \quad\vdots \\
   \o_{r-1}   &=& \e_1+\ldots+\e_{r-1} &&
   \pm\e_{r-2}&=& \pm(\o_{r-1}-\o_{r-2}) \\
   \o_r       &=& \frac12(\e_1+\ldots+\e_{r-1}+\e_r) &&
   \pm\e_r    &=& \pm(2\o_r-\o_{r-1})
  \end{array}
 $$
 Writing a dominant integral weight $\l\,=\,\l_1\o_1+\ldots+\l_r\o_r$ as a
 linear combination of fundamental weights with integers $\l_1,\,\ldots,\,
 \l_r\geq 0$ the criterion of Lemma \ref{relw} becomes:
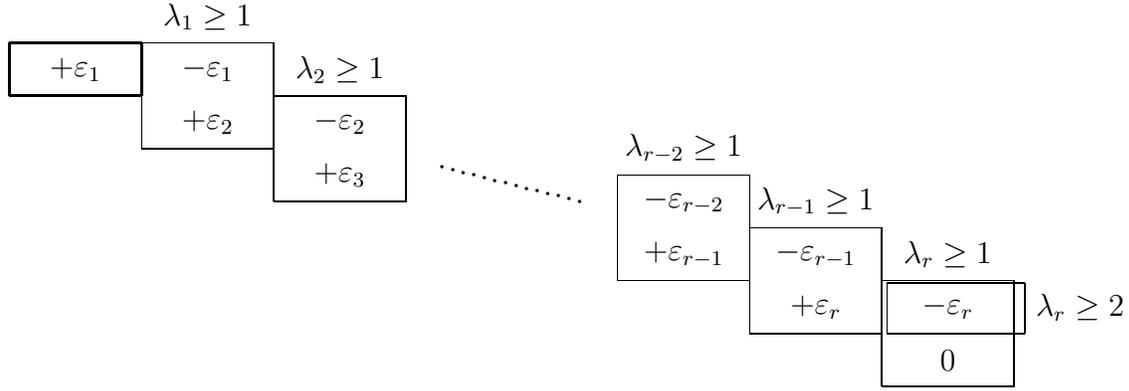
 \begin{center}\begin{picture}(430,150)(0,-20)
  \linethickness{.8pt}
  \rectangle(  0, 90)[ 50, 20]{Gray}
  \recttext (  0, 90)[ 50, 20]{$+\e_1$}
  \linethickness{.2pt}
  \rectangle( 50, 70)[ 50, 40]{White}
  \recttext ( 50,110)[ 50, 20]{$\l_1\geq 1$}
  \recttext ( 50, 90)[ 50, 20]{$-\e_1$}
  \recttext ( 50, 70)[ 50, 20]{$+\e_2$}
  \rectangle(100, 50)[ 50, 40]{White}
  \recttext (100, 90)[ 50, 20]{$\l_2\geq 1$}
  \recttext (100, 70)[ 50, 20]{$-\e_2$}
  \recttext (100, 50)[ 50, 20]{$+\e_3$}
  \multiput (162, 60)(4,-1){14}{$\cdot$}
  \rectangle(230, 20)[ 50, 40]{White}
  \recttext (230, 60)[ 50, 20]{$\l_{r-2}\geq 1$}
  \recttext (230, 40)[ 50, 20]{$-\e_{r-2}$}
  \recttext (230, 20)[ 50, 20]{$+\e_{r-1}$}
  \rectangle(280,  0)[ 50, 40]{White}
  \recttext (280, 40)[ 50, 20]{$\l_{r-1}\geq 1$}
  \recttext (280, 20)[ 50, 20]{$-\e_{r-1}$}
  \recttext (280,  0)[ 50, 20]{$+\e_r$}
  \rectangle(330,-20)[ 50, 40]{White}
  \recttext (330, 20)[ 50, 20]{$\l_r\geq 1$}
  \recttext (330, 00)[ 50, 20]{$-\e_r$}
  \recttext (330,-20)[ 50, 20]{$0$}
  \rectangle(332,  0)[ 52, 19]{White}
  \recttext (380,  0)[ 50, 20]{$\l_r\geq 2$}
 \end{picture}\end{center}

 \pfill
 Turning from the Riemannian case to the K\"ahler case $\g\,=\,\u_n$ we
 observe that the weights $\pm\e_1,\,\ldots,\,\pm\e_n$ of the defining
 standard representation $T\,=\,E\oplus\bar E$ form an orthonormal basis
 for an invariant scalar product on the dual $\t^*$ of a maximal torus
 $\t\subset\u_n$, but they become linearly dependent when projected to
 the dual of a maximal torus of the ideal $\su_n\subset\u_n$. In any
 case the fundamental weights and the weights of $T$ relate as
 $$
  \begin{array}{lclclcl}
   \o_1       &=& \e_1 &\qquad&
   \pm\e_1    &=& \pm\,\o_1 \\
   \o_2       &=& \e_1+\e_2 &&
   \pm\e_2    &=& \pm(\o_2-\o_1) \\
   \;\;\vdots & & \quad\vdots && \;\;\vdots && \quad\vdots \\
   \o_n       &=& \e_1+\ldots+\e_n &&
   \pm\e_n    &=& \pm(\o_n-\o_{n-1})
  \end{array}
 $$
 and the criterion of Lemma \ref{relw} becomes:
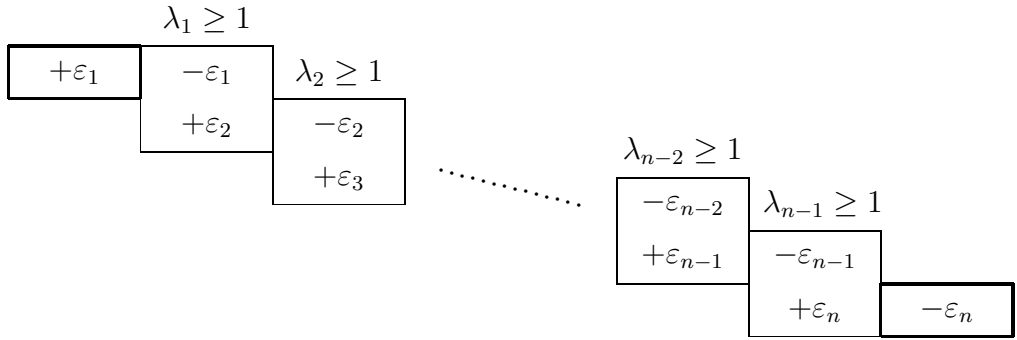
 \begin{center}\begin{picture}(380,130)(0,0)
  \linethickness{.8pt}
  \rectangle(  0, 90)[ 50, 20]{Gray}
  \recttext (  0, 90)[ 50, 20]{$+\e_1$}
  \linethickness{.2pt}
  \rectangle( 50, 70)[ 50, 40]{White}
  \recttext ( 50,110)[ 50, 20]{$\l_1\geq 1$}
  \recttext ( 50, 90)[ 50, 20]{$-\e_1$}
  \recttext ( 50, 70)[ 50, 20]{$+\e_2$}
  \rectangle(100, 50)[ 50, 40]{White}
  \recttext (100, 90)[ 50, 20]{$\l_2\geq 1$}
  \recttext (100, 70)[ 50, 20]{$-\e_2$}
  \recttext (100, 50)[ 50, 20]{$+\e_3$}
  \multiput (162, 60)(4,-1){14}{$\cdot$}
  \rectangle(230, 20)[ 50, 40]{White}
  \recttext (230, 60)[ 50, 20]{$\l_{n-2}\geq 1$}
  \recttext (230, 40)[ 50, 20]{$-\e_{n-2}$}
  \recttext (230, 20)[ 50, 20]{$+\e_{n-1}$}
  \rectangle(280,  0)[ 50, 40]{White}
  \recttext (283, 40)[ 50, 20]{$\l_{n-1}\geq 1$}
  \recttext (280, 20)[ 50, 20]{$-\e_{n-1}$}
  \recttext (280,  0)[ 50, 20]{$+\e_n$}
  \linethickness{.8pt}
  \rectangle(330,  0)[ 50, 20]{Gray}
  \recttext (330,  0)[ 50, 20]{$-\e_n$}
 \end{picture}\end{center}

 \pfill
 The quaternionic K\"ahler and hyperk\"ahler cases are more complicated,
 because the condition of being relevant has to be checked for both ideals
 $\sp(1)$ and $\sp(n)$ of $\g$. For a single ideal however the condition
 becomes simple again. The weights $\pm\e_1,\,\ldots,\,\pm\e_n$ of $E$ are
 again orthonormal for a suitable scalar product $\<,>$ on the dual $\t^*$
 of a maximal torus $\t\subset\sp(r)$ for $r=1$ or $r=n$ and relate to the
 fundamental weights by the formulas:
 $$
  \begin{array}{lclclcl}
   \o_1       &=& \e_1 &\qquad&
   \pm\e_1    &=& \pm\,\o_1 \\
   \o_2       &=& \e_1+\e_2 &&
   \pm\e_2    &=& \pm(\o_2-\o_1) \\
   \;\;\vdots & & \quad\vdots && \;\;\vdots && \quad\vdots \\
   \o_r       &=& \e_1+\ldots+\e_r &&
   \pm\e_r    &=& \pm(\o_r-\o_{r-1})
  \end{array}
 $$
 The graphical interpretation of Lemma \ref{relw} is given by the diagram:
 
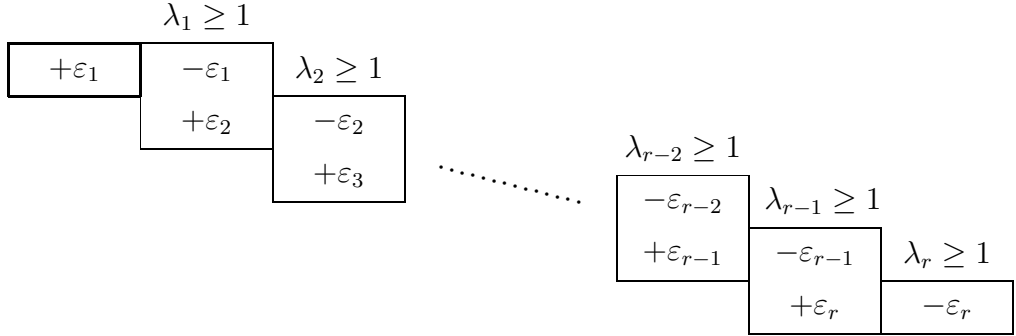
\begin{center}\begin{picture}(380,130)(0,0)
  \linethickness{.8pt}
  \rectangle(  0, 90)[ 50, 20]{Gray}
  \recttext (  0, 90)[ 50, 20]{$+\e_1$}
  \linethickness{.2pt}
  \rectangle( 50, 70)[ 50, 40]{White}
  \recttext ( 50,110)[ 50, 20]{$\l_1\geq 1$}
  \recttext ( 50, 90)[ 50, 20]{$-\e_1$}
  \recttext ( 50, 70)[ 50, 20]{$+\e_2$}
  \rectangle(100, 50)[ 50, 40]{White}
  \recttext (100, 90)[ 50, 20]{$\l_2\geq 1$}
  \recttext (100, 70)[ 50, 20]{$-\e_2$}
  \recttext (100, 50)[ 50, 20]{$+\e_3$}
  \multiput (162, 60)(4,-1){14}{$\cdot$}
  \rectangle(230, 20)[ 50, 40]{White}
  \recttext (230, 60)[ 50, 20]{$\l_{r-2}\geq 1$}
  \recttext (230, 40)[ 50, 20]{$-\e_{r-2}$}
  \recttext (230, 20)[ 50, 20]{$+\e_{r-1}$}
  \rectangle(280,  0)[ 50, 40]{White}
  \recttext (283, 40)[ 50, 20]{$\l_{r-1}\geq 1$}
  \recttext (280, 20)[ 50, 20]{$-\e_{r-1}$}
  \recttext (280,  0)[ 50, 20]{$+\e_r$}
  \rectangle(330,  0)[ 50, 20]{White}
  \recttext (330, 20)[ 50, 20]{$\l_r\geq 1$}
  \recttext (330,  0)[ 50, 20]{$-\e_r$}
 \end{picture}\end{center}

 \pfill
 Finally we consider the two exceptional cases $\g_2$ and $\spin_7$.
 Recall that the group $\G_2$ is the group of automorphisms of the octonions
 $\mathbb{O}$ as an algebra over $\R$. In this sense the holonomy
 representation $T_\R$ is the defining representation $\mathrm{Im}\,
 \mathbb{O}$ of $\g_2$ with complexification $T=[7]$. There are too
 many weights of the holonomy representation to be orthonormal for
 any scalar product on the dual $\t^*$ of a fixed maximal torus
 $\t\subset\g_2$, but at least we can choose an ordering of weights
 for $\t^*$ so that the weights of $T$ become totally ordered
 $+\e_1>+\e_2>+\e_3>0>-\e_3>-\e_2>-\e_1$. In this notation we
 have:
 $$
  \begin{array}{lclclcl}
   \o_1       &=& \e_1 &\qquad&
   \pm\e_1    &=& \pm\,\o_1 \\[1mm]
   \o_2       &=& \e_1+\e_2 &&
   \pm\e_2    &=& \pm(\o_2-\o_1) \\[1mm]
   && && \pm\e_3 &=& \mp (\o_2-2\o_1)
  \end{array}
 $$
 The scalar product of choice on $\t^*$ is specified by $\<\e_1,\e_1>=1=
 \<\e_2,\e_2>$ and $\<\e_1,\e_2>=\frac12$. Writing a dominant integral
 weight as $\l\,=\,a\o_1+b\o_2,\,a,\,b\,\geq\,0,$ we read Lemma \ref{relw} as:
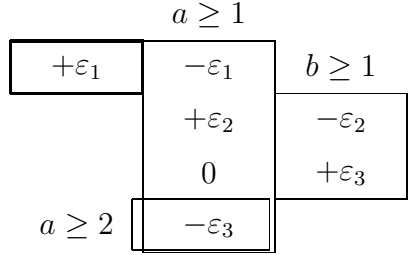
 \begin{center}\begin{picture}(150,100)(0,0)
  \linethickness{.8pt}
  \rectangle(  0, 60)[ 50, 20]{Gray}
  \recttext (  0, 60)[ 50, 20]{$+\e_1$}
  \linethickness{.2pt}
  \rectangle( 50,  0)[ 50, 80]{White}
  \recttext ( 50, 80)[ 50, 20]{$a\geq 1$}
  \recttext ( 50, 60)[ 50, 20]{$-\e_1$}
  \recttext ( 50, 40)[ 50, 20]{$+\e_2$}
  \recttext ( 50, 20)[ 50, 20]{$0$}
  \rectangle( 46,  1)[ 52, 19]{White}
  \recttext (  0,  0)[ 50, 20]{$a\geq 2$}
  \recttext ( 50,  0)[ 50, 20]{$-\e_3$}
  \rectangle(100, 20)[ 50, 40]{White}
  \recttext (100, 60)[ 50, 20]{$b\geq 1$}
  \recttext (100, 40)[ 50, 20]{$-\e_2$}
  \recttext (100, 20)[ 50, 20]{$+\e_3$}
 \end{picture}\end{center}

 \pfill
 The holonomy representation of the holonomy algebra $\g\,=\,\spin_7$ is
 the $8$--dimensional spinor representation $T=[8]$. It is convenient
 to write the weights $\pm\e_1,\,\ldots,\,\pm\e_4$ of $T$ and the fundamental
 weights $\o_1,\,\o_2$ and $\o_3$ in terms of the weights $\pm\eta_1,\,
 \pm\eta_2,\,\pm\eta_3,\,0$ of the $7$--dimensional defining representation
 of $\spin_7$, which form an orthonormal basis for a suitable scalar product
 on the dual $\t^*$ of the maximal torus. With this proviso the weights
 $\pm\e_1,\,\ldots,\,\pm\e_4$ of $T$ and the fundamental weights $\o_1,
 \,\o_2$ and $\o_3$ can be written as:
 $$
  \begin{array}{lclclclcl}
   \o_1    &=& \eta_1 &\qquad&
   \pm\e_1 &=& \pm\frac12(\eta_1+\eta_2+\eta_3)
           &=& \pm\,\o_3 \\[1mm]
   \o_2    &=& \eta_1+\eta_2 &&
   \pm\e_2 &=& \pm\frac12(\eta_1+\eta_2-\eta_3)
           &=& \pm(\o_2-\o_3) \\[1mm]
   \o_3    &=& \frac12(\eta_1+\eta_2+\eta_3) &&
   \pm\e_3 &=& \pm\frac12(\eta_1-\eta_2+\eta_3)
           &=& \pm(\o_3-\o_2+\o_1) \\[1mm]
   && &&
   \pm\e_4 &=& \pm\frac12(\eta_1-\eta_2-\eta_3)
           &=& \pm(\o_1-\o_3)
  \end{array}
 $$
 and Lemma \ref{relw} for a dominant integral weight $\l\,=\,a\o_1+b\o_2+c
 \o_3$ translates into:
 \begin{center}\begin{picture}(200,100)(-50,0)
  \linethickness{.8pt}
  \rectangle(  0, 60)[ 50, 20]{Gray}
  \recttext (  0, 60)[ 50, 20]{$+\e_1$}
  \linethickness{.2pt}
  \rectangle( 50,  0)[ 50, 80]{White}
  \recttext ( 50, 80)[ 50, 20]{$c\geq 1$}
  \recttext ( 50, 60)[ 50, 20]{$-\e_1$}
  \recttext ( 50, 40)[ 50, 20]{$+\e_2$}
  \recttext ( 50, 20)[ 50, 20]{$-\e_3$}
  \recttext ( 50,  0)[ 50, 20]{$+\e_4$}
  \rectangle(100, 20)[ 50, 40]{White}
  \recttext (100, 60)[ 50, 20]{$b\geq 1$}
  \recttext (100, 40)[ 50, 20]{$-\e_2$}
  \recttext (100, 20)[ 50, 20]{$+\e_3$}
  \rectangle(  0, 20)[ 98, 20]{White}
  \recttext (-50, 20)[ 50, 20]{$a\geq 1$}
  \recttext (  0, 20)[ 50, 20]{$-\e_4$}
 \end{picture}\end{center}
\section{The Space $\W(V)$ of Weitzenb\"ock Formulas}
\label{wspace}
 In this section we will define twistor operators, Weitzenb\"ock
 formulas and the space of Weitzenb\"ock formula
 with its different realizations. Then we will introduce the conformal
 weight operator, which in many cases generates all possible
 Weitzenb\"ock formulas. Finally we define the classifying
 endomorphism and study the corresponding eigenspace decomposition.
\subsection{Weitzenb\"ock Formulas}
 We consider parallel second order differential operators
 $P$ on sections of a vector bundle $VM$ over a Riemannian manifold
 $M$ with special holonomy $G$. By definition these are differential
 operators, which up to first order differential operators can always
 be written as the composition
 $$
  \Gamma(VM)
  \;\stackrel{\nabla^2}\longrightarrow\;\Gamma(T^*M\otimes T^*M\otimes VM)
  \;\stackrel\cong\longrightarrow\;\Gamma(TM\otimes TM\otimes VM)
  \;\stackrel F\longrightarrow\;\Gamma(VM)
 $$
 where $F$ is a parallel section of the vector bundle $\Hom(TM\otimes TM
 \otimes VM,VM)$ corresponding to a $G$--equivariant homomorphism $F\,\in\,
 \Hom_G(T\otimes T\otimes V,V)$. A particularly important example is
 the connection Laplacian $\nabla^*\nabla$ which arises from the linear
 map $a\otimes b\otimes\psi\longmapsto-\<a,b>\psi$. Note that we are
 only considering reduced holonomy groups $G$, which are connected by
 definition, so that $G$--equivariance is equivalent to $\g$--equivariance.
 Taking advantage of this fact we describe other parallel differential
 operators by means of the following identifications of spaces of invariant
 homomorphisms:
 $$
   \Hom_\g(\,T\otimes T\otimes V,\,V\,)
   \;\;=\;\;
   \Hom_\g(\,T\otimes T,\,\End\,V\,)
   \;\;=\;\;
   \End_\g(\,T\otimes V\,)
 $$
 Of course the identification $\Hom_\g(T\otimes T\otimes V,V)\,=\,\Hom_\g
 (T\otimes T,\End\,V)$ is the usual tensor shuffling $F(a\otimes b\otimes v)
 \,=\,F_{a\otimes b}v$ for all $a,\,b\,\in\,T$ and $v\,\in\,V$. The second
 identification $\Hom_\g(T\otimes T\otimes V,V)\,=\,\End_\g(T\otimes V)$
 depends on the existence of a $\g$--invariant scalar product on $T$ or
 the musical isomorphism $T\,\cong\,T^*$ via a summation
 $$
  F(\,b\otimes v\,)\;\;=\;\;\sum_\mu t_\mu\otimes F(t_\mu\otimes b\otimes v)
  \qquad
  F(a\otimes b\otimes v)\;\;=\;\;(\<a,\cdot>\ins\otimes\id)\,F(b\otimes v)
 $$
 over an orthonormal basis $\{t_\mu\}$. Under this identification the identity
 of $T \otimes V$ becomes the homomorphism $a\otimes b\otimes\psi\longmapsto
 \<a,b>\psi$ corresponding to the connection Laplacian $-\nabla^*\nabla$.
 The composition of endomorphisms turns $\End_\g(\,T\otimes V\,)$ and thus
 $\Hom_\g(T\otimes T,\,\End\,V)$ into an algebra, for $F,\,\tilde F\,\in\,
 \Hom_\g(T\otimes T,\End\,V)$ the resulting algebra structure reads:
 \begin{equation}\label{basi}
  (\,F\circ \tilde F\,)_{a\otimes b}
  \;\;=\;\;
  \sum_\mu \, F_{a\otimes t_\mu}\,\circ\,\tilde F_{t_\mu\otimes b }
 \end{equation}
 Last but not least we note that the invariance condition for $F\,\in\,
 \Hom_\g(\,T\otimes T,\End\,V\,)$ is equivalent to the identity $[\,X,\,
 F_{a\otimes b}\,]\,=\,F_{Xa\otimes b}+F_{a\otimes Xb}$ for all $X\in\g$ and
 $a,\,b\in T$.

 \pfill
 Assuming that $V=V_\l$ is irreducible of highest weight $\l$ we know from
 Lemma \ref{relw} that the isotypical components of $T\otimes V_\l$ are
 irreducible for non--symmetric holonomy groups. The algebra $\End_\g(\,
 T\otimes V_\l\,)$ is thus commutative and spanned by the pairwise orthogonal
 idempotents $\pr_\e$ projecting onto the irreducible subspaces $V_{\l+\e}$
 of $T\otimes V_\l$. In order to describe the corresponding second order
 differential operators we introduce first order differential operators
 $T_\e$ known as Stein--Weiss operators or generalized gradients by:
 $$
  T_\e:\quad\Gamma(\,V_\l M\,)\;\longrightarrow\;\Gamma(\,V_{\l + \e} M\,),
  \qquad
  \psi\;\longmapsto\;\pr_\e(\nabla \psi) \ .
 $$
 A typical example of a Stein--Weiss operator is the twistor operator of
 spin geometry, which projects the covariant derivative of a spinor onto
 the kernel of the Clifford multiplication. Straightforward calculations
 show that the second order differential operator associated to the
 idempotent $\pr_\e$ is the composition of $T_\e$ with its formal abjoint
 operator $T^*_\e:\,\Gamma(V_{\l + \e}M)\longrightarrow\Gamma(V_\l M)$ in
 the sense $\pr_\e(\nabla^2)\,=\,-T^*_\e T^{\phantom{*}}_\e$ compare
 \cite{uwe}. In consequence we can write the second order differential
 operator $F(\nabla^2)$ associated to $F\in\Hom_\g(\,T\otimes T,V_\l\,)$
 as a linear combination of the squares of Stein--Weiss operators:
 \begin{equation}\label{nabla}
  F(\;\nabla^2\;)\;\;=\;\;-\;\sum_\e f_\e\,T^*_\e T^{\phantom{*}}_\e\, .
 \end{equation}
 In fact with $\End_\g(\,T\otimes V_\l\,)\,=\,\Hom_\g(\,T\otimes T,\,
 \End\,V_\l\,)$ being spanned by the idempotents $\pr_\e$ every $F\in\End_\g
 (\,T\otimes V_\l\,)$ expands as $F\,=\,\sum_\e f_\e\,\pr_\e$ with coefficients
 $f_\e$ determined by $F|_{V_{\l+\e}}\,=\,f_\e\,\id$. A particular instance of
 equation (\ref{nabla}) is the identity $\nabla^*\nabla\,=\,\sum_\e T^*_\e
 T^{\phantom{*}}_\e$ associated to the expansion $\id_{T\otimes V_\l}\,=\,
 \sum_\e\pr_\e$. Motivated by this and other well--known identities of second
 order differential operators of the form (\ref{nabla}) we will in general
 call all elements $F\,\in\,\Hom_\g(\,T\otimes T\otimes V,\,V\,)\,=\,
 \End_\g(T\otimes V)$ Weitzenb\"ock formulas:

 \begin{Definition}[Space of Weitzenb\"ock Formulas on $VM$]
 \label{stiso}\hfill\break
  The Weitzenb\"ock formulas on a vector bundle $VM$ correspond bijectively
  to vectors in:
  $$
   \W(\,V\,)
   \;\;:=\;\;
   \Hom_\g(\,T\otimes T\otimes V,\,V\,)
   \;\;=\;\;
   \Hom_\g(\,T\otimes T,\,\End\,V\,)
   \;\;=\;\;
   \End_\g(\,T\otimes V\,)\,.
  $$
 \end{Definition}

 \noindent
 Of course we are mainly interested in Weitzenb\"ock formulas inducing
 differential operators of zeroth order or equivalently ``pure curvature
 terms''. Clearly a Weitzenb\"ock formula $F\in\Hom_\g(\,T\otimes T\otimes V,
 \,V\,)$ skew--symmetric in its two $T$--arguments will induce a pure curvature
 term $F(\nabla^2)$, because we can then reshuffle the summation in the
 calculation:
 \begin{equation}\label{cterm}
  F(\,\nabla^2 v\,)
  \;\;=\;\;
  \tfrac12\,\sum_{\mu\nu}F(\,t_\mu\otimes t_\nu\otimes(\nabla^2_{t_\mu,t_\nu}
  -\nabla^2_{t_\nu,t_\mu})v\,)
  \;\;=\;\;
  \tfrac12\,\sum_{\mu\nu}F_{t_\mu\otimes t_\nu}\,R_{t_\mu,\,t_\nu}v\ .
 \end{equation}
 Here and in the following we will denote with $\{t_\nu\}$ an orthonormal
 basis of $T$ and also a local orthonormal basis of the tangent bundle.
 Conversely the principal symbol of the differential operator $F(\nabla^2)$
 is easily computed to be $\sigma_{F(\nabla^2)}(\xi)v\,=\,F_{\xi^\flat\otimes
 \xi^\flat}v$ for every cotangent vector $\xi$ and every $v\in VM$. Hence
 the principal symbol vanishes identically exactly for the skew--symmetric
 Weitzenb\"ock formulas. Weitzenb\"ock formulas $F$ leading to a pure curvature
 term $F(\nabla^2)$ are thus completely characterized by being eigenvectors
 of eigenvalue $-1$ for the involution
 $$
  \tau:\qquad
  \W(\,V\,)\;\longrightarrow\;
  \W(\,V\,),\qquad F\;\longmapsto\;\tau(\,F\,)
 $$
 defined in the interpretation $\W(V)\,=\,\Hom_\g(\,T\otimes T\otimes V,\,
 V\,)$ as precomposition with the twist $\tau:\,T\otimes T\otimes V
 \longrightarrow T\otimes T\otimes V,\,a\otimes b\otimes v\longmapsto
 b\otimes a\otimes v$. In other words a Weitzenb\"ock formula $F$ will
 reduce to a pure curvature term if and only if $\tau(F)\,:=\,F\circ\tau
 \,=\,-F$.

 \pfill
 Considering the space of Weitzenb\"ock formulas $\W(V)$ as the algebra
 $\End_\g(T\otimes V)$ we can introduce additional structures on it: the
 unit $\1\,:=\,\id_{T\otimes V}\in\W(V)$, the scalar product
 $$
  \<F,\tilde F>
  \;\;:=\;\;
  \frac1{\dim\,V}\,\tr_{T\otimes V}(\,F\,\tilde F\,)
  \qquad\qquad F,\,\tilde F\;\in\;\W(V)
 $$
 satisfying $\,\<FG,\tilde F>\,=\,\<F,G\tilde F>\,$ and the trace
 $\,\tr\,F\,:=\,\<F,\1>$. Clearly the trace of the unit is given
 by $\tr\,\1\,=\,\dim\,T$. The definition of the trace can be
 rewritten in the form
 $$
  \tr\,F
  \;\;=\;\;
  \frac 1{\dim\,V}\,
  \tr_V\Big(\;v\;\longmapsto\;\sum_\mu F_{t_\mu\otimes t_\mu}v\;\Big)
 $$
 so that the trace is invariant under the twist $\tau$. A slightly more
 complicated argument using (\ref{basi}) shows that the scalar product
 is invariant under the twist, too. In particular the eigenspaces for $\tau$
 for the eigenvalues $\pm 1$ are orthogonal and all eigenvectors in the
 $(-1)$--eigenspace of $\tau$ have vanishing trace. Fom the definition
 of the trace we obtain that the trace of an element $F\,=\,\sum f_\e \pr_\e$
 of $\W(V_\l)$ in the irreducible case is given by
 \begin{equation}\label{trace}
  \tr\,F
  \;\;=\;\;
  \sum_\e f_\e \;\frac{\dim\,V_{\l+\e}}{\dim\,V_\l}
 \end{equation}
 in particular the idempotents $\pr_\e$ form an orthogonal basis of $\W(V_\l)$:
 $$
  \<\pr_\e,\pr_{\tilde\e}>\;\;=\;\;\delta_{\e\tilde\e}\;
  \frac{\dim V_{\l+\e}}{\dim V_\l}
 $$
 A different way to interprete the trace is to note that for every
 Weitzenb\"ock formula $F\in\W(V)$ considered as an equivariant homomorphism
 $F:\,T\otimes T\longrightarrow\End\,V$ the trace endomorphism $\sum_\mu
 F_{t_\mu\otimes t_\mu}\in\End_\g V$ is invariant under the action of
 $\g$. For an irreducible representation $V_\l$ it is thus the multiple
 $\sum_\mu F_{t_\mu\otimes t_\mu}=(\tr\,F)\,\id_{V_\l}$ of $\id_{V_\l}$
 by Schur's Lemma.
\subsection{The conformal weight operator}
 In order to study the fine structure of the algebra $\W(V)=\End_\g(
 T\otimes V)$ of Weitzenb\"ock formulas it is convenient to introduce
 the {\it conformal weight operator} $B\in\W(V)$ of the holonomy algebra
 $\g$ and its variations $B^\h\in\W(V)$ associated to the non--trivial
 ideals $\h\subset\g$ of $\g$. All these conformal weight operators
 commute and the commutative subalgebra of $\W(V)$ generated by them
 in the irreducible case $V=V_\l$ is actually all of $\W(V)$ for
 generic highest weight $\l$. In this subsection we work out some
 direct consequences of the description of Weitzenb\"ock formulas
 as polynomials in the conformal weight operators.

 \pfill
 Recall that the scalar product $\<,>$ on $T$ induces a scalar product
 on all exterior powers $\L^kT$ of $T$ via Gram's determinant. Using this
 scalar product on $\L^2T$ we can identify the adjoint representation
 $\so\,T$ of $\SO(n)$ with $\L^2T$ through $\<X,a\wedge b>\,=\,\<X a,b>$
 and hence think of the holonomy algebra $\g\,\subset\,\so\,T$ as a subspace
 of the euclidian vector space $\L^2T$:

 \begin{Definition}[Conformal Weight Operator]
 \label{defw}\hfill\break
  Consider an ideal $\h_\R\subset\g_\R$ in the real holonomy algebra.
  Its complexification $\h\,:=\,\h_\R\otimes_\R\C$ is an ideal in $\g$
  and a regular subspace $\h\subset\g\subset\L^2T$ in $\L^2T$ with associated
  orthogonal projection $\pr_\h:\,\L^2T\longrightarrow\h$. The conformal
  weight operator $B^\h\,\in\,\W(V)$ is defined by
  $$
   B^\h_{a \otimes b}\,v\;\;:=\;\;\pr_\h (a \wedge b)\;v
  $$
  in the interpretation of Weitzenb\"ock formulas as linear maps
  $T\otimes T\longrightarrow\End\,V$. Under the identification
  $\Hom_\g(T\otimes T,\End\,V)\,=\,\End_\g(T\otimes V)$ discussed
  above the conformal weight operator $B^\h$ becomes the following
  sum over an orthonormal basis $\{t_\mu\}$ of $T$:
  $$
   B^\h(\,b\otimes v\,)
   \;\;=\;\;
   \sum_\mu t_\mu\otimes\pr_\h(t_\mu\wedge b)\,v
  $$
  The notation $B\,:=\,B^\g$ will be used for the conformal weight
  operator of the algebra $\g$.
 \end{Definition}

 \noindent
 Most of the irreducible non--symmetric holonomy algebras $\,\g\,$ are simple
 and hence there is only one weight operator $\,B\,$ ( c.f. table (\ref{ht})).
 The exceptions are K\"ahler geometry $\,\g_\R=i\R\oplus\su_n\,$ with a
 one--dimensional center in dimension $\,2n\,$ and two commuting weight
 operators $\,B^{i\R}\,$ and $\,B^{\su}\,$ and quaternionic K\"ahler geometry
 $\,\g_\R=\sp(1)\oplus\sp(n)\,$ in dimension $\,4n,\,n\geq 2\,$ with two
 commuting weight operators $\,B^H\,$ and $\,B^E$.

 \begin{Lemma}[Fegan's Lemma \cite{fegan}]
 \label{fegan}\hfill\break
  The conformal weight operator $B^\h\in\W(V)$ of an ideal $\h\subset\g
  \subset\L^2T$ can be written
  $$
   B^\h\;\;=\;\;-\,\sum_\alpha\,X_\alpha\otimes X_\alpha
   \;\;\in\;\;\End_\g(\,T\otimes V\,)
  $$
  where $\{X_\alpha\}$ is an orthonormal basis of $\h$ for the scalar product
  on $\L^2T$ induced from $T$.
 \end{Lemma}

 \proof
 Let $\{t_\mu\}$ and $\{X_\alpha\}$ be orthonormal bases of $T$ and $\h$
 respectively. Using the characterization $\<X,a\wedge b>\,=\,\<X a,b>$
 of the identification $\so\,T\,=\,\L^2T$ we find:
 \begin{eqnarray*}
  B^\h(\,b\otimes v\,)
  &=& \sum_\mu t_\mu\otimes\pr_\h(\,t_\mu\wedge b\,)\,v\\
  &=& \sum_{\mu\alpha} t_\mu\otimes\<X_\alpha,t_\mu\wedge b>X_\alpha\,v
  \;\;=\;\;-\,\sum_\alpha X_\alpha b\otimes X_\alpha\,v\,.
  \hbox to0pt{\hskip17mm$\Box$\hss}
 \end{eqnarray*}

 \noindent
 A particularly nice consequence of Fegan's Lemma is that the conformal weight
 operators $B^\h$ and $B^{\tilde\h}$ associated to two ideals $\h,\,\tilde\h
 \subset\g$ always commute. In fact two disjoint ideals $\h\,\cap\,\tilde\h
 \,=\,\{0\}$ of $\g$ commute by definition, the general case follows easily.
 Hence the algebra structure on $\W(V)$ allows us to use the evaluation
 homomorphism
 \begin{equation}\label{evhom}
  \Phi:\qquad\C[\;\{\;B^\h\,|\;\h\textrm{\ irreducible ideal of\ }\g\;\}\;]
  \;\longrightarrow\;\W(\;V\;)
 \end{equation}
 from the polynomial algebra on abstract symbols $\{B^\h\}$ to the algebra
 $\End_\g(T\otimes V)$ for studying the fine structure of the space $\W(V)$
 of Weitzenb\"ock formulas.

 \pfill
 In order to turn Fegan's Lemma into an effective formula for the eigenvalues
 of the conformal weight operator $B^\h$ of an ideal $\h\subset\g$ we need to
 calculate the Casimir operator in the normalization given by the scalar
 product on $\L^2T$. Recall that the Casimir operator is defined as a sum
 $\Cas\,:=\,\sum_\alpha X^2_\alpha\,\in\,\U\h$ over an orthonormal basis
 $\{X_\alpha\}$ of $\h$ and is thus determined only up to a constant. Usually
 it is much more convenient to calculate the Casimir $\Cas$ with respect
 to a scalar product of choice and later normalize it to the Casimir
 $\Cas^{\L^2}$ with respect to the invariant scalar product induced
 on $\h\,\subset\,\L^2T$.

 For a given irreducible ideal $\h\subset\g$ in an irreducible holonomy
 algebra $\g$ the Casimir operator $\Cas$ for $\h$ is now real, symmetric
 and $\g$--invariant. Although the holonomy representation $T$ of $\g$
 may not be irreducible itself, it is the complexification of the
 irreducible real representation $T_\R$ so that we can still conclude that
 $\Cas$ acts as the scalar multiple $\Cas_T\,\id$ of the identity on $T$.
 The Casimir eigenvalue $\Cas^{\L^2}_{V_\l}$ of the properly normalized
 Casimir operator on a general irreducible representation $V_\l$ of $\g$
 of highest weight $\l$ can then be calculated from the Casimir $\Cas$ using
 \begin{equation}\label{cnor}
  \Cas^{\L^2}_{V_\l}
  \;\;=\;\;
  -\,2\,\frac{\dim \h}{\dim T}\;\frac{\Cas_{V_\l}}{\Cas_T}\,,
 \end{equation}
 where the ambiguity in the choice of normalization cancels out in
 the quotient $\frac{\Cas_{V_\l}}{\Cas_T}$. In fact the normalization
 (\ref{cnor}) is readily checked for the holonomy representation $V=T$
 $$
  \tr_T \,\Cas^{\L^2}
  \;\;=\;\;
  \dim T \cdot \Cas^{\L^2}_T
  \;\;=\;\;
  \sum_\alpha\tr_T\,X^2_\alpha\;\;=\;\; -\,2\,\dim \h\,,
 $$
 because the scalar product induced from $T$ on $\L^2T$ satisfies $\<X,Y>
 \,=\,-\frac12\tr_T\,XY$.

 \begin{Corollary}[Explicit Formula for Conformal Weights]
 \label{confW}\hfill\break
  Consider the tensor product $\,T \otimes V_\l\,=\,\oplus_{\e\subset\l}
  V_{\l+\e}\,$ of the holonomy representation $T$ with the irreducible
  representation $V_\l$ of highest weight $\l$. For an ideal $\h\,\subset
  \,\g$ let $\e_{max}$ be the highest weight of $T$ and $\rho$ be the
  half sum of positive weights of $\h$ in the dual $\t^*$ of a maximal
  torus $\t$. With respect to the basis $\{\pr_\e\}$ of idempotents the
  conformal weight operator $B^\h$ of the ideal $\h$ can be expanded $B^\h=
  \sum_{\e\subset\l}\,b^\h_\e\,\pr_\e$ with conformal weights
  $$
   b^\h_\e
   \;\;=\;\;
   2\,\frac{\dim \h}{\dim\,T}\,
   \frac{\<\l+\rho,\e>\;-\;\<\rho,\e_{max}>\;+\;
   \frac12\,(|\e|^2 - |\e_{max}|^2)}{\<\e_{max}+2\rho,\e_{max}>}\,,
  $$
  where $\<\cdot,\cdot>$ is an arbitrary scalar product on $\t^*$
  invariant under the Weyl group of $\h$.
 \end{Corollary}

 \proof
 According to Lemma~\ref{fegan} the conformal weight operator can be written
 as a difference $B^\h\,=\,-\frac12\,(\Cas^{\L^2} - \Cas^{\L^2}\otimes \id -
 \id\otimes\Cas^{\L^2})$ of properly normalized Casimir operators. In
 particular its restriction to the irreducible summand $V_{\l + \e}
 \subset T\otimes V_\l$ acts by multiplication with $b^\h_\e\,:=\,-\frac12
 (\Cas^{\L^2}_{V_{\l + \e}}-\Cas^{\L^2}_T - \Cas^{\L^2}_{V_\l})$. The
 conformal weights $b^\h_\e$ can thus be calculated using Freudenthal's
 formula $\,\Cas_{V_\l}\,=\,\<\l+2\rho,\l>\,$ for the Casimir eigenvalues
 of irreducible representations $V_\l$ and the normalization (\ref{cnor}).
 \qed

 \pfill
 It is clear from the definition that the conformal weight operator $B^\h\in
 \W(V)$ of an ideal $\h\subset\g$ of the holonomy algebra $\g$ is in the
 $(-1)$--eigenspace of the involution $\tau$ and thus induces a pure curvature
 term $B^\h(\,\nabla^2\,)$ on every vector bundle $VM$ associated to the
 holonomy reduction of $M$. Explicitly we can describe this curvature term
 using an orthonormal basis $\{X_\alpha\}$ of the ideal $\h$ for the scalar
 product induced on $\h\subset\L^2T$. Namely the curvature operator $R:\,
 \L^2TM\longrightarrow\g M\subset\L^2TM,\,a\wedge b\longmapsto R_{a,b},$
 associated to the curvature tensor $R$ of $M$ allows us to write down a
 well--defined global section
 \begin{equation}\label{qendo}
  q^\h(R)
  \;\;:=\;\;
  \sum_\alpha X_\alpha \,R(X_\alpha)
  \;\;\in\;\;
  \Gamma(\;\U^{\leq 2}\g M\;)
 \end{equation}
 of the universal enveloping algebra bundle associated to the holonomy
 reduction. Fixing a representation $G\longrightarrow\Aut(V)$ of the
 holonomy group the section $q^\h(R)$ in turn induces an endomorphism
 on the vector bundle $VM$ associated to $V$ and the holonomy reduction.

 A particularly important example of Weitzenb\"ock formulas
 is the classical formula of Weitzenb\"ock for the Laplace operator
 $\Delta = d^* \,d\,+\,d\,d^*$ acting on differential forms,~i.e.
 \begin{equation}\label{classical}
  \Delta \;\;=\;\;\nabla^* \, \nabla \;+\;q(R)
 \end{equation}
 The curvature term in this formula is precisely the curvature endomorphism
 for the full holonomy algebra $\g$, in particular $q(R)\,=\,\Ric$
 on the bundle of 1--forms on $M$. We recall that the curvature
 term in the Weitzenb\"ock formula (\ref{classical}) is known to be the
 Casimir operator of the holonomy algebra $\g$ on a symmetric space
 $M\,=\,\tilde G/G$, more precisely for every ideal $\h\subset\g$ the
 curvature term $q^\h(R)$ acts as the Casimir operator of the ideal $\h$
 on every homogeneous vector bundle $VM$ over a symmetric space $M$.

 A minor subtlety in the definition of the curvature terms $q^\h(R)$ should
 not pass unnoticed, in difference to the conformal weight operators $B^\h$
 the curvature terms $q^\h(R)$ associated to two ideals $\h,\,\tilde\h\subset
 \g$ do not in general commute. The problem is that the curvature operator
 $R:\,\L^2TM\longrightarrow\g M$ does not necessarily map the parallel
 subbundle $\h M\subset\L^2TM$ associated to an ideal $\h\subset\g$ to
 itself due to the presence of ``mixed terms'' in the curvature tensor
 $R$. In other words the section $q^\h(R)$ is not in general a section
 of $\U^{\leq2}\h M\subset\U^{\leq2}\g M$, so it is of no use that disjoint
 ideals centralize each other.

 Interestingly the problem of mixed terms is absent for symmetric holonomies
 as well as for the holonomies $\so_n,\,\g_2,\,\spin_7$ and the quaternionic
 holonomies $\sp(n)$ and $\sp(1)\,\sp(n)$. Mixed terms may however spoil
 commutativity of the curvature endomorphisms $q^\h(R)$ on a K\"ahler manifold
 $M$, in fact the curvature term associated to the center $i\R\subset\u_n$
 reads
 $$
  q^{i\R}(\;R\;)\;\;=\;\;-\,\frac1n\,I\,(\,I\Ric\,)\;\;\in\;\;
  \Gamma(\,\U^{\leq2}\u M\,)\ ,
 $$
 according to equation (\ref{ccop}), where $I$ is the parallel complex
 structure and $I\Ric$ is the composition of $I$ with the symmetric Ricci
 endomorphism of $TM$ thought of as a section of the holonomy bundle $\u M$.
 In consequence $q^{i\R}(R)$ is a section of $\U^{\leq2}\,i\R M$ if and only
 if $M$ is K\"ahler--Einstein, otherwise it will not in general commute with
 $q^{\su}(R)$. The central curvature term features prominently in the
 Bochner identity for K\"ahler manifolds.

 \begin{Lemma}
  $
   \hskip4cm
   B^\h(\;\nabla^2\;)\;\;=\;\;q^\h(\,R\,)
  $
 \end{Lemma}

 \proof
 Expanding the second covariant derivative $\nabla^2\psi=\sum t_\mu\otimes
 t_\nu\otimes\nabla^2_{t_\mu,\,t_\nu}\psi$ of the section $\psi$ with an
 orthonormal basis $\{t_\mu\}$ of $T$ and using the same resummation as
 in the derivation of equation (\ref{cterm}) we find for an orthonormal basis
 $\{X_\alpha\}$ of the ideal $\h$:
 \begin{eqnarray*}
  B^\h(\;\nabla^2\psi\;)
  &=&
  \frac12\,\sum_{\mu\nu}\,\pr_\h(t_\mu\wedge t_\nu)\,R^V_{t_\mu,\,t_\nu}\,\psi
  \\
  &=&
  \frac12\,\sum_{\alpha\mu\nu}\,\<t_\mu\wedge t_\nu,X_\alpha>\,X_\alpha\,
  R^V_{t_\mu,\,t_\nu}\,\psi
  \;\;=\;\;
  q^\h(\;R\;)
  \hbox to0pt{\hskip23mm$\Box$\hss}
 \end{eqnarray*}

 \noindent
 On the other hand Corollary~\ref{confW} tells us how to write the conformal
 weight operator $B^\h$ in terms of the basis $\,\{\pr_\e\}\,$ of projections
 onto the irreducible summands $V_{\l +\e}\subset T \otimes V_\l$. Using the
 identification of $B^\h(\nabla^2)$ with the universal curvature terms
 $q^\h(\,R\,)$ proved above we obtain some prime examples of Weitzenb\"ock
 formulas:

 \begin{Proposition}[Universal Weitzenb\"ock Formula]
 \hfill\label{universal}\break
  Consider a Riemannian manifold $M$ of dimension $n$ with holonomy group
  $G\subset\SO(n)$ and the vector bundle $V_\l M$ over $M$ associated to
  the holonomy reduction of $M$ and the irreducible representation $V_\l$
  of $G$ of highest weight $\l$. In terms of the Stein--Weiss operators
  $\,T_\e:\,\Gamma(V_\l M)\longrightarrow\Gamma(V_{\l+\e} M)$ arising from
  the decomposition $T \otimes V_\l\,=\,\oplus_{\e \subset \l} V_{\l+\e}$
  the action of the curvature endomorphisms $q^\h(R)$ can be written
  $$
   q^\h(R)
   \;\;=\;\;
   -\,\sum_{\e \subset \l}\,b^\h_\e
   \,\, T^*_\e\,T^{\phantom{*}}_\e \ ,
  $$
  where the $b^\h_\e$ are the eigenvalues of the conformal weight operator
  $B^\h\,\in\,\End_\g(T \otimes V_\l)$.
 \end{Proposition}

 \noindent
 As a direct consequence of Proposition~\ref{universal} and the
 classical Weitzenb\"ock formula (\ref{classical}) for the Laplace
 operator $\D =d d^* + d^*d$ on the bundle of differential forms we obtain
 $$
  \Delta
  \;\;=\;\;
  \sum_{\e \subset \l}\,(1\,-\,b_{\e})\,T^*_\e\,T^{\phantom{*}}_\e \ .
 $$
 In the case of Riemannian holonomy $G=\SO(n)$ the universal Weitzenb\"ock
 formula of Proposition \ref{universal} was considered in \cite{pg1} for the
 first time. The definition of the conformal weight operator and its expression
 in terms of the Casimir is taken from the same article. The conformal weight
 operator $B$ has been used for other purposes as well, see \cite{cgh} for
 example. Similar results can be found in~\cite{homma}.

 \pfill
 Considering $B$ as an element of the algebra $\W(V) $ all powers of $B$
 are $\g$--invariant endomorphisms. In the interpretation $\W(V)=\Hom_\g
 (T \otimes T, \End V)$ these powers read:
 \begin{equation}\label{wk}
  B^k_{a \otimes b}
  \;\;=\;\;
  \sum_{\mu_1,\ldots,\mu_{k-1}}\,
  \pr_\g(a            \wedge t_{\mu_1}) \,
  \pr_\g(t_{\mu_1}    \wedge t_{\mu_2})\,\ldots\,
  \pr_\g(t_{\mu_{k-2}}\wedge t_{\mu_{k-1}})\,
  \pr_\g(t_{\mu_{k-1}}\wedge b)\ .
 \end{equation}
 Recall now that in the irreducible case the trace $\sum F_{t_\mu\otimes
 t_\mu}=(\tr\,F)\id_{V_\l}$ of an element $F\in\W(V_\l)$ is a multiple of
 the identity of $V_\l$. Evidently the traces of the powers $B^k$ of $B$
 correspond to the action of the elements
 \begin{equation}\label{highercasimir}
  \Cas^{[k]}
  \;\;:=\;\;
  \sum_{\mu_1,\ldots,\mu_{k-1}}\,
  \pr_\g(t_{\mu_0}    \wedge t_{\mu_1}) \,
  \pr_\g(t_{\mu_1}    \wedge t_{\mu_2})\,\ldots\,
  \pr_\g(t_{\mu_{k-2}}\wedge t_{\mu_{k-1}})\,
  \pr_\g(t_{\mu_{k-1}}\wedge t_{\mu_0})\ .
 \end{equation}
 of the universal enveloping algebra $\U\g$ on $V$. The elements
 $\Cas^{[k]},\,k\geq2,$ all belong to the center of the universal enveloping
 algebra $\U \g$ and are called {\it higher Casimirs} since $\Cas^{[2]}=
 -2 \Cas^{\L^2}$ (c.f.~\cite{cgh}). A straightforward calculation shows
 \begin{equation}\label{kth}
  \Cas^{[k]}\;\;=\;\; \tr\,(B^k)\,\id_{V_\l} \;\;=\;\;
  \left(\sum_\e\,b_\e^k\,\frac{\dim\,V_{\l+\e}}
  {\dim\,V_\l}\right) \id_{V_\l}\ ,
 \end{equation}
 for an irreducible representation $V=V_\l$, where we use the second
 equation in (\ref{trace}) for computing the trace of $B^k=\sum b_\e^k
 \pr_\e$ explicitly. As an example we consider the equation $ \Cas^{[3]}=
 -\frac12\Cas^{\L^2}_\g \Cas^{\L^2}$, which follows from the recursion
 formula of Corollary~\ref{recursion2} or by direct calculation. Indeed
 $B^2-\frac14 \Cas^{\L^2}_\g B$ is an eigenvector of the involution $\tau$
 for the eigenvalue $+1$. Thus it is orthogonal to the eigenvector $B$ for
 the eigenvalue $-1$ and so:
 $$
  0 \;\;=\;\; \< B^2-\frac14 \Cas^{\L^2}_\g B, B> \;\;=\;\;
  \tr(B^3)-\frac14 \Cas^{\L^2}_\g\tr(B^2) \;\;=\;\; \tr(B^3)+\frac12
  \Cas^{\L^2}_\g\Cas^{\L^2} \ .
 $$
 From a slightly more general point of view the evaluation at the
 conformal weight operator defines an algebra homomorphism $\,\Phi:\;
 \C[B]\longrightarrow\End_\g(T \otimes V)$, whose kernel is generated
 by the minimal polynomial of $B$ as an endomorphism on $T\otimes V$.
 With $B$ being diagonalizable its minimal polynomial is the product
 $\mathrm{min}(B)\,=\,\prod_{b\,\in\,\{b_\e\}}(B\,-\,b)$ over all {\em
 different} conformal weights. In consequence the injective algebra
 homomorphism
 $$
  \Phi:\qquad \C[B]/\langle\mathrm{min}(\,B\,)\rangle
  \;\longrightarrow\;\End_\g(\;T \otimes\;V) \ .
 $$
 is an isomorphism as soon as all conformal weights are pairwise different.
 Indeed the dimension of $\End_\g(T \otimes V)$ is the number $N(G,\l)$
 of relevant weights or the number of conformal weights counted with
 multiplicity, while the number of different conformal weights determines
 the degree of $\mathrm{min}(B)$ and so the dimension of $\C[B]/\langle
 \mathrm{min}(B)\rangle$.

 In section~\ref{eigenvalues} we compute the $B$--eigenvalues in the cases
 $\g = \so_n, \g_2$ and $\spin_7$. It follows that they are pairwise different
 with the only exceptions of $\g= \so_{2r}$ and a representation of
 highest weight $\lambda = \lambda_1\omega_1 + \ldots +\lambda_r
 \omega_r$, with $\lambda_r = \lambda_{r-1}$, which is equivalent to
 $b_{\e_r} = b_{-\e_r}$ and $\g = \spin_7$ and a representation with
 highest weight $\,\l=a\o_1+b\o_2+c\o_3\,$ with $\,c=2a+1$, which is
 equivalent to $b_{-\e_4}= b_{\e_4}$. In these cases the degree of the
 minimal polynom is reduced by one and hence the image of $\Phi$ has
 codimension one. We thus have proved the following

 \begin{Proposition}[Structure of the Algebra of Weitzenb\"ock Formulas]
 \hfill\label{quotient}\break
  Let $G$ be one of the holonomy groups $\SO_{n},\,\G_2\,$ or $\Spin(7)$
  of non--symmetric manifolds. If $V_\l$ is irreducible, then $\Phi$ is
  an isomorphism
  $$
   \Phi:\;\;\C[B]/\langle\mathrm{min}(\,B\,)\rangle
   \;\stackrel\cong\longrightarrow\;\End_\g(\;T \otimes V_\l\;) \ ,
  $$
  with the only exception of the cases $\,G= \SO_{2r}\,$ and a highest
  weight $\,\l\,$ with $\,\l_{r-1} = \l_r$, or $\,G = \Spin(7)\,$ and a
  highest weight $\,\l=a\o_1+b\o_2+c\o_3\,$ with $\,c=2a+1$. In both
  cases the homomorphism $\Phi$ is not surjective and its image has
  codimension one.
 \end{Proposition}
\subsection{The Classifying Endomorphism}
 The decomposition of the space $\W(V)=\Hom_\g(T\otimes T,\End\,V)$ of
 Weitzenb\"ock formulas into the $(\pm 1)$--eigenspaces of the involution
 $\tau$ can be written as
 $$
  \Hom_\g(T\otimes T,\End\,V) \;\cong\;
  \Hom_\g(\L^2T,\End\,V)\,\oplus\,\Hom_\g(\S^2T,\End\,V) \ .
 $$
 However in general we have a further splitting of $T\otimes T$ leading to a
 further decomposition of the $\tau$--eigenspaces. Our aim is now to introduce
 an endomorphism $K$ on $\W(V)$ whose eigenspaces correspond to this finer
 decomposition.

 \begin{Definition}[The Classifying Endomorphism]
 \label{cls}\hfill\break
  The classifying endomorphism $K^\h$ of an ideal $\h_\R\subset\g_\R$ of
  the real holonomy algebra $\g_\R$ is the endomorphism $K^\h:\,\W(V)
  \longrightarrow\W(V)$ on the space of Weitzenb\"ock formulas defined
  in the interpretation $\W(V)\,=\,\Hom_\g(T\otimes T,\End\,V)$ by the
  formula
  $$
   K^\h(\,F\,)_{a\otimes b}\,v
   \;\;:=\;\;
   -\,\sum_\alpha\,F_{X_\alpha a\otimes X_\alpha b}v
  $$
  where $\{X_\alpha\}$ is an orthonormal basis for the scalar product induced
  on the ideal $\h\,\subset\,\L^2T$. As before we denote the classifying
  endomorphism of the ideal $\g$ simply by $K\,:=\,K^\g$.
 \end{Definition}

 \pfill
 Of course the definition of the classifying endomorphism $K^\h$ is motivated
 by Fegan's Lemma \ref{fegan} for the conformal weight operator $B^\h$.
 Note that for every $\g$--equivariant map $F:\;T\otimes T\longrightarrow
 \End\,V$ the map $K^\h(F):\;T\otimes T\longrightarrow\End\,V$ is again
 $\g$--equivariant, because we sum over an orthonormal basis $\{X_\alpha\}$
 of the ideal $\h\subset\g$ for a $\g$--invariant scalar product. In the
 interpretation $\W(V)\,=\,\End_\g(T\otimes V)$ the definition of $K^\h$ reads
 $$
  K^\h(\,F\,)(b\otimes v)
  \;\;=\;\;
  -\,\sum_{\mu\alpha}\,t_\mu\otimes F_{X_\alpha t_\mu\otimes X_\alpha b}\,v
  \;\;=\;\;
  \sum_{\mu\alpha}\,X_\alpha t_\mu\otimes F_{t_\mu\otimes X_\alpha b}\,v
 $$
 or more succinctly:
 \begin{equation}\label{ksym}
  K^\h(\,F\,)
  \;\;=\;\;
  \sum_\alpha\,(\,X_\alpha\otimes\id\,)\,F\,(\,X_\alpha\otimes\id\,)\ .
 \end{equation}
 In consequence the classifying endomorphisms $K^\h$ and $K^{\tilde\h}$
 for two ideals $\h,\,\tilde\h\subset\g$ commute on the space $\W(V)$
 of Weitzenb\"ock formulas similar to the conformal weights operators.
 The classifying endomorphisms will be extremely useful in finding the
 matrix of the twist $\tau:\;\W(V)\longrightarrow\W(V)$ in the basis of
 $\W(V)$ given by the orthogonal idempotents $\pr_\e$.

 \begin{Lemma}[Eigenvalues of the Classifying Endomorphism]
 \hfill\label{evc}\break
  Consider the decomposition of the tensor product $T\otimes T\,=\,
  \oplus_\alpha\,W_\alpha$ into irreducible summands. The classifying
  endomorphisms $K^\h$ are diagonalizable on $\Hom_\g(T\otimes T,\End V)$
  with eigenspaces $\Hom_\g(\,W_\alpha,\End V\,)\,\subset\,\Hom_\g
  (T\otimes T,\End V)$ for relevant $\alpha$ and eigenvalues:
  $$
   \kappa_{W_\alpha}
   \;\;=\;\;
   \frac12\,\Cas^{\L^2}_{W_\alpha}\;-\;\Cas^{\L^2}_T\ .
  $$
  In particular the classifying endomorphisms $K^\h$ acts as $K^\h(F)=
  \sum_\alpha\kappa_{W_\alpha}\,F|_{W_\alpha}$ on the space of Weitzenb\"ock
  formulas $\W(V)=\Hom_\g(T\otimes T,\End V)$.
 \end{Lemma}

 \proof
 From the very definition of $K^\h$ we see that it acts by precomposition
 with the map $-\sum X_\alpha \otimes X_\alpha$ in the interpretation
 $\W(V)\,=\,\Hom_\g(T\otimes T,\End V)$ of the space of Weitzenb\"ock
 formulas. The argument used in the proof of Corollary \ref{confW} shows
 that $K^\h$ is actually a difference of Casimir operators leading to the
 stated formula for its eigenspaces and eigenvalues.
 \qed

 \pfill
 In the case of the holonomies $\so_n, \g_2$ and $\spin_7$ we have
 $T\otimes T = \C \oplus \S^2_0T \oplus \g \oplus \g^\perp$ and using
 Lemma \ref{evc} we find the following $K$--eigenvalues:
 \begin{equation}
  \label{kt}
  \hbox{\begin{tabular}{|l|c|c|c|c|}
   \hline & & &  & \\[-4mm]
    & $\kappa_\C$ & $\; \kappa_{\S^2_0T}\;$
    & $\;\kappa_\g\;$ & $\;\kappa_{\g^\perp}\;$ \\[1mm]
   \hline & & & & \\[-3mm]
    $\so_n$ & $-(n-1)$ &
    $1$ & $-1$ & --- \\[1mm]
    $\g_2$ & $-4$ & $\frac{2}{3}$ & $0$ & $-2$ \\[1mm]
    $\spin_7$ & $-\frac{21}{4}$  & $\frac{3}{4}$ & $-\frac{1}{4}$&
    $-\frac{9}{4}$ \\[1mm]
   \hline
  \end{tabular}}
 \end{equation}
 Note that all these $K$--eigenvalues are different and consequently the twist
 $\tau$ is a polynomial in the classifying endomorphism $K$. Moreover a given
 invariant homomorphism $F\in\Hom_\g(T\otimes T, \End V)$ is an eigenvector
 of $K$ if and only if $F$ is different from zero one precisely one summand
 $W_\alpha \subset T\otimes T$, so $\1$ and $B$ are clearly $K$--eigenvectors.

 \begin{Lemma}[Properties of the Classifying Endomorphism]
 \label{elem}\hfill\break
  The classifying endomorphism $K:\W(V)\, \longrightarrow \,\W(V)$ is a
  symmetric endomorphism commuting with the twist map $\tau$ on the space
  $\W(V)$ of Weitzenb\"ock formulas equipped with the scalar product
  $\<F,\tilde F>\,:=\,\frac1{\dim\,V}\tr_{T\otimes V}\,F\tilde F$. The
  special endomorphisms $\,\1$ and $B$ for the same ideal $\h\,\subset\,\g$
  are $K$--eigenvectors:
  $$
   K(\,\1\,) \;\;=\;\;\Cas^{\L^2}_T\,\1
   \qquad\qquad
   K(\,B\,)  \;\;=\;\;(\,\Cas^{\L^2}_T\,-\,\frac12\,\Cas^{\L^2}_\h)\, B
  $$
 \end{Lemma}

 \proof
 The symmetry of $K$ is a trivial consequence of equation
 (\ref{ksym}) in the form
 $$
  \<K(\,F\,),\tilde F>\;\;=\;\;
  \frac1{\dim\,V}\,\sum_\nu\tr_{T\otimes V}
  \Big(\,(X_\nu\otimes\id)\,F\,(X_\nu\otimes\id)\,\tilde F\,\Big)
 $$
 and the cyclic invariance of the trace, moreover $K$ commutes with
 $\tau$ by definition. Coming to the explicit determination
 of $K(\1)$ and $K(B)$ we observe that the unit of $\End_\g(T\otimes V)$
 becomes the equivariant map $\1(a\otimes b)\,=\,\<a,b>\id_V$ in
 $\Hom_\g(T\otimes T,\,\End\,V)$ and so:
 $$
  (\,K\1\,)(a\otimes b)
  \;\;=\;\;
  -\sum_\nu\,\<X_\nu a,X_\nu b>\,\id_V
  \;\;=\;\;
  \sum_\nu\,\< a, \,X_\nu^2 b>\,\id_V
  \;\;=\;\;
  \Cas^{\L^2}_T\,\1(a\otimes b)\ .
 $$
 The conformal weight operator $B$ considered as an element of
 $\Hom_\g(T\otimes T, \End\,V)\,$ lives by definition in the
 eigenspace $\Hom_\g(\g,\End\,V)\,$ for the eigenvalue $-\frac12
 \,\Cas^{\L^2}_\g\,+\,\Cas^{\L^2}_T $ of $K$, where $\Cas^{\L^2}_\g$
 is the Casimir eigenvalue of the adjoint representation.
 \qed

 \pfill
 On a manifold with holonomy algebra $\g \subset \so_n \cong \L^2T $
 the Riemannian curvature tensor takes values in $\g$,~i.e. it can
 be considered as an element of $\S^2\g$. This fact has the following
 important consequence:

 \begin{Proposition}[Bochner Identities]
 \hfill\label{bi}\break
  Suppose $F\,\in\,\W(V)$ is an invariant homomorphism $T\otimes T
  \longrightarrow\End\,V$ factoring through the projection to the
  orthogonal complement $\, \g^\perp\subset\L^2T\subset T\otimes T$
  of $\g\subset\L^2T$. The curvature expression $\,F(\,\nabla^2\,)\,=\,0$
  vanishes regardless of what the curvature tensor $R$ is.
 \end{Proposition}

 \noindent
 We will call a Weitzenb{\"o}ck formula $F\in\W(V_\l)$ corresponding to
 an invariant homomorphism $T\otimes T\longrightarrow\End\,V$ factoring
 through $\,\g^\perp\,$ a {\it Bochner identity}. Writing such an invariant
 homomorphism $F$ in terms of the basis $\{\pr_\e\}$ as $F\,=\,\sum f_\e\,
 \pr_\e$ we get the following explicit form of the Bochner identity:
 $$
  \sum_\e \, f_\e\,\,T^*_\e \, T^{\phantom{*}}_\e \,\,=\,\,0 \ .
 $$
 The Bochner identities of $\G_2$-- and $\Spin(7)$--holonomy correspond
 to eigenvectors of the classifying endomorphism $K$ for the eigenvalues
 $-2$ and $\frac{9}{4}$ respectively. Since the zero weight space of
 $\g^\perp$ is in both cases one-dimensional it follows from
 Lemma~\ref{generic} that
 \begin{equation}\label{estimate1}
 \dim \Hom_\g(\g^\perp, \End V_\l)\;\; \le \;\; 1 \ ,
 \end{equation}
 i.~e.~there is at most one Bochner identity. Moreover the $K$--eigenvector
 $\1\in\End_\g(T\otimes V_\l)$ spans the $K$--eigenspace $\,\Hom_\g(\C,\End\,
 V_\l)\cong\C$. Because the zero weight space of $\g$ itself is the fixed
 Cartan subalgebra $\t\subset\g$, an application of Lemma \ref{generic}
 results in the estimates:
 $$
  \dim \Hom_{\g_2}(\g_2, \End V_\l)\;\; \le \;\; 2
  \qquad\qquad\quad
  \dim \Hom_{\spin_7}(\spin_7, \End V_\l)\;\; \le \;\; 3 \ .
 $$
\section{The Recursion Procedure for $\SO(n),\G_2$ and $\Spin(7)$}
 The definitions of the conformal weight operator $B$ and the classifying
 endomorphism $K$ given in the previous section are very similar. Given
 this similarity it should not come as a surprise that the actions of $B$
 and $K$ on the space $\W(V)$ of Weitzenb\"ock formulas obey a simple
 relation, which is the corner stone of the treatment of Weitzenb\"ock
 formulas proposed in this article. In the present section we first prove
 this relation and then use it to construct recursively an basis of
 $K$--eigenvectors of $\W(V_\l)$ for the holonomy groups $\SO(n)$, $\G_2$
 and $\Spin(7)$. A different interpretation of the same relation between
 $B$ and $K$ is studied in the next section concerning K\"ahler manifolds.
\subsection{The basic recursion procedure}
 Recall that the twist $\tau$ is defined in the interpretation $\W(V)\,=\,
 \Hom_\g(T\otimes T\otimes V,V)$ of the space of Weitzenb\"ock formulas as
 linear maps $T\otimes T\otimes V\longrightarrow V$ by precomposition with
 the endomorphism $\tau:\,a\otimes b\otimes v\longmapsto b\otimes a
 \otimes v$. Generalizing this precomposition we observe that $\W(V)$ is
 a right module over the algebra $\End_\g(T\otimes T\otimes V)$ containing
 $\tau$. Interestingly both the classifying endomorphism $K$ and
 the (right) multiplication by the conformal weight operator $B$ are induced
 by precomposition with elements in $\End_\g(T\otimes T\otimes V)$, too,
 say $K$ is the precomposition with the $\g$--invariant endomorphism
 $$
  K:\;\;T\otimes T\otimes V\;\longrightarrow\;T\otimes T\otimes V,
  \qquad
  a\otimes b\otimes v\;\longmapsto\;-\sum_\nu X_\nu a\otimes X_\nu b\otimes v
 $$
 while (right) multiplication by $B$ is precomposition with the
 $\g$--invariant endomorphism
 $$
  B:\;\;T\otimes T\otimes V\;\longrightarrow\;T\otimes T\otimes V,
  \qquad
  a\otimes b\otimes v\;\longmapsto\;-\sum_\nu a\otimes X_\nu b\otimes X_\nu v
 $$
 by Fegan's Lemma \ref{fegan}. From this description of the action of the
 classifying endomorphism $K$ and right multiplication of $B$ on $\W(V)$
 we immediately conclude:
 $$
  K\;+\;B\;+\;\tau B\tau\;\;=\;\;
  -\frac12(\;\Cas^{\L^2}\,-\,2\Cas^{\L^2}_T\,-\,\Cas^{\L^2}_{V_\l}\;)
 $$
 Applying once again Schur's Lemma this relation implies our basic
 Recursion Formula:

 \begin{Theorem}[Recursion Formula]
 \hfill\label{rec}\break
  Let $V_\l$ be an irreducible representation of the holonomy algebra
  $\g$. Then the action of $\,K,\,B$ and $\tau$ on $\W(V_\l)\,=\,\Hom_\g(
  T\otimes T\otimes V_\l,V_\l)$ by precomposition satisfies:
  $$
   K\;+\;B\;+\;\tau B\tau\;\;=\;\;
   \Cas^{\L^2}_T\;\;=\;\;-\,2\,\frac{\dim\,\h}{\dim\,T}
  $$
 \end{Theorem}

 \noindent
 We will now explain how this theorem yields a recursion formula for
 $K$--eigenvectors. In fact given an eigenvector $F\,\in\,\W(V)$ for
 the twist $\tau$ and the classifying endomorphism $K$ with eigenvalues
 $t$ and $\kappa$, i.e. $\tau F \,=\,tF$ and $KF\,=\,\kappa F$, the
 Recursion Formula allows us to produce a new $\tau$--eigenvector
 $\widehat F$ with eigenvalue $-t$. This simple prescription suffices
 to obtain a complete eigenbasis for $\W(V)$ of $\tau$-- and actually
 $K$--eigenvectors in general Riemannian geometry $\g\,=\,\so_n$ and,
 with some modifications, also in the exceptional cases $\g\,=\,\g_2$
 and $\g\,=\,\spin_7$. The quaternionic K\"ahler case can be dealt with
 similarly.

 \begin{Corollary}[Basic Recursion Procedure]
 \hfill\label{recursion2}\break
  Let $F\,\in\,\W(V)$ be an eigenvector for the involution $\,\tau\,$ and
  the classifying endomorphism $K$ of an ideal $\h\,\subset\,\g$, i.~e.~
  $ K(F)\,=\,\kappa F$ and $\tau(F)\,=\,\pm F$. The new Weitzenb\"ock formula
  $$
   F_{\mathrm{new}}
   \;\;:=\;\;
   (\,B\;-\;\frac{\Cas^{\L^2}_T-\kappa}2\,)\circ F
  $$
  is again a $\tau$--eigenvector in $\W(V)$ with $\,\tau(F_{\mathrm{new}})
  \,=\,\mp F_{\mathrm{new}}$. In particular we find:
  $$
   \1_{\mathrm{new}}\;\;=\;\;B\qquad\quad\mbox{and}\quad\qquad
   B_{\mathrm{new}}\;\;=\;\;B^2 \;-\; \frac14\Cas^{\L^2}_\g \, B\ .
  $$
 \end{Corollary}

 \noindent
 \proof
 We observe that the Recursion Formula \ref{rec} in the form $\tau B\tau
 \,=\,\Cas^{\L^2}_T-K-B$ implies under the assumptions $K(F)\,=\,\kappa F$
 and $\tau(F)\,=\,\pm F$ that
 $$
  \pm\tau(\,BF\,)\;\;=\;\;(\,\Cas^{\L^2}_T-\kappa\,)\,F\;-\;BF
 $$
 and consequently:
 $$
  \pm\tau(\,BF\,-\,\frac{\Cas^{\L^2}_T-\kappa}2\,F\,)
  \;\;=\;\;-\,(\,BF\,-\,\frac{\Cas^{\L^2}_T-\kappa}2\,F\,)\ .
 $$
 The formulas for $\,\1_{\mathrm{new}}\,$ and $\,B_{\mathrm{new}}\,$
 are immediate consequences of Lemma \ref{elem}.
 \qed

 \pfill
 Recall that a $K$--eigenvector is automatically a $\tau$--eigenvector.
 In general however the new Weitzenb\"ock formula $F_{\mathrm{new}}\,\in
 \,\W(V)$ does not need to be an eigenvector for $K$ again and there is
 no way to iterate the recursion. Nevertheless it is possible to avoid
 termination of the recursion procedure for most of the irreducible
 non--symmetric holonomy algebras by using appropriate projections.

 We note that any  $+1$--eigenvector of $\tau$ orthogonal to $\1$ is
 already $K$--eigenvector in $\Hom_\g(\S^2_0T,\End {V_\l})$. This is
 due to the fact that $\1$ spans the other summand of the
 $+1$--eigenspace of $\tau$. In particular the orthogonal projection
 of $ B_{\mathrm{new}}$ onto the orthogonal complement of $\1$, i.~e.~
 the polynomial  $B^2 -\frac14\Cas^{\L^2}_\g B +\frac2n\Cas^{\L^2}_{V_\l}$,
 is a $K$--eigenvector in $\Hom_\g(\S^2_0T, \End {V_\l})$. More generally
 we have:

 \begin{Corollary}[Orthogonal recursion procedure]
 \hfill\label{orp}\break
  Let $p_0(B), \ldots, p_k(B)$ be a sequence of polynomials obtained
  by applying the Gram Schmidt orthonormalization procedure to the powers
  $1,\,B,\,B^2,\,\ldots,\,B^k$ of the conformal weight operator $B$. If
  all these polynomials are $\tau$--eigenvectors and $p_k(B)$ is moreover
  a $K$--eigenvector, then the orthogonal projection $p_{k+1}(B)$ of
  $B^{k+1}$ onto the orthogonal complement of the span of $1,\,B,\,\ldots,
  \,B^k$ is a again a $\tau$--eigenvector.
 \end{Corollary}

 \proof
 Since $p_k(B)$ is a $K$--eigenvector the basic recursion procedure
 shows the existence of a polynomial in $B$ of degree $k+1$, which is a
 $\tau$--eigenvector. It follows that with $\mathrm{span}\{1,\,B,\,\ldots,
 \,B^k\} = \mathrm{span}\{p_0(B),\,\ldots,\,p_k(B)\}$ also $\mathrm{span}
 \{1,\,B,\,\ldots,\,B^{k+1}\}$ is $\tau$--invariant. Clearly the orthogonal
 projection $p_{k+1}(B)$ of $B^{k+1}$ onto the orthogonal complement of
 $\mathrm{span}\{1,\,B,\,\ldots,\,B^k\}$ is a polynomial in $B$ of degree
 $k+1$ with:
 $$
  \mathrm{span}\{1,\,B,\,B^2,\,\ldots,\,B^k\} \,\oplus\,
  \C\,p_{k+1}(\,B\,)
  \;\;=\;\;
  \mathrm{span}\{1,\,B,\,B^2,\,\ldots,B^k,B^{k+1}\} \ .
 $$
 Now the involution $\tau$ is symmetric with respect to the scalar
 product on $\End_\g(T\otimes V)$ and so the orthogonal complement of a
 $\tau$--invariant space is again $\tau$--invariant.
 \qed
\subsection{Computation of $B$--eigenvalues for $\SO(n),\G_2$ and $\Spin(7)$}
\label{eigenvalues}
 In this section we will compute the $B$--eigenvalues for the holonomies
 $\SO(n),\G_2$ and $\Spin(7)$ by applying the explicit formula of Corollary~\ref{confW}. In
 particular we will see that with only two exceptions all $B$--eigenvalues
 are pairwise different. This information is relevant in the proof of
 Proposition~\ref{quotient}.

 \pfill
 {\it The $\SO(n)$--case.} Recall that in Section \ref{tables} we fixed the
 notation for the fundamental weights $\o_1,\,\ldots,\,\o_r$ and the weights
 $\pm\e_1,\,\ldots,\,\pm\e_r$ of the defining representation $\R^n$ of $\SO(n)$
 with $r:=\lfloor\frac n2\rfloor$. Moreover the scalar product $\<,>$ on the
 dual of a maximal torus was chosen so that the weights $\pm\e_1,\,\ldots,\,
 \pm\e_r$ are an orthonormal basis. A highest weight can be written $\l\,=\,
 \l_1\o_1+\ldots+\l_r\o_r\,=\,\mu_1\e_1+\ldots+\mu_r\e_r$ with integral
 coefficients $\l_1,\,\ldots,\,\l_r\geq 0$ and coefficients $\mu_1,\,\ldots,
 \,\mu_r$, which are either all integral or all half--integral and decreasing.
 Independent of the parity of $n$ the conformal weights are
 $$
  b_{+\e_k} \;=\;  \mu_k \,-\, k       \,+\,1,\qquad
  b_{-\e_k} \;=\; -\mu_k \,-\, n\,+\,k \,+\,1,\qquad
  b_0       \;=\; -\,r
 $$
 according to Corollary \ref{confW}, where the zero weight only appears for
 $n$ odd. With only a few exceptions the conformal weights are totally ordered
 and thus pairwise different. In the case $n$ odd the coefficients $\mu_1,\,
 \ldots,\,\mu_r$ are decreasing in the sense $\mu_1\geq\mu_2\geq\ldots\geq
 \mu_r\geq 0$ so that we find strict inequalities
 $$
  b_{-\e_1}\,<\,b_{-\e_2}\,<\,\ldots\,<\,b_{-\e_r}
  \,\leq\,b_0\,<\,b_{+\e_r}\,<\,\ldots\,<\,b_{+\e_1}
 $$
 unless $\mu_r=0$ or equivalently $\l_r=0$. In the latter case $b_{-\e_r}\leq
 b_0$ happens to be an equality, however Lemma \ref{relw} tells us that the
 zero weight is irrelevant for highest weights $\l$ with $\l_r=0$. Without
 loss of generality we may thus assume all conformal weights to be different
 for $n$ odd. Similar considerations in the case of even $n$ based on the
 inequalities $\mu_1\geq\mu_2\geq\ldots\geq\mu_{r-1}\geq|\mu_r|$ satisfied
 by the coefficients $\mu_1,\,\ldots,\,\mu_r$ of $\l$ lead to
 $$
  b_{-\e_1}\,<\,b_{-\e_2}\,<\,\ldots\,<\,b_{-\e_{r-1}}\,<\,
  \{\;b_{-\e_r}\;,\;b_{+\e_r}\;\}\,<\,\b_{+\e_{r-1}}\,<\,\ldots
  \,<\,b_{+\e_2}\,<\,b_{+\e_1}
 $$
 where nothing specific can be said about the relation between $b_{-\e_r}$
 and $b_{+\e_r}$ due to $b_{+\e_r}-b_{-\e_r}\,=\,2\mu_r$. This should not
 come too surprising as the outer automorphism of $\SO(n)$ with $n$ even
 acts on $\t^*$ as a reflection along the hyperplane $\mu_r=0$.

 \pfill
 {\it The $\G_2$--case.} We write the highest weight as $\l = a\o_1+b\o_2$
 with integers $a,\,b\geq0$ and use the scalar product defined in Section
 \ref{tables} by setting $\<\e_1,\e_1>=1=\<\e_2,\e_2>$ and $\<\e_1,\e_2>=
 \frac12$, equivalently $\<\o_1,\o_1>=1$, $\<\o_2, \o_2>=3$ and $\<\o_1,
 \o_2>=\frac32$. According to Corollary \ref{confW} the conformal weight
 of the zero weight is given by $b_0=-\,2$, similarly:
 $$
  b_{\pm\e_1} \;=\; -\,(\tfrac53 \mp \tfrac53) \;\pm \;(\tfrac23 a + b),\quad
  b_{\pm\e_2} \;=\; -\,(\tfrac53 \mp \tfrac43) \;\pm \;(\tfrac13 a + b),\quad
  b_{\pm\e_3} \;=\; -\,(\tfrac53 \mp \tfrac13) \;\pm \;\tfrac13 a \ .
 $$
 Again all conformal weights or $B$--eigenvalues are pairwise different and
 totally ordered:
 $$
  b_{-\e_1}\,<\,b_{-\e_2}\,<\,b_{-\e_3}\,<\,b_0
  \,<\,b_{+\e_3}\,<\,b_{+\e_2}\,<\,b_{+\e_1} \ .
 $$

 \pfill
 {\it The $\Spin(7)$--case.} Using the fundamental weights $\o_1,\,\o_2,
 \,\o_3$ and the scalar product $\<,>$ introduced in Section \ref{tables}
 in terms of the weights $\pm\eta_1,\,\pm\eta_2,\,\pm\eta_3$ of the
 represention $\R^7$ we write the highest weight $\l=a\o_1+b\o_2+c\o_3$
 with integers $a,\, b,\,c\geq0$ and compute
 $$
  \begin{array}{lclclcl}
   b_{\pm \e_1}
   &=&
   -\,(\tfrac94 \mp \tfrac94)\;\pm\;(\tfrac12a+b+\tfrac34c)\ ,
   &\quad&
   b_{\pm \e_2}
   &=&
   -\,(\tfrac94 \mp \tfrac74)\;\pm\;(\tfrac12a+b+\tfrac14c)\ ,
   \\[2ex]
   b_{\pm \e_3}
   &=&
   -\,(\tfrac94 \mp \tfrac34)\;\pm\;(\tfrac12a\quad\;\;+\tfrac14c)\ ,
   &&
   b_{\pm \e_4}
   &=&
   -\,(\tfrac94 \mp \tfrac14)\;\pm\;(\tfrac12a\quad\;\;-\tfrac14c)\ .
  \end{array}
 $$
 In this case we obtain the inequalities
 $$
  b_{-\e_1}\,<\,b_{-\e_2}\,<\,b_{-\e_3}\,<\,\{\;b_{-\e_4}\;,\;
  b_{+\e_4}\;\}\,<\,b_{+\e_3}\,<\,b_{+\e_2}\,<\,b_{+\e_3}\ .
 $$
 however the difference $b_{+\e_4}-b_{-\e_4}\,=\,a-\frac12c+\frac12$ does
 not allow to draw conclusions about the relation between $b_{-\e_4}$ and
 $b_{+\e_4}$. In particular for a heighest weight $\l$ with $c\,=\,2a+1$
 the two conformal weights $b_{-\e_4}=b_{+\e_4}$ agree.
\subsection{Basic Weitzenb\"ock formulas for $\SO(n),\,\G_2$ and $\Spin(7)$}
 In this section we make the recursion procedure of Corollary~\ref{recursion2}
 explicit for the holonomy groups $\SO(n),\G_2$ and $\Spin(7)$. Let us start
 with the generic Riemannian holonomy algebra $\g\,=\,\so_n$ with only a
 single non--trivial ideal $\h\,=\,\g$. According to table (\ref{kt}) its
 classifying endomorphism $K$ has eigenvalues $(n-1),\;1\,$ and $\,-1$
 with eigenspaces $\C\1$, the orthogonal complement of $\1$ in the
 $\tau$--eigenspace for the eigenvalue $1$ and the $\tau$--eigenspace for
 the eigenvalue $-1$ respectively. The orthogonal projection of every
 $\tau$--eigenvector to the orthogonal complement of $\1$ is thus a
 $K$--eigenvector. Consequently we can modify the recursion procedure
 such that it associates to an eigenvector $F$ for $\tau$ of eigenvalue
 $-1$ the $K$--eigenvector
 $$
  F_{\mathrm{new}}\;\;:=\;\;(\,B\;+\;\frac{n-2}2\,)\,F\;-\;\frac1n\,
  \<BF,\1>\,\1
 $$
 for the eigenvalue $1$, while a $\tau$--eigenvector $F$ for the eigenvalue
 $+1$ orthogonal to $\1$ is mapped to the $K$--eigenvector:
 $$
  F_{\mathrm{new}}\;\;:=\;\;(\,B\;+\;\frac n2\,)\,F
 $$
 for the eigenvalue $-1$. For an irreducible representation $V_\l$
 of $\so_n$ we thus get a sequence $p_0(B),\,p_1(B),\,\ldots$ of
 $K$--eigenvectors in $\W(V_\l)$ defined recursively by $p_0(B)\,:=\,\1,
 \,p_1(B)\,:=\,B$ and $p_{k+1}(B)\,:=\,
 (p_k(B))_{\mathrm{new}}$ for $k\geq 1$.
 Evidently the different $p_k(B)$ are  polynomials of
 degree $k$ in $B$ and so the eigenvectors $p_0(B),\,\ldots,\,p_{d-1}(B)
 \,\in\,\W(V_\l)$
 with $d\,:=\,\mathrm{deg}\,\mathrm{min}\,B$ are necessarily linearly
 independent. According to Proposition~\ref{quotient} we always get a
 complete basis of $\tau$--eigenvectors with the exception of the
 case $\g= \so_{2r}$ and a representation $V_\l$ with
 $\l_{r-1}=\l_r$. Here we still have to add a $K$--eigenvector
 $F_{spin}$ spanning the orthogonal complement of the image of $\C[\,B\,]$ in
 $\W(V_\l)$.

 Note that the polynomials $p_{2k+1}(B),\, k=0,1,\ldots$ are in the
 $-1$--eigenspace of $\tau$. Hence the corresponding Weitzenb\"ock
 formulas give a pure curvature term. Let $N$ the number of
 irreducible components of $T\otimes V_\l$, then there are
 $\lfloor \frac{N}{2} \rfloor$ linearly independent equations of this
 type. This result, which is clear from our construction, was proved
 for the first time in \cite{branson}. The first  eigenvectors in this
 sequence are $\,p_0(B)\,=\,\1$ and $p_1(B)\,=\,B\,$ as well as:
 \begin{eqnarray}
  p_2(B)
  &=&
  B^2\;+\;\tfrac{n-2}2\,B\;+\;\tfrac2n\,\Cas^{\L^2}_{V_\l}
  \\[1.8ex]
  p_3(B)
  &=&
  B^3\;+\;(n-1)\,B^2\;+\;(\tfrac{2}n\,\Cas^{\L^2}_{V_\l} +
  \tfrac{n(n-2)}{4})\,B\;+\;\Cas^{\L^2}_{V_\l}
  \label{gleichung2}
 \end{eqnarray}

 \noindent
 Essentially the same procedure can be used in the case $\g\,=\,\g_2$ to
 compute a complete $K$--eigenbasis for the space $\W(V_\l)$ for an irreducible
 $\g_2$--representation $V_\l$. Again there is only one non--trivial ideal
 $\h\,=\,\g_2$ and hence only a single classifying endomorphism $K$.
 However the $\tau$--eigenspace in $\W(V_\l)$ for the eigenvalue $-1$
 decomposes into two $K$--eigenspaces. The recursion procedure gives the
 $K$--eigenvectors
 \begin{equation}\label{g2-1}
  p_0(B)
  \;\;=\;\; \1,\quad
  p_1(B) \;\;=\;\; B,\quad
  p_2(B) \;\;=\;\; B^2 + 2B + \tfrac27 \,\Cas^{\L^2}_{V_\lambda} \ .
 \end{equation}
 Using the recursion procedure again gives a polynom of degree 3 in
 $B$. Projecting it onto the orthogonal complement of $B$ we obtain
 \begin{equation}\label{g2-2}
  p_3(B)
  \;\;=\;\;
  B^3\,+\,\tfrac{13}{3}\,B^2
  \,+\,(\tfrac12\,\Cas^{\L^2}_{V_\lambda}+4)\,B
  \,+\,\tfrac23\,\Cas^{\L^2}_{V_\lambda} \ .
 \end{equation}
 We will see in Theorem~\ref{bochig2} that $p_3(B)$ is in fact a
 $K$--eigenvector for the eigenvalue $-2$, in other words $p_3(B)\,\in\,
 \Hom_{\g_2}(\g_2^\perp,V_\l)$ is a Bochner identity. Due to the estimate
 (\ref{estimate1}) any $\tau$--eigenvector orthogonal to $\1$ and $F_3$
 is a $K$--eigenvector and so we may obtain a complete eigenbasis $p_0(B),
 \,\ldots,\,p_6(B)$ in the $\G_2$--case by applying the Gram--Schmidt
 orthogonalization to the powers of $B$ and using Corollary~\ref{orp}.

 In order to make the generalized Bochner identity corresponding to the
 polynomial $p_3(B)$ explicit we recall that its coefficients as a linear
 combination of the basis projections $\pr_\e$ are the value of the polynomial
 $p_3$ at the corresponding $B$--eigenvalues $b_\e$. Substituting the explicit
 formulas for $b_\e$ and for $\Cas^{\L^2}_{V_\l}$ (c.f. Remark~\ref{g2cas})
 we obtain:
 \begin{equation}\label{g2bx}
  \begin{array}{rcccc}
   F_{\mathrm{Bochner}}
   \;\;:=\;\;27\,p_{3}(B)\;\;=\;\;
   +\!&\! a   \!&\!\!(a+3b+3)\!&\!\!(2a+3b+4)\!&\!\pr_{+\e_1}\cr
   -\!&\!(a+2)\!&\!\!(a+3b+5)\!&\!\!(2a+3b+6)\!&\!\pr_{-\e_1}\cr
   -\!&\!(a+2)\!&\!\!(a+3b+3)\!&\!\!(2a+3b+4)\!&\!\pr_{+\e_2}\cr
   +\!&\! a   \!&\!\!(a+3b+5)\!&\!\!(2a+3b+6)\!&\!\pr_{-\e_2}\cr
   -\!&\! a   \!&\!\!(a+3b+5)\!&\!\!(2a+3b+4)\!&\!\pr_{+\e_3}\cr
   +\!&\!(a+2)\!&\!\!(a+3b+3)\!&\!\!(2a+3b+6)\!&\!\pr_{-\e_3}\cr
   +\!&  6(\,a^2\,+\!\!\!&\!\!\!3b^2\,+\,3ab\,+\!\!\!&\!\!5a\,+
   9b\,+\,6)\!\!\!&\!\!\!\pr_0
  \end{array}
 \end{equation}

 \pfill
 Eventually let us discuss the example of $\g\,=\,\spin_7$. Here
 the modified recursion procedure gives the three $K$--eigenvectors
 \begin{equation}\label{s7-1}
  p_0(B)
  \;\;=\;\; \1,\quad
  p_1(B)
  \;\;=\;B,\quad
  p_2(B) \;\;=\;\; B^2+\tfrac52\,B +\tfrac14\,\Cas^{\L^2}_{V_\l} \ .
 \end{equation}
 and a $\tau$--eigenvector for the eigenvalue $-1$, which is of third
 order as a polynomial in $B$. After projecting it onto the orthogonal
 complement of $B$ we obtain:
 \begin{equation}\label{s7-2}
  p_3(B)
  \;\;=\;\;
  B^3
  \;+\;\tfrac{11}{2}\,B^2
  \;+\;\tfrac{1}{2\,\Cas^{\L^2}_{V_\l}}\,
   (\Cas^{[4]}_{V_\l}+\tfrac{55}{2}\,\Cas^{\L^2}_{V_\l} )\,B
  \;+\;\tfrac34\,\Cas^{\L^2}_{V_\l} \ ,
 \end{equation}
 where $\Cas_{V_\l}^{[4]}$ is the eigenvalue of the higher casimir
 $\Cas^{[4]}$ on the irreducible representation $V_\l$. Its explicit
 value is given in the appendix in Remark~\ref{spin7cas}.

 However in difference to the $\g_2$--case this is no
 $K$--eigenvector. Indeed in Section~\ref{kis} we will see
 that the space of polynomials in $B$ of degree at most 3 is not
 invariant under $K$. Hence there cannot be a further $K$--eigenvector
 expressible as a polynomial of order $3$ in $B$. In general the other
 $K$--eigenvectors are polynomials of degree 7 in $B$. They are too
 complicated to be written down, but surprisingly the $K$--eigenvector
 for the eigenvalue $-\frac94$, i.~e.~the Bochner identity, for a
 representation of highest weight $\l = a\o_1 +b\o_2 +c\o_3$ has the
 following simple explicit expression:
 \begin{equation}\label{s7b}
  \begin{array}{rcccc}
   F_{\mathrm{Bochner}}\;\;=\;\;
   +&\!  c  \!&\!(2b+c+2)\!&\!(2a+2b+c+4)\!&\!\pr_{+\e_1}\cr
   -&\!(c+2)\!&\!(2b+c+4)\!&\!(2a+2b+c+6)\!&\!\pr_{-\e_1}\cr
   -&\!(c+2)\!&\!(2b+c+2)\!&\!(2a+2b+c+4)\!&\!\pr_{+\e_2}\cr
   +&\!  c  \!&\!(2b+c+4)\!&\!(2a+2b+c+6)\!&\!\pr_{-\e_2}\cr
   -&\!  c  \!&\!(2b+c+4)\!&\!(2a+2b+c+4)\!&\!\pr_{+\e_3}\cr
   +&\!(c+2)\!&\!(2b+c+2)\!&\!(2a+2b+c+6)\!&\!\pr_{-\e_3}\cr
   +&\!(c+2)\!&\!(2b+c+4)\!&\!(2a+2b+c+4)\!&\!\pr_{+\e_4}\cr
   -&\!  c  \!&\!(2b+c+2)\!&\!(2a+2b+c+6)\!&\!\pr_{-\e_4}
  \end{array}
 \end{equation}
 This formula is proved in Theorem~\ref{bochis7}. Note that the
 coefficients of $\pr_{+\e_4}$ and $\pr_{-\e_4}$ are different.
 Hence in the critical case with $c = 2a+1$, i.~e.~where
 $b_{+\e_4} = b_{-\e_4}$, this $K$--eigenvector $F_{Bochner}$
 spans the space orthogonal to $\C[B]$ in $\End_\g(T \otimes V_\l)$.
\section{The Weitzenb\"ock Machine for K\"ahler Holonomies}
\label{kholon}
 The Weitzenb\"ock machine for the K\"ahler holonomy groups
 $\mathbf{U}(n)$ and $\SU(n)$ is remarkably different to the machine
 for the other non--symmetric irreducible holonomies. Perhaps
 the most distinctive characteristic of K\"ahler holonomy is
 that the complexified holonomy representation $T\,=\,\bar E\oplus E
 \,=\,T^{1,0}\oplus T^{0,1}$ is not irreducible. Consequently the set of
 generalized gradients $\{T_\e\}$ on sections of an irreducible vector bundle
 $V_\l M$ falls apart into two subsets, the sets of first order parallel
 differential operators factorizing over the complementary projections
 $\pr^{0,1}$ and $\pr^{1,0}$ from $T\otimes V_\l$ to $T^{0,1}\otimes V_\l$
 and $T^{1,0}\otimes V_\l$ respectively. In turn the space of Weitzenb\"ock
 formulas on $V_\l$ is the direct sum
 \begin{equation}\label{hah}
  \W(\,V_\l\,)
  \;\;=\;\;
  \W^{1,0}(\,V_\l\,)\;\oplus\;\W^{0,1}(\,V_\l\,)
 \end{equation}
 of the spaces $\W^{1,0}(V_\l):=\End_\g(\,T^{0,1}\otimes V_\l\,)$
 and $\W^{0,1}(V_\l):=\End_\g(\,T^{1,0}\otimes V_\l\,)$ of holomorphic
 and antiholomorphic Weitzenb\"ock formulas. Evidently the space
 $\W^{1,0}(V_\l)$ is spanned by the projections $\pr_\e$ to the
 holomorphic relevant weights $\e\subset\l^{1,0}$ defined as the
 relevant weights among the weights $-\e_1,\ldots,-\e_n$ of $T^{0,1}$,
 in particular $\pr^{0,1}$ is the sum $\pr^{0,1}=\sum_{\e\subset\l^{1,0}}
 \pr_\e$ over all relevant holomorphic weights. Similarly the space
 $\W^{0,1}(V_\l)$ is spanned by the projections to the relevant weights
 $\e\subset\l^{0,1}$ among the antiholomorphic weights $+\e_1,\ldots,+\e_n$
 of $T^{1,0}$ with $\pr^{1,0}=\sum_{\e\subset\l^{0,1}}\pr_\e$. The
 isomorphisms $T^{0,1}\,=\,T^{1,0*}$ and $T^{1,0}\,=\,T^{0,1*}$ are
 reponsible for this apparently skewed notation.

 In the Weitzenb\"ock machine we are eventually interested in the
 matrices of the twist $\tau$ and the classifying endomorphism $K$ in
 the basis $\pr_\e$ of idempotents of $\W(V_\l)$, and a simple observation
 drastically simplifies this task in the case of K\"ahler holonomies.
 The subspaces $T^{1,0}$ and $T^{0,1}$ are isotropic subspaces
 of $T$ so that $\Hom_{\u_n}(T^{1,0}\otimes T^{1,0}\otimes V_\l,V_\l)$
 and $\Hom_{\u_n}(T^{0,1}\otimes T^{0,1}\otimes V_\l,V_\l)$ are trivial.
 Hence the decomposition (\ref{hah}) becomes
 $$
  \W^{1,0}(\,V_\l\,)
  \;\;=\;\;
  \Hom_{\u_n}(T^{1,0}\otimes T^{0,1}\otimes V_\l,V_\l)
  \qquad
  \W^{0,1}(\,V_\l\,)
  \;\;=\;\;
  \Hom_{\u_n}(T^{0,1}\otimes T^{1,0}\otimes V_\l,V_\l)
 $$
 in the interpretation $\W(V_\l)\,=\,\Hom_{\u_n}(T\otimes T\otimes V_\l,V_\l)$
 of the space of Weitzenb\"ock formulas. In this interpretation the
 endomorphisms $\tau$ and $K$ are defined by precomposition with
 $a\otimes b\otimes\psi\longmapsto b\otimes a\otimes\psi$ and $a\otimes
 b\otimes\psi\longmapsto-\sum_\alpha X_\alpha a\otimes X_\alpha b \otimes\psi$
 respectively so that $K$ preserves the decomposition (\ref{hah}), while
 $\tau$ interchanges $\W^{1,0}(V_\l)$ and $\W^{0,1}(V_\l)$. In difference to
 the holonomies $\SO(n)$, $\G_2$ and $\Spin(7)$ discussed above the twist
 $\tau$ is thus not a polynomial in the classifying endomorphism $K$ for
 $\mathbf{U}(n)$ and $\SU(n)$.

 \pfill
 Before calculating the $B$--eigenvalues let us recall that the weights
 $\pm\e_1,\,\ldots,\,\pm\e_n$ of the holonomy representation $T$ form
 an orthonormal basis for a scalar product $\<,>$ on the dual $\t^*$ of
 a maximal torus $\t\subset\u_n$. The decomposition $\u_n\,=\,i\R\oplus
 \su_n$ of $\u_n$ into center $i\R$ and simple ideal $\su_n$ is accompanied
 by an analoguous decomposition $\l\,=\,\l^{i\R}+\l^{\su}$ of a weight
 $\l\in\t^*$ along complementary orthogonal projections defined on the
 basis by:
 $$
  \e_k^{i\R}\;\;:=\;\;\frac1n(\;\e_1\,+\,\ldots\,+\,\e_n\;)
  \qquad\qquad
  \e_k^{\su}\;\;:=\;\;\e_k\;-\;\frac1n(\;\e_1\,+\,\ldots\,+\,\e_n\;)\ .
 $$
 For the calculations below it is useful to keep the following formulas
 in mind for all $k,l$:
 $$
  \<\e_k^{i\R},\e_l^{i\R}>\;\;=\;\;\frac1n
  \qquad\qquad
  \<\e_k^{i\R},\e_l^{\su}>\;\;=\;\;0
  \qquad\qquad
  \<\e_k^{\su},\e_l^{\su}>\;\;=\;\;\delta_{kl}\;-\;\frac1n\ .
 $$
 Interestingly the half sum $\rho\,=\,\frac{n-1}2\e_1+\frac{n-3}2\e_2+\ldots
 +\frac{1-n}2\e_n$ of positive roots of $\u_n$ is not equal to the sum of
 the fundamental weights $\o_1:=\e_1,\,\o_2:=\e_1+\e_2,\,\o_3:=\e_1+\e_2+\e_3$
 etc.~specified in Section~\ref{tables}. Nevertheless $\rho=\rho^{\su}$ is
 orthogonal to a weight $\l^{i\R}\in\t^*$ due to
 $$
  \rho\;\;=\;\;\o^{\su}_1\;+\;\o^{\su}_2\;+\;\ldots\;+\;\o^{\su}_n\,,
 $$
 in particular $\<\e_k^{\su},2\rho>=\<\e_k,2\rho>=n-2k+1$. Having established
 this much of the notation we need to verify that the chosen scalar product
 $\<,>$ is suitable for calculating Casimir eigenvalues. The problem is that
 there are two linearly independent scalar products on $\u_n$ and so we have
 to check the normalization condition (\ref{cnor}) for the Casimir eigenvalues
 for both ideals $i\R$ and $\su_n$ of $\u_n$ separately. Recalling
 that the holonomy representation $T$ is the complexification of a real
 irreducible representation $T_\R$ we observe that the Casimir eigenvalues
 of the normalized Casimirs $\Cas^{\L^2,i\R}$ and $\Cas^{\L^2,\su}$ of the
 two ideals $i\R$ and $\su_n$ are the same on $T, T^{1,0}$ and $T^{0,1}$.
 Taking the irreducible subspace $T^{1,0}$ of highest weight $\e_1$ for
 convenience we verify (\ref{cnor})
 $$
  \begin{array}{lclclcccl}
   \Cas^{\L^2,i\R}_T
   &=& \Cas^{\L^2,i\R}_{T^{1,0}}
   &=& -\;2\,{\displaystyle\frac{\dim\,i\R}{\dim\,T}}
   &=& -\,{\displaystyle\frac1n}
   &=& -\,\<\e_1^{i\R},\e_1^{i\R}>
   \\[3mm]
   \Cas^{\L^2,\su}_T
   &=& \Cas^{\L^2,\su}_{T^{1,0}}
   &=& -\;2\,{\displaystyle\frac{\dim\,\su}{\dim\,T}}
   &=& -\,{\displaystyle\frac{n^2-1}n}
   &=& -\,\<\e_1^{\su}+2\rho,\e_1^{\su}>
  \end{array}
 $$
 for the normalized Casimir:
 \begin{equation}\label{ckk}
  \Cas^{\L^2}_{V_\l}
  \;\;=\;\;-\;\<\l+2\rho,\l>
  \;\;=\;\;-\;\<\l^{i\R},\l^{i\R}>\;-\;\<\l^{\su}+2\rho,\l^{\su}>\,.
 \end{equation}
 Corollary \ref{confW} allows us to calculate the eigenvalues of the
 conformal weight operator $B$ aka conformal weights. With respect to
 the orthonormal basis $\e_1,\,\ldots,\,\e_n$ of $\t^*$ the
 highest weight of an irreducible representation $V_\l$ can be written
 $\l=\mu_1\e_1+\ldots+\mu_n\e_n$ with decreasing integral coefficients
 $\mu_1\geq\ldots\geq\mu_n$ so that the conformal weights are
 \begin{equation}\label{kaehler}
  b_{-\e_k}\;\;=\;\;-\mu_k\,+\,k\,-\,n
  \qquad\qquad
  b_{+\e_k}\;\;=\;\;\mu_k\,-\,k\,+\,1
 \end{equation}
 for the conformal weight operator $B$ of the full holonomy algebra $\u_n$.
 In particular both the holomorphic and antiholomorphic conformal weights
 are separately ordered
 $$
  b_{-\e_n} \,>\, b_{-\e_{n-1}} \,>\, \ldots \,>\, b_{-\e_1}
  \qquad\qquad
  b_{+\e_1} \,>\, b_{+\e_2} \,>\, \ldots \,>\, b_{+\e_n}\,,
 $$
 but we can not exclude the equality between the conformal weights for a
 holomorphic and an antiholomorphic weight, moreover there are certainly
 dominant integral weights $\l$ with $b_{-\e_k}=b_{+\e_{n-k}}$ for all
 $k=1,\ldots,n$. For this reason it is prudent to study the central conformal
 weight operator $B^{i\R}$ besides the conformal weight operator $B$ of the
 holonomy algebra $\u_n$, whose eigenvalues or central conformal weights
 $$
  b^{i\R}_{\e}
  \;\;=\;\;
  \frac12\,\Big(\,\<\l^{i\R}+\e^{i\R},\l^{i\R}+\e^{i\R}>
  \,-\,\<\l^{i\R},\l^{i\R}>\,-\,\<\e_1^{i\R},\e_1^{i\R}>\,\Big)
  \;\;=\;\;
  \<\l^{i\R},\e^{i\R}>
 $$
 equal $\pm\frac1n(\mu_1+\ldots+\mu_n)$ for holomorphic ($-$) and
 antiholomorphic ($+$) weights respectively. The algebra homomorphism
 $\Phi:\;\C[B^{i\R},B^{\su}]\longrightarrow\End_{\u_n}(T\otimes V_\l)$
 is thus surjective for all dominant integral weights $\l$ with
 $\l^{i\R}\neq 0$.

 In order to find the eigenspaces of the classifying endomorphism $K$ and
 their dimensions we recall that $K$ preserves the orthogonal subspaces
 $\W^{1,0}(V_\l)$ and $\W^{0,1}(V_\l)$ of holomorphic and antiholomorphic
 Weitzenb\"ock formulas. The eigenvalues of $K$ on $\W^{1,0}(V_\l)$ agree
 with the conformal weights for the relevant antiholomorphic weights for
 $T^{0,1}=V_{-\e_n}$ with associated projections in $T^{1,0}\otimes T^{0,1}
 \subset T\otimes T^{0,1}$. The eigenvalues of $K$ on the complementary
 subspace $\W^{0,1}(V_\l)$ correspond similarly to the conformal weights for
 the relevant holomorphic weights for $T^{1,0}=V_{+\e_1}$. Using the decision
 criterion of Lemma \ref{relw} we find the relevant antiholomorphic weights
 $+\e_1$ and $+\e_n$ for $\l=-\e_n$ and the relevant holomorphic weights
 $-\e_1$ and $-\e_n$ for $\l=+\e_1$ with corresponding decompositions
 $$
  T^{1,0}\otimes V_{-\e_n}
  \;\;=\;\;
  (\su_n\otimes_\R\C)\;\oplus\;\C
  \qquad\qquad
  T^{0,1}\otimes V_{+\e_1}
  \;\;=\;\;
  \C\;\oplus\;(\su_n\otimes_\R\C)
 $$
 where $\su_n\otimes_\R\C\,=\,V_{+\e_1-\e_n}$ is the adjoint representation.
 Using the explicit values (\ref{kaehler}) for the conformal weights we
 conclude that $K$ has eigenvalues $-n$ and $0$ with multiplicities $1$
 and $n-1$ respectively on $\W^{1,0}(V_\l)$, in other words $K$ vanishes
 on the orthogonal complement of the one--dimensional eigenspace of eigenvalue
 $-n$ in $\W^{1,0}(V_\l)$. The distinguished projections $\pr^{0,1}\in
 \W^{0,1}(V_\l)$ and $\pr^{1,0}\in\W^{0,1}(V_\l)$ are natural candidates
 for eigenvectors of $K$ for the eigenvalue $-n$ and in fact we find in
 the interpretation $\W^{1,0}(V_\l)=\Hom_{\u_n}(T^{1,0}\otimes T^{0,1}
 \otimes V_\l,\,V_\l)$
 \begin{eqnarray*}
  (\,K\pr^{0,1}\,)(a\otimes b\otimes v) &=&
  -\sum_\alpha\<X_\alpha a,\pr^{0,1}(X_\alpha b)>\,v \;\;=\;\;
  \sum_\alpha\< X_\alpha\,X_\alpha a,\pr^{0,1}b>\,v \\
  &=& \Cas^{\L^2}_{T^{1,0}}\;\pr^{0,1}(a\otimes b\otimes v)
  \quad\,\;\;=\;\;-\,n\;\pr^{0,1}(a\otimes b\otimes v)\,,
 \end{eqnarray*}
 where the $X_\alpha$ are an orthonormal basis for $\u_n$ with respect to
 the scalar product induced from $\u_n\subset\L^2T$. Turning to antiholomorphic
 Weitzenb\"ock formulas in $\W^{0,1}(V_\l)$ a similar argument implies that
 the classifying endomorphism $K$ vanishes on the orthogonal complement of
 the eigenvector $\pr^{1,0}$ with eigenvalue $-n$. Put differently this
 result reads:

 \begin{Lemma}[Explicit Form of the Classifying Endomorphism]
 \hfill\label{cle}\break
  On the basis vectors $\{\pr_\e\}$ of the space $\W(V_\l)$ of Weitzenb\"ock
  formulas on an irreducible representation $V_\l$ the classifying
  endomorphism $K$ of the holonomy algebra $\u_n$ acts by
  $$
   K\,\pr_\e\;\;=\;\;-\,\frac{\dim\,V_{\l+\e}}{\dim\,V_\l}\,\pr^{0,1}
   \qquad\textrm{or}\qquad
   K\,\pr_\e\;\;=\;\;-\,\frac{\dim\,V_{\l+\e}}{\dim\,V_\l}\,\pr^{1,0}
  $$
  depending on whether $\e$ is a holomorphic or antiholomorphic weight
  of $T$ respectively.
 \end{Lemma}

 \proof
 Taking the stated formula as a definition of a linear map $K:\W(V_\l)
 \longrightarrow\W(V_\l)$ we apply it to a holomorphic Weitzenb\"ock
 formula $F^{1,0}\,=\,\sum_{\e\subset\l^{1,0}}f_\e\pr_\e$ and find
 $$
  K(\;\sum_{\e\subset\l^{1,0}}f_\e\,\pr_\e\;)
  \;\;=\;\;
  -\,\Big(\sum_{\e\subset\l^{1,0}}f_\e\,
  \frac{\dim\,V_{\l+\e}}{\dim V_\l}\Big)\,\pr^{0,1}
  \;\;=\;\;
  -\,\<\pr^{0,1},F^{1,0}>\,\pr^{0,1}
 $$
 so that $K\pr^{0,1}=-n\,\pr^{0,1}$ and $KF^{1,0}=0$ for $F^{1,0}\in\W^{1,0}
 (V_\l)$ orthogonal to $\pr^{0,1}$. Hence the endomorphism $K$ implicitly
 defined above agrees with the classifying endomorphism on the subspace
 $\W^{1,0}(V_\l)$. Mutatis mutandis we check agreement of $K$ with the
 classifying endomorphism on the subspace $\W^{0,1}(V_\l)$.
 \qed

 \pfill
 In the case of K\"ahler holonomies the twist $\tau$ is not a polynomial
 on the classifying endomorphism $K$ and the explicit description of $K$
 in Lemma \ref{cle} does not seem to be of any use in understanding $\tau$.
 It turns out however that the Recursion Formula \ref{rec} considered as
 a matrix equation for the twist $\tau$ already has a unique solution.
 More precisely the Recursion Formula \ref{rec} can be put into the form
 \begin{equation}\label{recz}
  \tau\,K\;\;=\;\;-\,(\;\tau\,B\;+\;B\,\tau\;+\;n\,\tau\;)
 \end{equation}
 because $\tau$ is an involution and $-2\frac{\dim\,\u_n}{\dim\,T}=-n$.
 Moreover $\tau$ interchanges the subspaces $\W^{1,0}(V_\l)$ and $\W^{0,1}
 (V_\l)$ of holomorphic and antiholomorphic Weitzenb\"ock formulas, thus
 $$
  \tau\,\pr_{-\e}\;\;=\;\;\sum_{+\tilde\e}\tau_{-\e,+\tilde\e}\,\pr_{+\tilde\e}
  \qquad\qquad
  \tau\,\pr_{+\e}\;\;=\;\;\sum_{-\tilde\e}\tau_{+\e,-\tilde\e}\,\pr_{-\tilde\e}
 $$
 with as yet unknown rational matrix coefficients $\tau_{+\e,-\tilde\e}$ and
 $\tau_{-\e,+\tilde\e}$. Here and in the following the notation $-\e$ refers
 exclusively to holomorphic, $+\e$ to antiholomorphic weights. Besides
 interchanging $\W^{1,0}(V_\l)$ and $\W^{0,1}(V_\l)$ the twist interchanges
 $\pr^{0,1}$ and $\pr^{1,0}$:
 $$
  (\,\tau\,\pr^{0,1}\,)(a\otimes b\otimes v)
  \;\;=\;\;\<b,\pr^{0,1}a>\,v\;\;=\;\;\<a,\pr^{1,0}b>\,v
  \;\;=\;\;\pr^{1,0}(a\otimes b\otimes v)\,.
 $$
 In consequence Lemma \ref{cle} tells us the matrix coefficients of the
 operator $\tau\,K$
 $$
  (\,\tau\,K\,)\,\pr_{-\e}
  \;\;=\;\;-\,\frac{\dim\,V_{\l-\e}}{\dim\,V_\l}\,\tau\,\pr^{0,1}
  \;\;=\;\;-\,\sum_{+\tilde\e}\frac{\dim\,V_{\l-\e}}{\dim\,V_\l}\,
  \pr_{+\tilde\e}
 $$
 while $B$ is diagonal in the basis $\pr_\e$ with eigenvalues $b_\e$. All
 in all equation (\ref{recz}) becomes:
 $$
  (\,\tau\,K\,)_{-\e,+\tilde\e}
  \;\;=\;\;
  -\,\frac{\dim\,V_{\l-\e}}{\dim\,V_\l}
  \;\;=\;\;
  -\,(\;b_{-\e}\;+\;b_{+\tilde\e}\;+\;n\;)\,\tau_{-\e,+\tilde\e}
 $$
 This equation has a unique solution $\tau_{-\e,+\tilde\e}$, because $b_{-\e}+
 b_{+\tilde\e}+n$ is never zero. In fact
 $$
  b_{-\e_k}\;+\;\b_{+\e_l}\;+\;n
  \;\;=\;\;
  (\,\mu_k\,-\,\mu_l\,)\;+\;(\,l\,-\,k\,)\;+\;1
  \qquad\quad
  1\;\leq\;k,\,l\;\leq\;n
 $$
 equals $1$ for $k=l$ and is positive for $k<l$ as the coefficients $\mu_1
 \geq\ldots\geq\mu_n$ are decreasing. Similarly $(\mu_k-\mu_l)+(l-k)+1$ is
 negative for $k>l$ unless $k=l+1$ and $\mu_l=\mu_k$, in which case Lemma
 \ref{relw} ensures that $-\e_l$ is irrelevant by $\l_{l+1}=0$. Mutatis
 mutandis the same argument applies to the remaining matrix coefficients
 $\tau_{+\e,-\tilde\e}$ of the twist operator.

 \begin{Theorem}[Matrix Coefficients of the Twist Operator]
 \hfill\label{tm}\break
  The matrix coefficients $\tau_{-\e,+\tilde\e}$ and $\tau_{+\e,-\tilde\e}$
  of the twist operator in the basis $\pr_\e$ of projections in $\W(V_\l)$
  defined by $\tau\pr_{-\e}=\sum_{+\tilde\e}\tau_{-\e,+\tilde\e}\,
  \pr_{+\tilde\e}$ and $\tau\pr_{+\e}=\sum_{-\tilde\e}\tau_{+\e,-\tilde\e}
  \,\pr_{-\tilde\e}$ are given by:
  $$
   \tau_{-\e,+\tilde\e}
   \;\;=\;\;
   \frac1{b_{-\e}+b_{+\tilde\e}+n}\,\frac{\dim\,V_{\l-\e}}{\dim\,V_\l}
   \qquad
   \tau_{+\e,-\tilde\e}
   \;\;=\;\;
   \frac1{b_{+\e}+b_{-\tilde\e}+n}\,\frac{\dim\,V_{\l+\e}}{\dim\,V_\l}
  $$
 \end{Theorem}

 \noindent
 Another interesting consequence of the decomposition of the space $\W(V_\l)$
 into the subspaces $\W^{1,0}(V_\l)$ and $\W^{0,1}(V_\l)$ of holomorphic
 and antiholomorphic Weitzenb\"ock formulas should not pass unnoticed. As
 the twist $\tau$ interchanges $\W^{1,0}(V_\l)$ and $\W^{0,1}(V_\l)$
 every eigenvector for $\tau$ is determined by its summand in the subspace
 $\W^{1,0}(V_\l)$. In fact an arbitrarily chosen $F^{1,0}\in\W^{1,0}
 (V_\l)$ gives rise to an eigenvector
 $$
  F\;\;:=\;\;F^{1,0}\;\pm\;\tau(\,F^{1,0}\,)
 $$
 for $\tau$ of eigenvalue $\pm1$. Of course we are primarily interested in
 the eigenvalue $-1$:

 \begin{Corollary}[Pure Curvature Terms in K\"ahler Holonomy]
 \hfill\label{purec}\break
  Recall that the eigenvectors $F=\sum_\e f_\e\pr_\e$ of the twist operator
  $\tau$ for the eigenvalue $-1$ in the space $\W(V_\l)$ of Weitzenb\"ock
  formulas on an irreducible representation $V_\l$ correspond to Weitzenb\"ock
  formulas reducing to a pure curvature term, compare (\ref{nabla}) and
  (\ref{cterm}):
  $$
   F(\;\nabla^2\;)
   \;\;=\;\;
   -\,\sum_\e\,T^*_\e\,T^{\hphantom{*}}_\e
   \;\;=\;\;
   \frac12\,\sum_{\mu,\nu}F_{t_\mu\otimes t_\nu}\,R_{t_\mu,\,t_\nu}\ .
  $$
  In K\"ahler holonomies there is a bijection $F^{1,0}\longmapsto
  F^{1,0}-\tau(F^{1,0})$ between $\W^{1,0}(V_\l)$ and the eigenspace
  of $\tau$ for eigenvalue $-1$ or pure curvature term Weitzenb\"ock
  formulas.
 \end{Corollary}

 \noindent
 On a K\"ahler manifold $M$ the holonomy algebra bundle $\u M$ can be
 identified with the bundle of skew--symmetric endomorphisms of the
 tangent bundle commuting with the parallel complex structure $I$, it
 has a parallel subbundle $i\R M\subset\u M$ associated to the center
 $i\R\subset\u_n$ of the holonomy algebra and spanned by $I$ in every
 point. It is well--known that the curvature tensor $R$ of a K\"ahler
 manifold is a $(1,1)$--form in the sense $R_{IX,IY}Z\,=\,R_{X,Y}Z$ for
 all $X,Y,Z\in TM$, hence the Ricci endomorphism is not only a symmetric
 endomorphism of the tangent bundle, but commutes with $I$, too:
 $$
  \Ric(IZ)
  \;\;=\;\;\sum_\mu R_{IZ,t_\mu}t_\mu
  \;\;=\;\;-\sum_\mu R_{Z,It_\mu}t_\mu
  \;\;=\;\;\sum_\mu R_{Z,t_\mu}It_\mu
  \;\;=\;\;I\,(\Ric\,Z)\,.
 $$
 In turn the skew--symmetric endomorphism $I\Ric$ commuting with $I$ can
 be thought of as a section of the holonomy algebra bundle $\u M$ with a
 natural action on every vector bundle $VM$ associated to the holonomy
 reduction. A different way to understand this action is to note that
 $I\Ric\in\Gamma(\u M)$ is up to sign the result of applying the curvature
 operator $R:\,\L^2TM\longrightarrow\u M$ to the complex structure, because
 we find the action
 $$
  R(\,I\,)\,Z
  \;\;=\;\;
  \frac12\sum_\mu R_{t_\mu,It_\mu}Z
  \;\;=\;\;
  -\frac12\sum_\mu\Big(\,R_{Z,t_\mu}It_\mu\,-\,R_{Z,It_\mu}t_\mu\,\Big)
  \;\;=\;\;
  -I\,\Ric\,Z
 $$
 on $TM$. Using $\<I,I>=-\frac12I^2=n$ we conclude that the curvature
 term $q^{i\R}(R)$ reads:
 \begin{equation}\label{ccop}
  q^{i\R}(\,R\,)
  \;\;=\;\;\frac1n\,I\,R(\,I\,)
  \;\;=\;\;-\,\frac1n\,I\,(I\,\Ric)
 \end{equation}
\section{A Matrix Presentation of the Twist Operator $\tau$}
\label{matrix}
 In the last section we have studied the Weitzenb\"ock machine for
 K\"ahler holonomies and succeeded in describing the twist $\tau$
 and the classifying endomorphism $K$ explicitly in the natural basis
 $\pr_\e$ of the space of Weitzenb\"ock formulas. The principal idea
 of this calculation was to read the Recursion Formula \ref{rec},
 which we used in Section \ref{eigenvalues} for a recursive construction
 of an eigenbasis for the classifying endomorphism, directly as a matrix
 equation for the unknown matrix of $\tau$. It turns out that the same
 idea allows us to derive a closed matrix expression for the twist operator
 $\tau$ and the classifying endomorphism $K$ in the three other holonomies
 $\SO(n)$, $\G_2$ and $\Spin(7)$ considered in this article.

 \pfill
 In the case of $\SO(n)$--holonomy the classifying endomorphism $K$
 and the twist differ only on the span of the unit $\1\in\W(V_\l)$,
 which becomes the identity $\1=\id_{T\otimes V_\l}$ in the interpretation
 $\W(V_\l)=\End_{\so_n}(T\otimes V_\l)$ and the connection Laplacian
 $\1(\nabla^2)=-\nabla^*\nabla$ in the interpretation as a second order
 parallel differential operator. The norm of $\1$ is simply the dimension
 $\<\1,\1>=n$ of $T$ and so we conclude
 \begin{equation}\label{kqmod}
  K\;\;=\;\;\tau\;-\;\1\otimes\1
 \end{equation}
 where $\1\otimes\1$ denotes the endomorphism:
 $$
  \1\otimes\1:\qquad\W(V_\l)\;\longmapsto\;\W(V_\l),
  \qquad F\;\longmapsto\;\<\1,F>\,\1\,.
 $$
 In fact in the interpretation $\W(V_\l)=\Hom_{\so_n}(T\otimes T
 \otimes V_\l,\,V_\l)$ both $K$ and $\tau$ act by precomposition
 with an endomorphism of $T\otimes T$ extended to $T\otimes T\otimes
 V_\l$ in such a way that the eigenvalues of $K$ and $\tau$ on the
 orthogonal complement $\S^2_0T\oplus\L^2$ of $\C$ agree according
 to Table \ref{kt}. Hence we need only check equation (\ref{kqmod})
 on $\1$, where we can use the trivial statement $\tau(\1)=\1$ and
 $K(\1)=-(n-1)\1$ from Lemma \ref{elem}. In light of equation
 (\ref{kqmod}) the Recursion Formula \ref{rec} can be rewritten
 in the form
 \begin{equation}\label{recso}
  \tau\,K
  \;\;=\;\;\id\;-\;\1\,\otimes\,\1
  \;\;=\;\;-\,(\;\tau\,B\;+\;B\,\tau\;+\;(n-1)\,\tau\;)
 \end{equation}
 in the case of $\SO(n)$--holonomy, note that $\tau$ is an involution
 with $\tau\circ\1\otimes\1=\1\otimes\1$. In order to turn equation
 (\ref{recso}) into a matrix equation for the unknown matrix coefficients
 $$
  \tau\,\pr_\e\;\;=:\;\;\sum_{\tilde\e}\tau_{\e,\tilde\e}\,\pr_{\tilde\e}
  \qquad\qquad
  \tau_{\e,\tilde\e}\;\in\;\mathbb{Q}
 $$
 of the twist $\tau$ we observe that by its very definition $(\1\otimes\1)
 \pr_\e\,=\,\<\1,\pr_\e>\1\,=\,\frac{\dim\,V_{\l+\e}}{\dim\,V_\l}\1$.
 With this additional piece of information equation (\ref{recso}) becomes
 $$
  \delta_{\e,\tilde\e}\;-\;\frac{\dim\,V_{\l+\e}}{\dim\,V_\l}
  \;\;=\;\;
  -\,(\;b_\e\;+\;b_{\tilde\e}\;+\;n\;-\;1\;)\,\tau_{\e,\tilde\e}
 $$
 as again the operator $B$ is diagonal in the basis $\pr_\e$. Unlike the
 K\"ahler case we can not exclude the possibility that $b_\e+b_{\tilde\e}
 +n-1$ is zero. If we assume $\e$ to be relevant and thus $\dim\,V_{\l+\e}
 \neq 0$ however, then this can only happen on the diagonal $\e=\tilde\e$
 under the assumption $\dim\,V_{\l+\e}\,=\,\dim\,V_\l$. For all relevant,
 off--diagonal weights $\e\neq\tilde\e$ the matrix coefficients of $\tau$
 are thus given by
 \begin{equation}\label{som}
  \tau_{\e,\tilde\e}
  \;\;=\;\;
  -\,\frac1{b_\e\,+\,b_{\tilde\e}\,+\,n\,-\,1}\,
  \Big(\;\delta_{\e,\tilde\e}\;-\;\frac{\dim\,V_{\l+\e}}{\dim\,V_\l}\;\Big)\ .
 \end{equation}
 and this formula is still true on the diagonal $\e=\tilde\e$ provided
 $\dim\,V_{\l+\e}\neq\dim\,V_\l$. In order to eliminate these rather
 unpleasant restrictions we note that the equation $\tau(\1)=\1$ allows
 us to recover the diagonal coefficients $\tau_{\e,\e}$ of $\tau$ with
 $\dim\,V_{\l+\e}=\dim\,V_\l$. Moreover the matrix coefficients $\tau_{\e,
 \tilde\e}$ are rational functions in the highest weight $\l$, because
 the conformal weights $b_\e$ and $b_{\tilde\e}$ are linear in $\l$ and
 the dimensions of $V_\l$ and $V_{\l+\e}$ are polynomials in $\l$ due to
 Weyl's Dimension Formula. Taken together these arguments imply that the
 troublesome matrix coefficients $\tau_{\e,\tilde\e}$ left undefined by
 equation (\ref{som}) agree with the values of the analytic continuation
 considered as a rational function of $\l$ to values to the discriminant
 set, where zeroes of the numerator and the denominator cancel out.

 \begin{Theorem}[The Twist Operator in $\SO(n)$--Holonomy]
 \hfill\label{matw}\break
  The matrix coefficients $\tau_{\e,\tilde\e}$ of the twist $\tau:\,
  \W(V_\l)\longrightarrow\W(V_\l)$ in the basis of projections $\pr_\e
  \in\W(V_\l)$ are defined by $\tau\,\pr_\e=\sum_{\tilde\e}\tau_{\e,\tilde\e}
  \,\pr_{\tilde\e}$. In $\SO(n)$--holonomy they are given by
  $$
   \tau_{\e,\tilde\e}
   \;\;=\;\;
   -\,\frac1{b_\e\,+\,b_{\tilde\e}\,+\,n\,-\,1}\,
   \Big(\;\delta_{\e,\tilde\e}\;-\;\frac{\dim\,V_{\l+\e}}{\dim\,V_\l}\;\Big)
  $$
  provided we interprete $\tau_{\e,\tilde\e}$ as a rational function in
  the highest weight $\l$ and cancel zeroes of the numerator and denominator
  for values of $\e$ and $\tilde\e$, for which $b_\e+b_{\tilde\e}+n-1=0$.
 \end{Theorem}

 \pfill
 Turning from the case of $\SO(n)$--holonomy to $\G_2$ we note that the
 relation between the twist operator $\tau$ and the classifying endomorphism
 $K$ becomes more complicated than equation (\ref{kqmod}) due to the presence
 of the Bochner identity (\ref{g2bx}). For definiteness we let $\beta:=
 F_{\mathrm{Bochner}}$ be the $K$--eigenvector for the eigenvalue $-2$
 specified in (\ref{g2bx}), whose coefficients $\beta=\sum_\e\beta_\e\,\pr_\e$
 are given explicitly as polynomials in the highest weight $\l=a\o_1+b\o_2,\,
 a,\,b\geq 0$. The argument in the $\SO(n)$--case extends to the relation
 $$
  K\;\;=\;\;
  \frac13\,\tau\;+\;\frac13\;-\;\frac23\,\1\otimes\1
  \;-\;\frac2{\<\beta,\beta>}\beta\otimes\beta
 $$
 between $K$ and $\tau$ in the $\G_2$--case so that the Recursion Formula
 \ref{rec} becomes
 $$
  \tau\,K
  \;\;=\;\;
  \frac13\,\id\;+\;\frac13\,\tau\;-\;\frac23\,\1\otimes\1
  \;+\;\frac2{\<\beta,\beta>}\beta\otimes\beta
  \;\;=\;\;
  -\,(\;\tau\,B\;+\;B\,\tau\;+\;4\,\tau\;)\ ,
 $$
 because $\beta$ is an eigenvector for $\tau$ of eigenvalue $-1$. Moreover
 the calculation
 $$
  (\,\beta\otimes\beta\,)\,\pr_\e
  \;\;=\;\;\sum_{\tilde\e}
  \Big(\,\beta_\e\,\frac{\dim\,V_{\l+\e}}{\dim\,V_\l}\,\beta_{\tilde\e}\,\Big)
  \pr_{\tilde\e}
 $$
 tells us the matrix coefficients of $\beta\otimes\beta$ and so we end up
 with the formula
 \begin{equation}\label{g2om}
  \tau_{\e,\tilde\e}
  \;\;=\;\;
  -\,\frac1{3\,b_\e\;+\;3\,b_{\tilde\e}\;+\;13}\,
  \left(\;\delta_{\e\tilde\e}\;-\;2\,\Big(\;1\,-\,3\,
  \frac{\beta_\e\,\beta_{\tilde\e}}{\<\beta,\beta>}\;\Big)
  \,\frac{\dim\,V_{\l+\e}}{\dim\,V_\l}\;\right)
 \end{equation}
 for the matrix coefficients $\tau_{\e,\tilde\e}$ of $\tau$ defined by
 $\tau\,\pr_\e=\sum_{\tilde\e}\tau_{\e,\tilde\e}\pr_{\tilde\e}$. Again
 this formula has to be read as a rational function in the highest weight
 $\l$ to get the right values for $\tau_{\e,\tilde\e}$ for index tuples
 $\e,\,\tilde\e$ satisfying $3b_\e+3b_{\tilde\e}+13=0$ by cancelling the
 zero of the denominator with the apparent zero of the numerator.

 \begin{Theorem}[The Twist Operator in $\G_2$--Holonomy]
 \hfill\label{matg2}\break
  The matrix coefficients $\tau_{\e,\tilde\e}$ of the twist $\tau:\,
  \W(V_\l)\longrightarrow\W(V_\l)$ in the basis of projections $\pr_\e
  \in\W(V_\l)$ are defined by $\tau\,\pr_\e=\sum_{\tilde\e}\tau_{\e,\tilde\e}
  \,\pr_{\tilde\e}$. In $\G_2$--holonomy they are given by
  $$
   \tau_{\e,\tilde\e}
   \;\;=\;\;
   -\,\frac1{3\,b_\e\;+\;3\,b_{\tilde\e}\;+\;13}\,
   \left(\;\delta_{\e\tilde\e}\;-\;2\,\Big(\;1\,-\,3\,
   \frac{\beta_\e\,\beta_{\tilde\e}}{\<\beta,\beta>}\;\Big)
   \,\frac{\dim\,V_{\l+\e}}{\dim\,V_\l}\;\right)
  $$
  where the coefficients $\beta_\e$ of the Bochner identity $\beta=\sum_\e
  \beta_\e\pr_\e$ are given by (\ref{g2bx}).
 \end{Theorem}

 \pfill
 Eventually we want to sketch briefly the corresponding argument in the
 case of exceptional holonomy $\Spin(7)$. Chosing the Bochner identity
 $\beta:=F_{\mathrm{Bochner}}$ of equation (\ref{s7b}) with coefficients
 $\beta=\sum_\e\beta_\e\,\pr_\e$ we can cast the relation between
 $K$ and $\tau$ into the form
 $$
  K\;\;=\;\;
  \frac12\,\tau\;+\;\frac14\;-\;\frac34\,\1\otimes\1
  \;-\;\frac2{\<\beta,\beta>}\beta\otimes\beta
 $$
 so that the Recursion Formula \ref{rec} is equivalent to the identity:
 $$
  \tau\,K
  \;\;=\;\;
  \frac12\,\id\;+\;\frac14\,\tau\;-\;\frac34\,\1\otimes\1
  \;+\;\frac2{\<\beta,\beta>}\beta\otimes\beta
  \;\;=\;\;
  -\,(\;\tau\,B\;+\;B\,\tau\;+\;\frac{21}4\,\tau\;)\ ,
 $$
 All in all we end up with the following formula for the matrix coefficients
 \begin{equation}\label{s7om}
  \tau_{\e,\tilde\e}\;\;=\;\;
  -\,\frac1{2\,b_\e\;+\;2\,b_{\tilde\e}\;+\;11}\,
  \left(\;\delta_{\e\tilde\e}\;-\;\Big(\;\frac32\,-\,4\,
  \frac{\beta_\e\,\beta_{\tilde\e}}{\<\beta,\beta>}\;\Big)
  \,\frac{\dim\,V_{\l+\e}}{\dim\,V_\l}\;\right)
 \end{equation}
 of the twist operator $\tau$ defined by $\tau\,\pr_\e=\sum_{\tilde\e}
 \tau_{\e,\tilde\e}\,\pr_{\tilde\e}$. Recall that the twist $\tau$
 classifies all Weitzenb\"ock formulas $F\in\W(V_\l)$, which reduce
 to a pure curvature term. In consequence the matrix expressions for
 the the twist $\tau$ in the basis $\{\pr_\e\}$ of $\W(V_\l)$ allows
 us to check this condition effectively for every given Weitzenb\"ock
 formula on all K\"ahler manifolds $M$ and on all manifolds $M$ with
 holonomy $\SO(n)$, $\G_2$ or $\Spin(7)$.
\section{Examples}
\label{examples}
 In this section we will present a few examples of how to obtain for
 a given representation $V_\l$ the possible Weitzenb\"ock formulas on
 sections of the associated bundle $V_\l M$. The general procedure is
 as follows: we first determine the relevant weights $\e$ using the
 diagrams of Section \ref{tables}. This gives the decomposition of $T
 \otimes V_\l$ into irreducible summands and defines the twistor
 operators $T_\e$. Next we compute the $B$--eigenvalues $b_\e $, e.g.
 using the general formula of Corollary~\ref{confW} and obtain the
 universal Weitzenb\"ock formulas of Proposition \ref{universal}.
 Other
 Weitzenb\"ock formulas correspond to the $B$--polynomials constructed
 in the preceding section. If $F= F(B)$ is such a polynomial then the
 coefficient of $T^*_\e T^{\phantom{*}}_\e$ is given as $-F(b_\e)$.

 \pfill
 As a first example we consider the bundle of $p$--forms on a
 Riemannian manifold $(M^n,g)$ for simplicity we assume $n=2r+1$ and
 $p\le r-1$, i.e. $\g = \so_{2r+1}$ and $\l= \o_p$, i.e. the
 associated bundle is the bundle of $p$--forms. The relevant weights
 according to the tables of Section~\ref{eigenvalues} are
 $\,\e_1,\,-\e_p$ and $\e_{p+1}$ with the decomposition
 $$
  T \otimes V_\l
  \;\;=\;\;
  V_{\l+\e_1} \,\oplus\, V_{\l -\e_p}
  \,\oplus\, V_{\l + \e_{p+1}}
  \;\;\cong\;\;
  V_{\l+\e_1} \,\oplus\, \L^{p-1} \,\oplus\, \L^{p+1} \
 $$
 and twistor operators $T_{\e_1},\,T_{-\e_p}$ and $T_{\e_p}$. To
 compare our twistor operators with differential and codifferential
 $\,d, d^*\,$ we have to embed $\Lambda^{p-1}$ resp. $\Lambda^{p+1}$
 into the tensor product $T \otimes \Lambda^p$. This leads to the
 following formula
 $$
  T_{-\e_p}^*\,T_{-\e_p}^{\phantom{*}}
  \;\;=\;\;
  \tfrac{1}{n-p+1}\;d\,d^*,
  \qquad
  T_{+\e_p}^*T_{+\e_p}^{\phantom{*}}
  \;\;=\;\;
  \tfrac{1}{p+1}\;d^*\,d
 $$
 Next we take the relevant $B$--eigenvalues from
 Section~\ref{eigenvalues} they are:
 $$
  b_{+\e_1} = 1,\qquad b_{-\e_p} = -n + p, \qquad b_{e_{p+1}} = -p \ .
 $$
 Since we have only three summands in the decomposition of $T\otimes
 V_\l$ we obtain only one Weitzenb\"ock formula with a pure curvature
 term, which is the formula given in Proposition \ref{universal}:
 $$
  \begin{array}{rcccccccc}
   q(R)\quad=
   &\!\!-\!\!&T_{+\e_1}^*\,T_{+\e_1}^{\phantom{*}}
   &\!\!+\!\!&(n-p)&\!\!T_{-\e_p}^*\,T_{-\e_p}^{\phantom{*}}
   &\!\!+\!\!&  p  &\!\!T_{+\e_p}^*\,T_{+\e_p}^{\phantom{*}}
   \\[.5cm]
   =
   &\!\!-\!\!&T_{+\e_1}^*\,T_{+\e_1}^{\phantom{*}}
   &\!\!+\!\!&\frac{n-p}{n-p+1}&d\,d^*\quad
   &\!\!+\!\!&\frac{p}{p+1}    &d^*\,d\quad
  \end{array}
 $$
 If we add the Weitzenb\"ock formula~(\ref{nabla}) for $\nabla^*\nabla$ to
 this expression for $q(R)$ we obtain the classical Weitzenb\"ock formula
 for the Laplacian on $p$--forms:
 $$
  \Delta
  \;\;=\;\;
  \nabla^*\nabla\;+\;q(R)
  \;\;=\;\;
  (n-p+1)\,T_{-\e_p}^*\,T_{-\e_p}^{\phantom{*}}
  \;+\;
  (p+1)\,T_{+\e_p}^*\,T_{+\e_p}^{\phantom{*}}
  \;\;=\;\;
  d\,d^* + d^*\,d
 $$

 \pfill
 Let $(M^{2r},g)$ be a Riemannian spin manifold with spinor bundle
 $S = S_+ \oplus S_-$. We consider the two bundles defined by the
 Cartan summand in $S_\pm\otimes T$ with highest weights $\l_+ =
 \o_1 + \o_{r-1}$ and $\l_- = \o_1 + \o_r$. Using the tables of
 Section~\ref{tables} we find the relevant weights $+\e_1, \,-\e_1$
 and $+\e_2$ for both $\l_\pm$ and in addition $-\e_r$ or $+\e_r$
 for $\l_+$ or $\l_-$ respectively. The corresponding tensor product
 decomposition is
 $$
  T \otimes  V_{\l_\pm}
  \;\;=\;\;
  V_{\l_\pm +\e_1} \,\oplus\,V_{\l_\pm-\e_1}\,\oplus\,V_{\l_\pm + \e_2}
  \,\oplus\,V_{\l_\pm \mp  \e_r}\ .
 $$
 Note that $\l_\pm - \e_1$ is the defining representation for the
 bundles $S_\pm$ and that $\l_\pm \mp \e_r = \l_\mp$. Projecting the
 covariant derivative of a section of $T \otimes V_{\l_\pm}$ onto one
 of these summands defines four twistor operators.  The fourth operator
 $T_{\mp \e_r}: \Gamma(V_{\l_\pm})\longrightarrow\Gamma(V_{\l_\mp})$ is
 usually called the {\it Rarita--Schwinger operator}. A solution of the
 Rarita--Schwinger equation is by definition a section of $\psi\in\Gamma
 (V_{\l_\pm})$  with both $T_{\mp\e_r}\psi = 0$ and $T_{-\e_1}\psi = 0$.

 The $B$--eigenvalues for $\so_n$--representations
 were computed in Section~\ref{eigenvalues} and in particular:
 $$
  b_{+\e_1} \;=\; \tfrac32,\quad b_{-\e_1} \;=\; -2r + \tfrac12,\quad
  b_{+\e_2} \;=\; -\,\tfrac12,\quad b_{\pm \e_r} \;=\; -\,r +\tfrac12
 $$
 Since the decomposition of $T \otimes V_{\l_\pm}$ has four summands
 we will obtain two Weitzenb\"ock formulas with a pure curvature
 term. The first one is again the universal Weitzenb\"ock formula
 of Proposition \ref{universal} corresponding to $B$, whereas the
 second corresponds to $p_3(B)$, the degree 3 polynomial of the recursion
 procedure defined in (\ref{gleichung2}). Its coefficients are the values
 $p_3(b_\e)$ for the relevant weights $\e$. The Casimir operator of an
 irreducible $\so_n$--representation $V_\l$ with highest weight $\l$ is
 computed as $\Cas_{V_\l} =-\<\l+ 2\rho, \, \l>$, where $\,\<\cdot,\cdot>\,$
 is the standard scalar product on $\R^r$. In particular we have
 $\Cas_{V_{\l_\pm}}=-\,\frac14\,r\,(2r+7)$. Eventually we obtain
 the following two Weitzenb\"ock formulas on sections of $\,V_{\l_\pm}$:
 $$
  \begin{array}{rcl}
   q(R)
   &=&
   -\frac32\,T_{+\e_1}^*\,T_{+\e_1}^{\phantom{*}} \;+\;
   (2r-\frac12)\,T_{-\e_1}^*\,T_{-\e_1}^{\phantom{*}} \;+\;
  \frac12\,T_{+\e_2}^*\,T_{+\e_2}^{\phantom{*}}
   \;+\;(r-\frac12)\,T_{\pm\e_r}^*\,T_{\pm\e_r}^{\phantom{*}}\\[.7cm]
   p_3(B)(\nabla^2)
   &=&
   -(\frac32+r)(r-1)\,T_{+\e_1}^*\,T_{+\e_1}^{\phantom{*}} \;\;+\;\;
   (2r-1)(r^2-1)\,T_{-\e_1}^*\,T_{-\e_1}^{\phantom{*}} \\[.4cm]
   &&
   \qquad\qquad\qquad \qquad\qquad \quad\qquad\qquad\;+\;\;
   (r-\frac12)(r+1)\,T_{+\e_2}^*\,T_{+\e_2}^{\phantom{*}}
   \;+\;T_{\pm\e_r}^{\phantom{*}}
  \end{array}
 $$
 Note that similar Weitzenb\"ock formulas were obtained in \cite{branson}.
 More precisely their curvature terms $Z_1$ and $Z_2$ are related to $B$
 and $p_3(B)$ by the following equations:
 $$
  Z_1
  \;\;=\;\;
  \tfrac{(2r+3)(r-1)}{r(2r+1)}\,B\;-\;\,\tfrac{3}{r(2r+1)}\,p_3(B)
  \qquad
  Z_2
  \;\;=\;\;
  -\,\tfrac{(2r-1)(r+1)}{r(2r+1)}\,B \;+\;\,\tfrac{1}{r(2r+1)}\,p_3(B)\ ,
 $$
 whereas the operators are related by: $T_{+\e_1}= G_Z,\;
 T_{-\e_1}= G_{\Sigma},\;T_{+\e_2}= G_Y,\; T_{\mp\e_r}=G_T$.

 \pfill
 In the last part of this section we want to describe for
 $\G_2$-- and $\Spin(7)$--holonomy all pure curvature Weitzenb\"ock
 formulas on parallel subbundles of the form bundle. In particular we
 will present the form Laplacian $\Delta = d^*d + d d^* =
 \nabla^*\nabla + q(R)$ as a linear combination of the operators $\,
 T^*_\e T^{\phantom{*}}_\e$ and discuss the existence of harmonic
 forms.

 We start with the case of $\G_2$--holonomy. Let $\Gamma_{a,\,b}$ be
 the irreducible $\G_2$--representation with highest weight $\,a\o_1+b\o_2,
 \,a,b \ge 0$, e.g. $\Gamma_{0,\,0} =\C$ is the trivial representation,
 $\Gamma_{1,\,0}=T$ and $\Gamma_{0,\,1}= \L^2_{14} \cong \g_2$. Recall
 that up to dimension 77 irreducible $\G_2$--representations are uniquely
 determined by their dimension. However there are two different
 irreducible representations in dimension $77$, one of them is $[77]^-:=
 \Gamma_{3,\,0}$, the other is $\Gamma_{0,2}$, the space of $\G_2$-curvature
 tensors. Moreover
 $$
  \dim \Gamma_{2,\,0} \;=\;27,\qquad\dim \Gamma_{1,\,1} \;=\; 64 \ .
 $$
 The spaces of 2-- and 3--forms have the following decompositions:
 \begin{equation}\label{deco4}
  \L^2T \;\cong\; \Lambda^5T \;\cong\;T \oplus
  \L^2_{14},\qquad\qquad \Lambda^3T \;\cong\; \Lambda^4T
  \;\cong\; \C \oplus T \oplus \Lambda^3_{27} \ ,
 \end{equation}
 where the subscripts denote the dimension of the representation.
 Next we give the relevant weights for the representations
 $\Gamma_{1,\,0},\,\Gamma_{0,\,1}$ and $\Gamma_{2,\,0}$. We start
 with $\l=\o_2$, i.e. the representation $V_\l = \Gamma_{0,\,1} =
 \Lambda^2_{14}$, here the relevant weights are $\e=-\e_2,\,\e_3,\,\e_1$
 with $\l+\e= \o_1,\,2\o_1,\,\o_1+\o_2$ and the corresponding decomposition
 reads
 $$
  T \,\otimes\, \Gamma_{0,\,1}
  \;\;=\;\;
  \Gamma_{1,\,0}\,\oplus\,\Gamma_{2,\,0}\,\oplus\,\Gamma_{1,\,1}
  \;\;=\;\;
  T^* \,\oplus\, \Lambda^3_{27}T^* \,\oplus\, [64] \ .
 $$
 In this case the universal Weitzenb\"ock formula of Proposition
 \ref{universal}) is the only pure curvature Weitzenb\"ock formula.
 With the explicit $B$--eigenvalues given in Section~\ref{eigenvalues}
 we find on sections of $\Lambda^2_{14}T^*M$:
 $$
  \begin{array}{ccrlrlrl}
   q(R) &=
   &&4       \,T^*_{-\e_2}\,T^{\phantom{*}}_{-\e_2}
   &+&\tfrac43\,T^*_{+\e_3}\,T^{\phantom{*}}_{+\e_3}
   &-&        \;T^*_{+\e_1}\,T^{\phantom{*}}_{+\e_1}
   \\[.5cm]
   \Delta &=
   &&5       \,T^*_{-\e_2}\,T^{\phantom{*}}_{-\e_2}
   &+&\tfrac73\,T^*_{+\e_3}\,T^{\phantom{*}}_{+\e_3}
   &&
  \end{array}
 $$
 Hence a form $\,\psi\,$  in $\,\Lambda^2_{14}\subset \Lambda^2T\,$
 is harmonic if and only if $T_{-\e_2}\psi = 0 = T_{+\e_3}\psi$ (the
 manifold is assumed to be compact), i.e. if and and only if $\,
 \nabla\psi = T_{+\e_1}\psi$ or equivalently if and only if $\nabla
 \psi $ is a section of $\,\Gamma_{11} = [64]$.  This statement
 corresponds to the fact that on a compact manifold a form is
 harmonic if and only if it is closed and coclosed.

 \pfill
 For $\lambda = \o_1$, i.e. the representation $V_\l =
 \Gamma_{1,0} = T$ the relevant weights are determined as $\e=
 -\e_1,\,0,\,+\e_2,\,+\e_1$ with $\l+\e = 0,\,\o_1,\,\o_2,\,2\o_1\,$
 leading to the decomposition:
 $$
  T \,\otimes\, \Gamma_{1,0}
  \;\;=\;\;
  \Gamma_{0,\,0}\,\oplus\,\Gamma_{1,\,0}\,\oplus\,
  \Gamma_{0,\,1}\,\oplus\,\Gamma_{2,\,0}
  \;\;=\;\;
  \C \,\oplus\, T^*\,\oplus\,\Lambda^2_{14}T^*\,\oplus\,\Lambda^3_{27}T^* \ .
 $$
 Here we have two pure curvature Weitzenb\"ock formulas. In fact both
 curvature terms are zero, since $q(R)=\Ric = 0$ on $\Gamma_{1,0} = T$.
 In addition to the universal Weitzenb\"ock formula of Proposition
 \ref{universal}) we have the equation corresponding to the polynomial
 $F_3$ given in (\ref{g2bx}) with $\,a=1, b=0$. After substituting the
 $B$--eigenvalues we obtain the following Weitzenb\"ock formulas on
 $1$--forms:
 $$
  \begin{array}{ccrlrlrlrl}
   0 &=
   &&4       \,T^*_{-\e_1}\,T^{\phantom{*}}_{-\e_1}
   &+&2       \,T^*_{0}    \,T^{\phantom{*}}_{0}
   & &
   &-&\tfrac23\,T^*_{+\e_1}\,T^{\phantom{*}}_{+\e_1}
   \\[.5cm]
   0 &=
   &-&\tfrac{16}3\,T^*_{-\e_1}\,T^{\phantom{*}}_{-\e_1}
   &+&\tfrac83   \,T^*_{0}    \,T^{\phantom{*}}_{0}
   &-&\tfrac83   \,T^*_{+\e_2}\,T^{\phantom{*}}_{+\e_2}
   &+&\tfrac89   \,T^*_{+\e_1}\,T^{\phantom{*}}_{+\e_1}
   \\[.5cm]
   \Delta &=
   &&5       \,T^*_{-\e_1}\,T^{\phantom{*}}_{-\e_1}
   &+&3       \,T^*_{0}    \,T^{\phantom{*}}_{0}
   &+&        \;T^*_{+\e_2}\,T^{\phantom{*}}_{+\e_2}
   &+&\tfrac13\,T^*_{+\e_1}\,T^{\phantom{*}}_{+\e_1}
  \end{array}
 $$
 It follows that $\,\D\ge \frac 13 \, \nabla^*\nabla$, i.~e.~
 there are no non--parallel harmonic 1--forms, which is of course
 Bochners theorem in the case of  $\G_2$--manifolds.

 \pfill
 Next we consider the case $\l=2\o_1$, i.~e.~the
 representation $\,V_\l = \Gamma_{2,0} = \Lambda^3_{27}$. Here the
 relevant weights are $\,\e = -\e_1,\,\,0,\,-\e_3,\,\e_2,\,\e_1\,$
 with $\,\l+\e = \o_1,\,2\o_1,\,\o_2,\,\o_1+\o_2,\,3\o_1\,$ and with
 the decomposition
 $$
  \begin{array}{rccccccccc}
   T \otimes \Gamma_{2,\,0}
   \quad=\;
   &\Gamma_{1,\,0} &\!\!\!\oplus\!\!\!
   &\Gamma_{2,\,0} &\!\!\!\oplus\!\!\!
   &\Gamma_{0,\,1} &\!\!\!\oplus\!\!\!
   &\Gamma_{1,\,1} &\!\!\!\oplus\!\!\!
   &\Gamma_{3,\,0}
   \\[10pt]
   =\;
   &T^*          &\!\!\!\oplus\!\!\!
   &\L^3_{27}T^* &\!\!\!\oplus\!\!\!
   &\L^2_{14}    &\!\!\!\oplus\!\!\!
   &[\,64\,]     &\!\!\!\oplus\!\!\!
   &[\,77\,]^-
  \end{array}
 $$
 Hence we have two pure curvature Weitzenb\"ock formulas on sections
 of $\L^3_{27}$. The first one is the formula for $q(R)$
 corresponding to $B$, while the second corresponds to
 $\,-\frac{27}{240}\,p_3(B)$:
 $$
  \begin{array}{ccrlrlrlrlrl}
   q(R) &=
   &\!\!&\tfrac{14}{3}\,T^*_{-\e_1}\,T^{\phantom{*}}_{-\e_1}
   &+\!\!&2            \,T^*_{0}    \,T^{\phantom{*}}_{0}
   &+\!\!&\tfrac83     \,T^*_{-\e_3}\,T^{\phantom{*}}_{-\e_3}
   &-\!\!&\tfrac{1}{3} \,T^*_{+\e_2}\,T^{\phantom{*}}_{+\e_2}
   &-\!\!&\tfrac{4}{3} \,T^*_{+\e_1}\,T^{\phantom{*}}_{+\e_1}
   \\[8pt]
   0 &=
   &-\!\!&\;\tfrac{7}{6}\,T^*_{-\e_1}\,T^{\phantom{*}}_{-\e_1}
   &+\!\!&\tfrac12     \,T^*_{0}    \,T^{\phantom{*}}_{0}
   &+\!\!&\tfrac56     \,T^*_{-\e_3}\,T^{\phantom{*}}_{-\e_3}
   &-\!\!&\tfrac{2}{3} \,T^*_{+\e_2}\,T^{\phantom{*}}_{+\e_2}
   &+\!\!&\tfrac{1}{3} \,T^*_{+\e_1}\,T^{\phantom{*}}_{+\e_1}
   \\[8pt]
   \Delta &=
   &\!\!&\;\tfrac92   \,T^*_{-\e_1}\,T^{\phantom{*}}_{-\e_1}
   &+\!\!&\tfrac72     \,T^*_{0}    \,T^{\phantom{*}}_{0}
   &+\!\!&\tfrac92     \,T^*_{-\e_3}\,T^{\phantom{*}}_{-\e_3}
   &&&&
  \end{array}
 $$
 It follows that a form $\,\psi\,$  in $\,\L^3_{27}\subset\L^3T\,$
 is harmonic if and only if $\nabla \psi$ is a section of $\Gamma_{1,1}
 \oplus \Gamma_{3,0}$. Note that the expression for $\D$ was obtained
 by adding the Bochner identity, i.~e.~the second Weitzenb\"ock formula,
 to the equation for $\nabla^*\nabla+q(R)$.

 \pfill
 Finally we turn to the case of $\Spin(7)$--holonomy. Irreducible
 $\Spin(7)$--representations are parametrized as $\,\Gamma_{a,\,b,\,c}
 = a\o_1 +b\o_2 + c\o_2$. Again $\Gamma_{0,\,0,\,0} = \C$ is the
 trivial representation and $\Gamma_{0,\,0,\,1} = T\,$ denotes the
 8--dimensional holonomy representation. We want to describe the
 twistor operators for the parallel subbundles of the form bundle.
 For this we need the following representations, which are also
 uniquely determined by their dimension:
 $$
  \begin{array}{lclclclclclclcl}
   \dim \Gamma_{1,\,0,\,0} &\!\!=\!\!&  7 &\;\;&
   \dim \Gamma_{0,\,1,\,0} &\!\!=\!\!& 21 &\;\;&
   \dim \Gamma_{1,\,0,\,1} &\!\!=\!\!& 48 &\;\;&
   \dim \Gamma_{1,\,1,\,0} &\!\!=\!\!& 105
   \\[3pt]
   \dim \Gamma_{2,\,0,\,0} &\!\!=\!\!& 27  &&
   \dim \Gamma_{0,\,0,\,2} &\!\!=\!\!& 35  &&
   \dim \Gamma_{2,\,0,\,1} &\!\!=\!\!& 168 &&
   \dim \Gamma_{1,\,0,\,2} &\!\!=\!\!& 189 \ .
  \end{array}
 $$
 In dimension $112$ there are two different irreducible representations
 denoted by $[112]^a:=\Gamma_{0,\,1,\,1}\,$ and $[112]^b:=
 \Gamma_{0,\,0,\,3}$. Using Poincar\'e duality we can decompose the
 differential forms
 \begin{equation}\label{deco5}
  \begin{array}{ccccccccccc}
   \L^2T^*
   &\cong&
   \L^2_7T^* &\!\!\oplus\!\!& \L^2_{21}T^* &&&&
   &\cong&
   \L^6T^*
   \\[3pt]
   \L^3T^*
   &\cong&
   T^*       &\!\!\oplus\!\!& \L^3_{48}T^* &&&&
   &\cong&
   \L^5T^*
   \\[3pt]
   \L^4T^*
   &\cong&
   \C           &\!\!\oplus\!\!& \L^4_{7}T^*  &\!\!\oplus\!\! &
   \L^4_{27}T^* &\!\!\oplus\!\!& \L^4_{35}T^*
   &&
  \end{array}
 \end{equation}
 into irreducible subspaces, where again the subscripts refer to the
 dimension.

 \pfill
 We start with the representation $V_\l\,=\,\Gamma_{1,\,0,\,0}\,=\,
 \L^2_7$ of highest weight $\l=\o_1$. The relevant weights are $+\e_1$
 and $-\e_4$ with $B$--eigenvalues $\,b_{+\e_1} =\frac12\,$ and
 $\,b_{-\e_4}= -3\,$ leading to the decomposition:
 $$
  T \otimes \Gamma_{1,\,0,\,0}
  \;\;=\;\;
  \Gamma_{1,\,0,\,1} \,\oplus\,\Gamma_{0,\,0,\,1}
  \;\;=\;\;
  \Lambda^3_{48} \,\oplus\, T \ .
 $$
 Because the bundle defined by $\,\Gamma_{1,\,0,\,0}\,$ can be
 considered as the subbundle of the spinor bundle orthogonal to the
 parallel spinor, the curvature endomorphism $q(R)$ is a multiple of
 the scalar curvature and hence vanishes. Thus we obtain on sections
 of $\,\L^2_7\,$ the only Weitzenb\"ock formula
 $$
  0 \;\;=\;\;
  -\tfrac12\,T^*_{+\e_1}\,T^{\phantom{*}}_{+\e_1} \;+\;\;
  3\,T^*_{-\e_4}\,T^{\phantom{*}}_{-\e_4}
 $$
 It follows that $\,\Delta \ge \frac12\, \nabla^*\nabla$, i.~e.~there
 are no non-parallel harmonic forms in $\L^2_7$.

 \pfill
 For the second component of the space of 2--forms $\L^2_{21}\,=\,
 \Gamma_{0,\,1,\,0}\,$ the relevant weights are $+\e_1,\,-\e_2,\,e_3$
 with $B$--eigenvalues $\,1,\,-5,\,-\frac32\,$ in the decomposition:
 $$
  T \otimes \Gamma_{0,\,1,\,0}
  \;\;=\;\;
  \Gamma_{0,\,1,\,1} \,\oplus\,\Gamma_{0,\,0,\,1} \,\oplus\,\Gamma_{1,\,0,\,1}
  \;\;=\;\;
  [112]^a\,\oplus\,T\,\oplus\, \L^3_{48}\ .
 $$
 On sections of $\,\L^2_{21}\,$ we have the Weitzenb\"ock formula
 $$
  q(R)
  \;\;=\;\;
  -\,T^*_{+\e_1}\,T^{\phantom{*}}_{+\e_1}
  \;+\;5\,T^*_{-\e_2}\,T^{\phantom{*}}_{-\e_2}
  \;+\;\tfrac32\,T^*_{+\e_3}\,T^{\phantom{*}}_{+\e_3}\ .
 $$
 Hence a form $\psi$ in $\L^2_{21}$ is harmonic if and only if
 $\nabla \psi$ is a section of $\Gamma_{0,\,1,\,1}$.

 \pfill
 The last parallel subbundle of the form bundle with only one pure
 curvature Weitzenb\"ock formula is $V_\l =\Gamma_{2,\,0,\,0}\,=\,
 \L^4_{27}$. The relevant weights are $+\e_1,\,-\e_4$ with
 $B$--eigenvalues $1,\,-\frac72$ and the decomposition
 $$
  T \otimes \Gamma_{2,\,0,\,0} \;\;=\;\; \Gamma_{2,\,0,\,1}
  \,\oplus\,\Gamma_{1,\,0,\,1} \;\;=\;\; [168]
  \,\oplus\,\Lambda^3_{48} \ ,
 $$
 with $[168]:=\Gamma_{2,0,1}$. On sections of $\,\L^4_{27}\,$ we
 have the Weitzenb\"ock formula:
 $$
  q(R)
  \;\;=\;\;
  -\,T^*_{+\e_1}\,T^{\phantom{*}}_{+\e_1}
  \;+\;\tfrac72\,T^*_{-\e_2}\,T^{\phantom{*}}_{-\e_2} \ .
 $$
 Hence a form $\psi$ in $\L^4_{27}$ is harmonic if and only if $\nabla \psi$
 is a section of $\Gamma_{2,\,0,\,1}$. For the remaining subbundles we have
 at least two pure curvature Weitzenb\"ock formulas one of which is a Bochner
 identity, i.~e.~with a zero curvature term.

 \pfill
 We start with the representation $T\,=\,\Gamma_{0,\,0,\,1}$ describing
 $1$-- and $6$--forms on $M$. The relevant weights are $+\e_1,\,-\e_1,
 +\e_2,+\e_4$ conformal weights or $B$--eigenvalues $\frac34,\,-\frac{21}{4},
 \,-\frac14$ and $-\frac94$ respectively. The corresponding decomposition
 of the representation $T\otimes T$ reads:
 $$
  T \otimes \Gamma_{0,\,0,\,1}
  \;\;=\;\;
  \Gamma_{0,\,0,\,2} \,\oplus\, \Gamma_{0,\,0,\,0} \,\oplus\,
  \Gamma_{0,\,1,\,0} \,\oplus\, \Gamma_{1,\,0,\,0}
  \;\;=\;\;
  \L^4_{35}T^* \,\oplus\, \C        \,\oplus\,
  \L^2_{21}T^* \,\oplus\, \L^2_7T^* \ .
 $$
 In this case we have the universal Weitzenb\"ock formula and the Bochner
 identity (\ref{s7b}) for $(a,b,c) = (0,0,1)$. Since $q(R)\,=\,\Ric $ on
 the tangent bundle we obtain two zero curvature Weitzenb\"ock formulas
 on sections of $T$:
 $$
  \begin{array}{ccrlrlrlrl}
   0
   &=
   &-&\tfrac34    \,T^*_{+\e_1}\,T^{\phantom{*}}_{+\e_1}
   &+&\tfrac{21}4 \,T^*_{-\e_1}\,T^{\phantom{*}}_{-\e_1}
   &+&\tfrac14    \,T^*_{+\e_2}\,T^{\phantom{*}}_{+\e_2}
   &+&\tfrac94    \,T^*_{+\e_4}\,T^{\phantom{*}}_{+\e_4}
   \\[5pt]
   0
   &=
   &+&\!\!\!  15\,T^*_{+\e_1}\,T^{\phantom{*}}_{+\e_1}
   &-&\!\!\! 105\,T^*_{-\e_1}\,T^{\phantom{*}}_{-\e_1}
   &-&\!\!\!  45\,T^*_{+\e_2}\,T^{\phantom{*}}_{+\e_2}
   &+&\!\!\!  75\,T^*_{+\e_4}\,T^{\phantom{*}}_{+\e_4}
  \end{array}
 $$
 Evidently the first equation tells us $\Delta \geq \frac14 \nabla^*\nabla$
 so that every harmonic $1$--form is necessarily parallel. Of course this is
 Bochner's theorem reproved in the case of $\Spin(7)$--manifolds.  Another
 direct consequence is the well--known fact that any Killing vector field
 on a compact $\Spin(7)$--manifold has to be parallel. Indeed Killing vector
 fields are vector fields $X\in\Gamma(TM)$, for which $\nabla X^\sharp\in
 \Gamma(T^*M\otimes T^*M)$ is skew and thus a 2--form. On $\Spin(7)$--manifolds
 this implies $T_{+\e_1}X = 0 = T_{-\e_1}X$ and so all twistor operators
 vanish on $X$.

 \pfill
 Next we consider the representation $\L^4_{35}T^* = \Gamma_{0,\,0,\,2}$
 with relevant weights $\e_1,\,-\e_1,\,\e_2,\,\e_4$, conformal weights or
 $B$--eigenvalues $\frac32,\,-6,\,0,\,-\frac52$ and decomposition
 $$
  T \otimes \Gamma_{0,\,0,\,2}
  \;\;=\;\;
  \Gamma_{0,\,0,\,3} \,\oplus\,\Gamma_{0,\,0,\,1} \,\oplus\,
  \Gamma_{0,\,1,\,1} \,\oplus\,\Gamma_{1,\,0,\,1}
  \;\;=\;\;
  [112]^b\,\oplus\,  T^*\,\oplus\,
  [112]^a \,\oplus\, \Lambda^3_{48}T^* \ .
 $$
 Thus there are two pure curvature Weitzenb\"ock formulas on sections of
 $\L^4_{35}$. For the second we take $\frac1{96}F_{Bochner}$ and obtain
 $$
  \begin{array}{ccrlrlrlrl}
   q(R)
   &=
   &-&\tfrac32 \,T^*_{+\e_1}\,T^{\phantom{*}}_{+\e_1}
   &+&6        \,T^*_{-\e_1}\,T^{\phantom{*}}_{-\e_1}
   &&
   &+&\tfrac52 \,T^*_{+\e_4} \,T^{\phantom{*}}_{+\e_4}
   \\[8pt]
   0
   &=
   &&\tfrac12 \,T^*_{+\e_1}\,T^{\phantom{*}}_{+\e_1}
   &-&2        \,T^*_{-\e_1}\,T^{\phantom{*}}_{-\e_1}
   &-&\phantom2\,T^*_{+\e_2}\,T^{\phantom{*}}_{+\e_2}
   &+&\tfrac32 \,T^*_{+\e_4}\,T^{\phantom{*}}_{+\e_4}
   \\[8pt]
   \Delta
   &=
   &&
   &&5        \,T^*_{-\e_1} \,T^{\phantom{*}}_{-\e_1}
   &&
   &+&5        \,T^*_{+\e_4}\,T^{\phantom{*}}_{+\e_4}
  \end{array}
 $$
 Note that in order to obtain the optimal expression for the operator
 $\Delta$ it was not sufficient to take its definition $\Delta=\nabla^*
 \nabla + q(R)$, we still had to add a multiple of the Bochner identity.
 Last but not least we consider the representation $\L^3_{48}T^*=\Gamma_{1,
 \,0,\,1}$. According to the table at the end of Section \ref{tables} the
 relevant weights are $+\e_1,\,-\e_1,+\e_2,-\e_3,\e_4,-\e_4$ with conformal
 weights or $B$--eigenvalues $\,\frac54,\,-\frac{23}{4},\,\frac14,\,
 -\frac{15}{4},\,-\frac74$ and $\,-\frac{11}{4}\,$ respectively leading
 to
 $$
  \begin{array}{rccccccccccc}
   T \otimes \Gamma_{1,\,0,\,1}
   \quad=\;
   &\Gamma_{1,\,0,\,2} &\!\!\!\oplus\!\!\!
   &\Gamma_{1,\,0,\,0} &\!\!\!\oplus\!\!\!
   &\Gamma_{1,\,1,\,0} &\!\!\!\oplus\!\!\!
   &\Gamma_{0,\,1,\,0} &\!\!\!\oplus\!\!\!
   &\Gamma_{2,\,0,\,0} &\!\!\!\oplus\!\!\!
   &\Gamma_{0,\,0,\,2}
   \\[10pt]
   =\;
   &[\,189\,]    &\!\!\!\oplus\!\!\!
   &\L^2_7T^*    &\!\!\!\oplus\!\!\!
   &[\,105\,]    &\!\!\!\oplus\!\!\!
   &\L^2_{21}T^* &\!\!\!\oplus\!\!\!
   &\L^4_{27}T^* &\!\!\!\oplus\!\!\!
   &\L^4_{35}T^*
  \end{array}
 $$
 On sections of the associated bundle $\L^3_{48}T^*M\subset\L^3T^*M$
 one has three curvature Weitzen\-b\"ock formulas, the formula corresponding
 to $B$ and the Bochner identity $\frac{1}{84}F_{Bochner}$:
 $$
  \begin{array}{ccrlrlrlrlrlrlrl}
   q(R)
   &\!=\;
   &\!\!\!-\!\!\!\!&\tfrac54      \,T^*_{ \e_1} T^{\phantom{*}}_{ \e_1}
   &\!\!\!+\!\!\!\!&\tfrac{23}{4} \,T^*_{-\e_1} T^{\phantom{*}}_{-\e_1}
   &\!\!\!-\!\!\!\!&\tfrac{1}{4}  \,T^*_{ \e_2} T^{\phantom{*}}_{ \e_2}
   &\!\!\!+\!\!\!\!&\tfrac{15}{4} \,T^*_{-\e_3} T^{\phantom{*}}_{-\e_3}
   &\!\!\!+\!\!\!\!&\tfrac74      \,T^*_{ \e_4} T^{\phantom{*}}_{ \e_4}
   &\!\!\!+\!\!\!\!&\tfrac{11}{4} \,T^*_{-\e_4} T^{\phantom{*}}_{-\e_4}
   \\[.5cm]
   0
   &\!=\;
   &\!\!\!\!\!\!\!&\tfrac14      \,T^*_{ \e_1} T^{\phantom{*}}_{ \e_1}
   &\!\!\!-\!\!\!\!&\tfrac{45}{28}\,T^*_{-\e_1} T^{\phantom{*}}_{-\e_1}
   &\!\!\!-\!\!\!\!&\tfrac{3}{4}  \,T^*_{ \e_2} T^{\phantom{*}}_{ \e_2}
   &\!\!\!+\!\!\!\!&\tfrac{27}{28}\,T^*_{-\e_3} T^{\phantom{*}}_{-\e_3}
   &\!\!\!+\!\!\!\!&\tfrac{5}{4}  \,T^*_{ \e_4} T^{\phantom{*}}_{ \e_4}
   &\!\!\!-\!\!\!\!&\tfrac{9}{28} \,T^*_{-\e_4} T^{\phantom{*}}_{-\e_4}
   \\[.5cm]
   \Delta
   &\!=\;
   &\!\!\! \!\!\!\!&
   &\!\!\!\!\!\!\!&\tfrac{36}{7} \,T^*_{-\e_1} T^{\phantom{*}}_{-\e_1}
   &\!\!\! \!\!\!\!&
   &\!\!\!+\!\!\!\!&\tfrac{40}{7} \,T^*_{-\e_3} T^{\phantom{*}}_{-\e_3}
   &\!\!\!+\!\!\!\!&4             \,T^*_{ \e_4} T^{\phantom{*}}_{ \e_4}
   &\!\!\!+\!\!\!\!&\tfrac{24}{7} \,T^*_{-\e_4} T^{\phantom{*}}_{-\e_4}
  \end{array}
 $$
 Consequently a $3$--form $\psi\,\in\,\Gamma(\L^3_{48}T^*M)$ of type
 is harmonic if and only if its covariant derivative $\nabla \psi$
 takes values in $([189]\oplus[105])M\subset T^*M\otimes\L^3_{48}T^*M$
 everywhere.
\section{Bochner Identities in $\G_2$-- and $\Spin(7)$--Holonomy}
\label{kis}
 Aim of this section is to provide a proof of the Bochner identities
 for the holonomies $\g_2$ and $\spin_7$ and thus to complete the
 description of the space of Weitzenb\"ock formulas in these cases.
 Interestingly it seems necessary to introduce a fairly more abstract
 point of view of Weitzenb\"ock formulas in order to get to this point.
\subsection{Universal Weitzenb\"ock Classes and the Kostant Theorem}
 The essential additional twist we will employ in this section is that
 we will basis the study of Weitzenb\"ock formulas on the study of the
 action of central elements of the universal enveloping algebra.
 As a byproduct we get a explicit formula for a
 central element of order $4$ in the universal enveloping algebra of
 $\spin_7$ and a central element of order $4$ in the universal enveloping
 algebra $\U\g_2$.

 Recall that the universal enveloping algebra $\U\g$ of a Lie algebra
 $\g$ is the associative algebra with $1$ generated freely by the vector
 space $\g$ subject only to the commutator relation $X\,Y-Y\,X\,=\,[X,Y]$.
 Thus $\U\g$ is spanned by monomials of the form $X_1\,\ldots\,X_r$ in
 elements $X_1,\,\ldots,\,X_r$ of $\g$ and the filtration $\U^{\leq\bullet}
 \g$ by the degree $r$ of these monomials makes $\U\g$ a filtered algebra.
 Even more important for our purposes is the Hopf algebra structure of
 $\U\g$ with the cocommutative comultiplication
 \begin{equation}\label{comult}
  \D:\;\;\U\g\;\longrightarrow\;\U\g\otimes\U\g,
  \qquad \Q\;\longmapsto\;\sum\,\D_L\Q\otimes\D_R\Q
 \end{equation}
 defined as the unique algebra homomorphism sending $X\in\g$ to $\D X
 \,:=\,X\otimes 1+1\otimes X$ in $\U\g\otimes\U\g$. Defining $\D$ in
 this way clearly implies for all $d,\,r\,\geq\,0$:
 \begin{equation}\label{dft}
  \D(\;\U^{\leq d+r}\g\;)\;\;\subset\;\;\U^{<d}\g\otimes\U\g
  \;+\;\U\g\otimes\U^{\leq r}\g\,.
 \end{equation}
 An integral part of the structure of the universal enveloping algebra
 $\U\g$ is the algebra homomorphism $\U\g\longrightarrow\End\,V$ associated
 to a representation $V$ of $\g$. For finite--dimensional representations
 $V$ the images of these algebra homomorphisms are easily characterized.

 \begin{Lemma}[Bicommutant Theorem]
 \label{bicom}\hfill\break
  Consider a finite dimensional representation $V$ of a semisimple Lie
  algebra $\g$ over $\C$ and the induced representation $\U\g\longrightarrow
  \End\,V$ of $\U\g$. The image of this algebra homomorphism is precisely
  the commutant of the algebra $\End_\g V$ of $\g$--invariant endomorphisms
  $$
   \im(\,\U\g\longrightarrow\End\,V\,)
   \;\;=\;\;
   \{\;A\in\End\,V\;|\;\;
   [\,A,\,F\,]\,=\,0\textrm{\ for all\ }F\,\in\,\End_\g V\;\}\,.
  $$
  In particular the map $\Zent\U\g\longrightarrow\Zent\End_\g V$ is surjective
  for $V$ finite dimensional.
 \end{Lemma}

 \noindent
 The Bicommutant Theorem is actually a special motivating example of von
 Neumann's Bicommutant Theorem, observe that every $*$--subalgebra of
 $\End\,V$ is necessarily von Neumann for a finite dimensional vector
 space $V$. The image of $\U\g$ in $\End\,V$ is the subalgebra generated
 by the $*$--closed subspace $\g$ of $\End\,V$ and thus von Neumann with
 commutant $\End_\g V$. We will give a more geometric proof of this theorem
 based on the Peter--Weyl Theorem in Appendix \ref{biapp}.

 Coming back to Weitzenb\"ock formulas we conclude that for irreducible
 representations $V_\l$ the algebra homomorphism $\U\g\longrightarrow\End\,
 V_\l$ is surjective and hence the same is true for the algebra homomorphism
 $$
  \Phi\,:\qquad\Hom_\g(\,T\otimes T,\,\U\g\,)
  \;\longrightarrow\;\Hom_\g(T\otimes T,\End\,V_\l)
  \;\;=\;\;\W(V_\l)
 $$
 where $\Hom_\g(T\otimes T,\End\,V_\l)$ is one of the interpretation of
 the space $\W(V_\l)$ of Weitzen\-b\"ock formulas on $V_\l M$. Motivated
 by this surjection we will call $\Hom_\g(T\otimes T,\U\g)$ the space of
 universal Weitzenb\"ock formulas. With the universal enveloping algebra
 $\U\g$ being a module over its center the space $\Hom_\g(T\otimes T,\U\g)$
 of universal Weitzenb\"ock formulas is naturally a module for $\Zent\U\g$,
 too, and the filtration $\Hom_\g(T\otimes T,\U^{\leq\bullet}\g)$ turns it
 into a filtered module for the filtration $\Zent^{\leq\bullet}\U\g\,:=\,
 \Zent\U\g\cap\U^{\leq\bullet}\g$ of the center.

 \begin{Definition}[Universal Weitzenb\"ock Classes]
 \label{uwe}\hfill\break
  The space of universal Weitzenb\"ock formulas $\W^{\leq\bullet}\,:=\,
  \Hom_\g(T\otimes T,\U^{\leq\bullet}\g)$ is a filtered module over the
  center $\Zent^{\leq\bullet}\U\g$ of the universal enveloping algebra
  $\U\g$. It splits into the direct sum of filtered $\Zent\,\U\g$--submodules
  called universal Weitzenb\"ock classes:
  $$
   \W^{\leq\bullet}
   \;\;=\;\;\bigoplus_{\alpha}\,\W^{\leq\bullet}_{W_\alpha}
   \;\;:=\;\;\bigoplus_{\alpha}\,\Hom_\g(\,W_\alpha,\,\U^{\leq\bullet}\g\,)\,.
  $$
 \end{Definition}

 \noindent
 It is clear from the definition that with $F\in \W^{\le k}$ also
 $\left.F\right|_{W_\alpha}\in \W^{\le k}_{W_\alpha}$. Moreover the
 powers $B^k$ of the conformal weight operator $B$ are in the image
 of $\W^{\leq k}$ under the surjection $\Phi$. Indeed $B^k$ is the
 image of the nvariant map $p_k:T\otimes T\longrightarrow\U^{\leq k}
 \g$ defined by
 $$
  p_k(\;a\otimes b\;)
  \;\; =\;\;
  \sum_{\mu_1,\ldots,\mu_{k-1}}
  \pr_\g(\,a\wedge t_{\mu_1}\,)\,\pr(\,t_{\mu_1}\wedge t_{\mu_2}\,)\,
  \ldots\,\pr_\g(t_{\mu_{k-1}}\wedge b\,)\, ,
 $$
 where $\pr:\,T\otimes T\longrightarrow\g\subset\Lambda^2T$ is the same
 orthogonal projection used before in the definition of $B$. Under the
 vector space identification $\U\g\cong\S\,\g$ we may consider $p_k(a
 \otimes b)$ as the polynomial $p_k(a\otimes b)[X]=\<X^ka,b>$ on $\g$.
 Important for our considerations below is that the space of universal
 Weitzenb\"ock formulas is a free module over $\Zent\,\U\g$:

 \begin{Theorem}[Kostant's Theorem]\label{kostant}
 \hfill\break
  For every finite dimensional representation $V$ the space $\Hom_\g(\,V,
  \,\U\g\,)$ is a free $\Zent\,\U\g$--module, whose rank over $\Zent\,\U\g$
  agrees with the multiplicity of the zero weight in $V$:
  $$
   \Hom^{\leq\bullet}_\g(\,V,\,\U\g\,)\;\;\cong\;\;
   \Zent\U\g\,\otimes\,\Hom^{\bullet}_\t(\,V,\,\C\,)
  $$
  In particular the module $\Hom^{\leq\bullet}_\g(\,\g,\,\U\g\,)\,\cong\,
  \Zent\,\U\g\otimes\mathrm{Prim}^{\bullet+1}\,\g$ is generated freely as
  a filtered $\Zent\U\g$--module by the primitive elements of $\Zent\,\U\g$
  with degrees shifted by $-1$.
 \end{Theorem}

 \noindent
 As an example we consider holonomy $\g_2$ and the spaces which are
 mapped under $\Phi$ onto the $K$--eigenspaces. We refer to the
 appendix for the other holonomies. Then
 $$
  \begin{array}{lclcl}
   \Hom_{\g_2}(&\!\!\!\C\!\!\!\!&,\;\U\g_2\;)
   & \cong & \Zent\U\g_2
   \\[1ex]
   \Hom_{\g_2}(&\!\!\!\S^2_0T\!\!\!\!&,\;\U\g_2\;)
   & \cong & \Zent\U\g_2\,\<F_2,\,F_4,\,F_6>
   \\[1ex]
   \Hom_{\g_2}(&\!\!\!\g_2\!\!\!\!&,\;\U\g_2\;)
   & \cong & \Zent\U\g_2\,\<F_1,F_5>
   \\[1ex]
   \Hom_{\g_2}(&\!\!\!\g_2^\perp\!\!\!\!&,\;\U\g_2\;)
   & \cong & \Zent\U\g_2\,\langle G_3\rangle
  \end{array}
 $$
 where $F_1,\,F_2,G_3,F_4,F_5$ and $F_6$ are free generators of degree
 $1,2,\ldots,6$. The numbers of generators, i.~e.~the dimension of the
 corresponding zero weight space, can be read off Table~(\ref{wt}). The
 degree of the generators, also called generalized exponents, can be
 obtained by decomposing $\S^k\g_2$ into irreducible components (e.~g.~using
 the program LiE) and by determining the multiplicity of $W_\alpha$ in this
 decomposition for sufficiently many $k$. As mentioned in the theorem:
 the degrees of the generators $F_1,\,F_5$ are the degrees of the generators
 $C_2,\,C_6$ of $\Zent\,\U\g_2$ shifted by one.

 \pfill
 It follows from Kostant's theorem that a basis in the eigenspace
 $\W_{W_\alpha}(V_\l)$ of the classifying endomorphism $K$ may be obtained
 as the image under the surjective representation map $\W_{W_\alpha}
 \longrightarrow\W_{W_\alpha}(V_\l)$ of certain free generators for the
 universal Weitzenb\"ock classes $\W_{W_\alpha}$. Indeed the module
 multiplication with $\Q\,\in\,\Zent\U\,\g$ in $\W$ turns in $\W(V_\l)$
 into multiplication with the value of the central character for $\l$ on
 $\Q$, because:
 $$
  \left.\Q\right|_{V_\l}\;\;=:\;\;\chi_\l(\,\Q\,)\,\id_{V_\l}\, .
  \qquad\qquad
  \chi_\l(\;\Q\;)\;\;=\;\;\tfrac1{\dim\,V_\l}\,\tr_{V_\l}\,\Q
 $$
 In general the value of the central character on $\Q\,\in\,\Zent\,\U\g$ is
 a polynomial in the highest weight $\l$ invariant under the Weyl group of
 $\g$. At least in principle we know the central characters of the higher
 Casimirs $\Cas^{[k]}\,\in\,\Zent^{\leq k}\U\g$ defined in equation
 (\ref{highercasimir}) as traces of the powers of the conformal weight
 operator, since equation (\ref{kth}) implies:
 \begin{equation}\label{ev1}
  \chi_\l(\,\Cas^{[k]}\,)
  \;\;=\;\;
  \sum_\e\,b_\e^k\,\;\frac{\dim\,V_{\l+\e}}{\dim\,V_\l}\ .
 \end{equation}
 In order to proceed we use the diagonal $\D$ of the Hopf algebra $\U\g$
 together with the representation of $\U\g$ on the euclidian vector space
 $T$ to define an algebra homomorphism
 $$
  \D \,:\,\Zent\,\U\g
  \;\stackrel\D\longrightarrow\;(\,\U\g\otimes\U\g\,)^\g
  \;\longrightarrow\;\Hom_\g(\,T\otimes T,\,\U\g\,)\;\;=\;\;\W
 $$
 by $(\D\Q)_{a\otimes b}\,=\,\sum \<a,(\D_L\Q)b>\D_R\Q$
 for all $\Q\,\in\,\Zent\,\U\g$. A particularly nice property of $\D$ is
 that the image of $\D\Q\,\in\,\W$ under the representation map $\Phi:
 \W\,\longrightarrow\,\W(V_\l)$ can be written in the following way:
 \begin{equation}\label{ev}
  \D\Q
  \;\;=\;\;
  \sum_{\e\subset\l}\,\chi_{\l+\e}(\,\Q\,)\,\pr_\e\;\;\in\;\;\W(V_\l)
 \end{equation}
 In fact working our way through the identification $\Hom_\g(T\otimes T,
 \End\,V_\l)\,=\,\End_\g(T\otimes V_\l)$ we find the usual tensor product
 action of $\Q\,\in\,\Zent\,\U\g$ on $T\otimes V_\l$:
 \begin{eqnarray*}
  \Phi(\,\D\Q\,)(\,b\otimes v\,)
  &=&
  \sum_\mu\,t_\mu\otimes (\,\D\Q\,)_{t_\mu\otimes b}\,v
  \;\;=\;\,
  \sum_\mu\,t_\mu\otimes\<t_\mu,(\,\D_L\Q\,)\,b>\,\D_R\Q\,v\\
  &=&
  \sum\,(\,\D_L\Q\,)\,b\otimes(\,\D_R\Q\,)\,v\;\;=\;\;\D\Q(\,b\otimes
  v\,)\;\;=\;\Q(\,b\otimes v\,)
 \end{eqnarray*}
 The last $\D$ is the restriction of the comultiplication to $\Zent \U\g
 \subset\U \g$, which is precisely the action of $\Q\in \Zent\U\g$ on
 $T\otimes V$. Hence we have the following commutative diagram
 $$
  \begin{CD}
   \Zent\U\g & @>{\D}>> & \Hom_\g(\,T\otimes T,\,\U\g\,)
   &@.{\quad\cong\quad} & \bigoplus_{\alpha}\;\Hom_\g(\,W_\alpha,\,\U\g\,) \\
   @V{\Delta}V{\mathrm{}}V && @VV{\mathrm{\Phi}}V && @VV{}V \\
   \End_\g(T\otimes V_\l)& @>{\cong}>> & \Hom_\g(\,T\otimes T,\,\End V_\l\,)
   &@.{\quad\cong\quad} & \bigoplus_{\alpha}\;\Hom_\g(\,W_\alpha,\,\End V_\l\,)
  \end{CD}
 $$
 where the right vertical arrow is the restriction of the
 representation map $\Phi$ onto $\W_{W_\alpha}$. Recall that the left
 square consists of algebra and the right square of
 $\Zent\U\g$--module homomorphisms and that moreover all vertical arrows
 are surjective maps.

 \begin{Example}[Conformal Weight Operator]
 \hfill\label{cwo}\break
  A special case of this construction is the relation between the conformal
  weight operator $B$ and the Casimir. The image of the Casimir $\Cas^{\L^2}
  \,\in\,\Zent\U\g$ becomes
  $$
   \D(\Cas^{\L^2})\;\;=\;\;-2\,B\;+\;(\,\Cas^{\L^2}_T\,+\,\Cas^{\L^2}_{V_\l}\,)
  $$
  because $\D(X^2)\,=\,X^2\otimes 1+2X\otimes X+1\otimes X^2$ for every
  $X\in\g$ and Fegan's Lemma \ref{fegan}. Moreover since $\Delta$ is an algebra
  homomorphism we conclude
  $$
   p(B)
   \;\;=\;\;
   \D\,p(\,-\,{\textstyle\frac12}\,\Cas^{\L^2}
         \,+\,{\textstyle\frac12}(\Cas^{\L^2}_T+\Cas^{\L^2}_{V_\l})\,)
  $$
  for every polynomial $p(B)$ in the conformal weight operator $B$. In
  particular the space of polynomials in $B$ is in the image under
  $\D$ of the subalgebra generated by $\Cas^{\L^2}$.
 \end{Example}

 \noindent
 The crucial additional information we get from introducing the universal
 Weitzenb\"ock classes is the filtration degree of the generators of the
 $\Zent\,\U\g$--modules $\W_{W_\alpha}$. In order to prove the Bochner
 identities for holonomy $\g_2$ and $\spin_7$ we still need the following

 \begin{Lemma}[Filtration Property of $\Delta$]
 \label{dfilt}\hfill\break
  Consider the Weitzenb\"ock class $\W_{W_\alpha}$ associated to an irreducible
  subspace $W_\alpha\,\subset\,T\otimes T$. If there is no non--trivial,
  $\g$--equivariant map from $W_\alpha$ to $\U^{<d}\g$ for some $d\geq1$,
  i.~e.~if $\,\W_{W_\alpha}^{<d}\,=\,\{0\}$, then the composition of $\Delta$
  with the restriction $\mathrm{res}_{W_\alpha}$ to $\W_\alpha$ is filtered
  $$
   \mathrm{res}_{W_\alpha}\,\circ\,\Delta:
   \quad
   \Zent^{\leq d+\bullet}\U\g\;\longrightarrow\;\W^{\leq\bullet}_{W_\alpha},
   \qquad\Q\;\longmapsto\;\left.\Delta\Q\right|_{W_\alpha}
  $$
  of degree $-d$. In particular the restriction $\left.\D\Q\right|_{W_\alpha}
  \,=\,0$ vanishes for all $\Q\,\in\,\Zent^{<2d}\U\g$.
 \end{Lemma}

 \proof
 By the filtration property (\ref{dft}) of the comultiplication we can write
 the diagonal $\D\Q$ of an element $\Q\,\in\,\Zent^{\leq d+r}\U\g$ in a not
 necessarily unique way as a sum of two terms $\D\Q\,=\,\D\Q^{<d}\,+\,
 \D\Q^{\leq r}$ satisfying $\D\Q^{<d}\,\in\,(\U^{<d}\g\otimes\U\g)^\g$
 and $\D\Q^{\leq r}\,\in\,(\U\g\otimes\U^{\leq r}\g)^\g$ respectively.
 Both these summands give rise to $\g$--equivariant, linear maps
 $T\otimes T\,\longrightarrow\,\U\g$ through the pairing of $T\otimes T$
 with the left $\U\g$--factor, explicitly
 $$
  (\,\D\Q^{<d}\,)_{a\otimes b}\;\;:=\;\;
  \sum\<a,\D\Q^{<d}_L\,b>\;\D\Q^{<d}_R
 $$
 with essentially the same formula for $\D\Q^{\leq r}$. By construction
 the map $T\otimes T\,\longrightarrow\,\U\g$ associated to $\D\Q^{\leq r}$
 maps into $\U^{\leq r}\g$, while the map associated to $\D\Q^{<d}$ vanishes
 upon restriction to $W_\alpha\,\subset\,T\otimes T$, because by assumption
 there is no non--trivial, $\g$--invariant, pairing of $W_\alpha$ with the
 left image of $\D\Q^{<d}$ defined by:
 $$
  \mathrm{span}\{\;\D\Q^{<d}_L\;|\;\;\D\Q^{<d}\,=\,
  \sum\,\D\Q^{<d}_L\,\otimes\,\D\Q^{<d}_R\;\;\}
  \;\;\subset\;\;
  \U^{<d}\g
 $$
 For the second statement we note $\left.\D\Q\right|_{W_\alpha}\,\in\,
 \W^{<d}_\alpha=\{0\}$ for all $\Q\in\Zent^{\le 2d-1}\U\g$.
 \qed
\subsection{Proof of the Bochner Identities in Holonomy $\g_2$ and
 $\spin_7$}
 Let us now discuss the details of the proof of the additional Bochner
 identity in $\G_2$--holonomy. Applying the Gram--Schmidt orthogonalization
 process \ref{orp} to the powers $1,\,B,\,B^2$ and $B^3$ of the conformal
 weight operator we obtained in Equations~(\ref{g2-1}) and
 (\ref{g2-2}) a sequence $p_0(B), p_1(B), p_2(B), p_3(B)$ of
 $\tau$--eigenvectors. In order to proceed with the recursion procedure
 it remains to be shown that $p_3(B)$ is a $K$--eigenvector.

 We know that $p_3(B)$ is a $-1$ eigenvector of $\tau$, orthogonal to
 $B$ and expressible as polynomial in $B$ of degree $3$. Thus $p_3(B)$
 is an element in the image of $\W^{\leq 3}$ in $\W(V_\l)$ and  can be
 written as a sum $p_3(B)\,=\,p_3(B)_{\g_2}+ p_3(B)_{\g_2^\perp}$ of
 two vectors $p_3(B)_{\g_2}$ and $p_3(B)_{\g_2^\perp}$ in the image of
 $\W^{\leq3}_{\g_2}$ and $\W^{\leq 3}_{\g_2^\perp}$ in $\W(V_\l)$ respectively.
 However the image of $\W^{\leq3}_{\g_2}$ in $\W(V_\l)$ is spanned by
 $B$, because the filtered $\Zent\,\U\g_2$--module $\W_{\g_2}$ is
 generated by two elements in degrees $1$ and $5$ and the representation
 $\W_{\g_2}\longrightarrow\W_{\g_2}(V_\l)$ turns module multiplication
 into multiplication by the central character $\chi_\l$. Consequently
 the vector $p_3(B)$ is orthogonal to the image of $\W^{\leq 3}_{\g_2}$
 in $\W(V_\l)$ and lies in the eigenspace $\W_{\g_2^\perp}(V_\l)$ of the
 classifying endomorphism $K$:

 \begin{Theorem}[Bochner Identity in $\G_2$--Holonomy]
 \hfill\label{bochig2}\break
  The following cubic polynomial in the conformal weight operator defines
  an eigenvector for the classifying endomorphisms $K$ of eigenvalue $-2$:
  $$
   p_3(B)
   \;\;:=\;\;
   B^3\;+\;\tfrac{13}3\,B^2\;+\;(\,\tfrac12\Cas^{\L^2}_{V_\l}\,+\,4\,)\,B
   \;+\;\tfrac23\,\Cas^{\L^2}_{V_\l}
  $$
  Inserting the eigenvalues or conformal weights $b_\e$ of $B$ we arrive at
  equation \ref{g2bx}.
 \end{Theorem}

 \pfill
 In the last part of this section we will prove the Bochner identity
 for holonomy $\spin_7$. As in the $\g_2$--case we apply the Gram--Schmidt
 orthogonalization process \ref{orp} to the powers $1,\,B,\,B^2$ and
 $B^3$ and obtain in Equations~(\ref{s7-1}) and (\ref{s7-2}) a sequence
 $p_0(B)$, $p_1(B)$, $p_2(B)$, $p_3(B)$ of $\tau$--eigenvectors. Again
 $p_3(B)$ is a $(-1)$--eigenvector of $\tau$, orthogonal to $B$ and
 expressible as polynomial in $B$ of degree $3$ so that the summands
 in the decomposition $p_3(B)\,=\,p_3(B)_{\spin_7}+p_3(B)_{\spin_7^\perp}$
 are in the image of $\W^{\leq3}_{\spin_7}$ and $\W^{\leq 3}_{\spin_7^\perp}$
 in $\W(V_\l)$ respectively. Of course we want to extract the Bochner identity
 $p_3(B)_{\spin_7^\perp}$ from $p_3(B)$. At this point the argument in
 the $\Spin(7)$--case becomes more complicated, because the $\Zent\,\U
 \spin_7$--module $\W_{\spin_7}$ has generators in degree $1$, $3$ and
 $5$ so that the image of $\W_{\spin_7}^{\leq3}$ in $\W_{\spin_7}(V_\l)$
 has dimension two. Even with $p_3(B)$ orthogonal to $B$ we may thus not
 conclude that the component $p_3(B)_{\spin_7}=0$ vanishes. The idea to
 cope with this complication is to construct an element $Q_\l\in \Zent\U
 \spin_7$ depending polynomially on the highest weight $\l$ such that
 $\D Q_\l\in\W_{\spin_7}(V_\l)$ is orthogonal to $B$. The Bochner identity
 is then the projection of $p_3(B)$ onto the orthogonal complement of
 $\D Q_\l$.

 \pfill
 In the 4-dimensional space $\Zent^{\leq4}\U\spin_7$ we look for an element
 $\Q_\l$ as a linear combination of the base vectors $\1,\,\Cas,\,\Cas^2$
 and $\Cas^{[4]}$ with unknown coefficients. We know $\D \Cas$ and $\D \Cas^2$
 from Example~\ref{cwo} and $\D \Cas^{[4]}$ from Equations~(\ref{ev1}) and
 (\ref{ev}) so that
 $$
  \<\D Q_\l,\1>\;\;=\;\;0
  \qquad\qquad
  \<\D Q_\l ,B >\;\;=\;\;0
  \qquad\qquad
  \<\D Q_\l ,B^2 >\;\;=\;\;0
 $$
 turn into three linear independent equations for the four unknown
 coefficients. Using a computer algebra system to do the necessary
 calculations we find the convenient solution
 \begin{eqnarray*}
  \Q_\l\!
  &=&
  2\,\Cas^{\L^2}_{V_\l}\,\Cas^{[4]}
  \;-\;160\,\Cas^{\L^2}_{V_\l}\,(\Cas^{\L^2})^2\\
  &&
  +\,(\,320\,(\Cas^{\L^2}_{V_\l})^2
  \;-\;1184\,\Cas^{\L^2}_{V_\l}\;-\;4\,\Cas^{[4]}_{V_\l}\,)\,\Cas^{\L^2}\\
  &&
  +\,(\,-\,160\,(\Cas^{\L^2}_{V_\l})^3\,+\,2\,
  \Cas^{\L^2}_{V_\l}\,\Cas^{[4]}_{V_\l}\,+\,1712\,(\Cas^{\L^2}_{V_\l})^2
  \,-\,9408\,\Cas^{\L^2}_{V_\l}\,-\,21\,\Cas^{[4]}_{V_\l}\,)
 \end{eqnarray*}
 in $\,\Zent^{\leq4}\U\spin_7$, where we denote the eigenvalue of the
 central element $\Cas^{[4]}\in\U\spin_7$ on the irreducible representation
 $V_\l$ by $\Cas^{[4]}_{V_\l}$ in analogy to the eigenvalues of $\Cas^{\L^2}$.

 By construction $\D\Q_\l$ is orthogonal to $1,\,B$ and $B^2$ and we conclude
 that $\D Q_\l$ is indeed an eigenvector for the classifying endomorphism $K$.
 In fact the component $\left.\D\Q_\l\right|_{\spin_7^\perp}\,=\,0$ vanishes
 according to Lemma \ref{dfilt}, because $\Q_\l$ has degree $4$ and there is
 no non--trivial equivariant map $\spin_7^\perp\longrightarrow\U^{<3}\spin_7$.
 Similarly the components $\left.\D\Q_\l\right|_{\S^2_\circ T}\,=\,0$
 and $\left.\D\Q_\l\right|_{\C}\,=\,0$ are trivial, since the image of
 $\W^{\leq2}_{\S^2_\circ T}$ in $\W(V_\l)$ has dimension $1$ spanned by
 $p_2(B)$ while $\W_\C(V_\l)$ is spanned by $p_0(B)$. With $\D\Q_\l\,=\,
 \left.\D\Q_\l\right|_{\spin_7}$ being an eigenvector of $K$ orthogonal
 to $p_1(B)\,=\,B$ the problematic component $p_3(B)_{\spin_7}$ of $p_3(B)$
 must be a multiple of $\D\Q_\l$. In consequence the complementary component
 $p_3(B)_{\spin_7^\perp}$ of $p_3(B)$ is the projection of $p_3(B)$ onto the
 orthogonal complement of $\D Q_\l$ and may serve as the $\spin_7$--Bochner
 identity. Using again a computer algebra system for the necessary
 calculations we find that this projection of $p_3(B)$ to the orthogonal
 complement of $\Delta\Q_\l$ agrees with the endomorphism $F_{\mathrm{Bochner}}
 \in\W(V_\l)$  specified in equation (\ref{s7b}):

 \begin{Theorem}[Bochner Identity in $\Spin(7)$--Holonomy]
 \hfill\label{bochis7}\break
  The endomorphism $F_{\mathrm{Bochner}}\,\in\,\W(V_\l)$ specified in
  equation (\ref{s7b}) with components
  $$
   \begin{array}{rcccc}
    F_{\mathrm{Bochner}}\;\;=\;\;
    +&\!     c  \!&\!(2b+c+2)\!&\!(2a+2b+c+4)\!&\!\pr_{+\e_1}\cr
    -&\!(c+2)\!&\!(2b+c+4)\!&\!(2a+2b+c+6)\!&\!\pr_{-\e_1}\cr
    -&\!(c+2)\!&\!(2b+c+2)\!&\!(2a+2b+c+4)\!&\!\pr_{+\e_2}\cr
    +&\!     c  \!&\!(2b+c+4)\!&\!(2a+2b+c+6)\!&\!\pr_{-\e_2}\cr
    -&\!     c  \!&\!(2b+c+4)\!&\!(2a+2b+c+4)\!&\!\pr_{+\e_3}\cr
    +&\!(c+2)\!&\!(2b+c+2)\!&\!(2a+2b+c+6)\!&\!\pr_{-\e_3}\cr
    +&\!(c+2)\!&\!(2b+c+4)\!&\!(2a+2b+c+4)\!&\!\pr_{+\e_4}\cr
    -&\!     c  \!&\!(2b+c+2)\!&\!(2a+2b+c+6)\!&\!\pr_{-\e_4}
   \end{array}
  $$
  is an eigenvector of the classifying endomorphism $K$ for the eigenvalue
  $\,-\frac94$.
 \end{Theorem}
\begin{appendix}
\section{Geometric Proof of the Bicommutant Theorem}\label{biapp}
 Consider a finite--dimensional representation $V$ of a semisimple Lie group
 $G$ with Lie algebra $\g$ and the algebra homomorphism $\U\g\longrightarrow
 \End\,V$ associated to the infinitesimal representation of $\g$ on $V$.
 Evidently the image of $\U\g$ under this homomorphism is the subalgebra
 $\mathcal{A}_\g V$ of $\End\,V$ generated by the image of $\g$. An
 alternative characterization of this subalgebra can be given using
 the notion of the centralizer subalgebra or more succinctly the commutant
 $\mathrm{Comm}\,\mathcal{A}$ of a subset $\mathcal{A}\,\subset\,\End\,V$:
 $$
  \mathrm{Comm}\,\mathcal{A}\;\;:=\;\;
  \{\;F\,\in\,\End\,V\;|\;\;[\,F,\,A\,]\textrm{\ for all\ }A\,\in\,
  \mathcal{A}\;\}
 $$
 Von Neumann's famous Bicommutant Theorem states that the commutant of
 the commutant of a subset $\mathcal{A}\,\subset\,\End\,V$ of the star
 algebra of bounded operators on a complex Hilbert space $V$ is precisely
 the von Neumann subalgebra generated by $\mathcal{A}\,\cup\,\mathcal{A}^*$.
 Every subalgebra of $\End\,V$ for a finite--dimensional Hilbert space $V$
 is von Neumann and with the image of $\g\,=\,\g^*$ being closed under taking
 adjoints we can characterize the subalgebra $\mathcal{A}_\g V$ generated
 by $\g$ as the commutant of the subalgebra $\End_\g V$ of $\g$--invariant
 endomorphisms on $V$:

 \begin{Lemma}[Bicommutant Theorem]
 \label{double}\hfill\break
  Let $V$ be a finite--dimensional representation of a complex semisimple
  Lie algebra $\g$. The image of the representation homomorphism $\U\g
  \longrightarrow\End\,V$ is precisely the commutant $\mathrm{Comm}\,\End_\g
  V$ of the subalgebra $\End_\g V\,\subset\,\End\,V$ of $\g$--invariant
  endomorphisms on $V$.
 \end{Lemma}

 \proof
 With $\g$ being a complex semisimple Lie algebra it has a compact real form
 $\g_\R$ defined by a real structure $X\longmapsto\overline{X}$ on $\g$. Being
 a compact real Lie algebra $\g_\R$ is the Lie algebra of a compact simply
 connected Lie group $G_\R$ and hence there exists a $\g_\R$--invariant
 hermitian form $(,)$ on the complex vector space $V$ satisfying $X^*\,
 =\,-\overline{X}$ for all $X\,\in\,\g$. Thus the image of $\g$ in
 $\End\,V$ is closed under taking adjoints and by von Neumann's Theorem we
 conclude that the algebra $\mathcal{A}_\g V$ generated by $\g$ is the
 bicommutant $\mathrm{Comm}^2\g\,=\,\mathrm{Comm}\,\End_\g V$ of $\g$.

 Actually alluding to von Neumann's theorem may seem somewhat strange for
 this purely representation theoretic lemma, and for this reason we want
 to sketch a more elementary proof making use of Schur's Lemma and the
 Theorem of Peter--Weyl instead. In this way the reader will presumably
 get a better understanding of the crux of the statement, which interestingly
 enough was the main observation leading von Neumann to his Bicommutant
 Theorem in the first place. As $V$ is a representation of a semisimple
 Lie algebra $\g$ defined over $\C$ we can decompose $V$ completely into
 irreducible subrepresentations $V_\l\,\subset\,V$. With respect to
 the general form of this decomposition
 $$
  V\;\;=\;\;\bigoplus_\l\Hom_\g(\,V_\l,\,V\,)\,\otimes\,V_\l
 $$
 the commutant of the image of $\g$ in $\End\,V$ is precisely the subalgebra
 $$
  \End_\g V\;\;=\;\;
  \bigoplus_\l\End\,\Hom_\g(\,V_\l,\,V\,)\,\otimes\,\id_{V_\l}
 $$
 by Schur's Lemma, which is a direct sum of matrix algebras one for each
 non--trivial isotypic component. The crucial observation underlying both
 the theorem above and von Neumann's theorem is the fact that the commutant
 of such a subalgebra is given by:
 $$
  \mathrm{Comm}\,\End_\g V\;\;=\;\;
  \bigoplus_\l\id_{\Hom_\g(V_\l,V)}\,\otimes\,\End\,V_\l
 $$
 Skipping the details of this argument using only elementary linear algebra
 we conclude that without loss of generality all isotypic components of $V$
 can be assumed irreducible. In other words we may assume that $V\,=\,
 \bigoplus V_\l$ is the direct sum of a finite set of irreducible, pairwise
 non--isomorphic representations $V_\l$ and we want to show that the
 resulting representation homomorphism is surjective:
 $$
  \U\g\;\longrightarrow\;\bigoplus\End\,V_\l\;\;\subset\;\;
  \End\left(\,\bigoplus V_\l\,\right)\ .
 $$
 In order to proceed we reformulate this statement thinking
 of $\U\g$ as the algebra of left--invariant scalar differential operators
 on the complex--valued functions $C^\infty(\,G_\R,\,\C\,)$ on $G_\R$, in
 particular there exists a linear map $\U\g\longrightarrow C^\infty(\,G_\R,
 \,\C\,)^*,\;\Q\longmapsto\mathrm{ev}_e\circ\Q$ by composing the action of the
 differential operator $\Q\,\in\,\U\g$ with the evaluation at the identity.
 On the other hand the Theorem of Peter--Weyl identifies the direct sum
 $\bigoplus\End\,V_\l$ with a subspace of $C^\infty(\,G_\R,\,\C\,)$ as all
 $V_\l$ are irreducible and non--isomorphic. Working out the details of
 this identification we see that the resulting linear map reads:
 $$
  \U\g\;\longrightarrow\;\left(\;\bigoplus\End\,V_\l\;\right)^*,
  \qquad \Q\;\longmapsto\;\left(\;\eta\otimes e\;\longmapsto\;
  \eta(\,\Q e\,)\;\right)
 $$
 Now all $\End\,V_\l\,\cong\,(\,\End\,V_\l\,)^*$ are real representations
 with invariant symmetric bilinear form given by the trace and $\eta(\Q e)
 \,=\,\tr(\,\Q\circ(\eta\otimes e)\,)$ implies that the linear map constructed
 above is just another interpretation of the representation homomorphism.
 Moreover it is well--known that all matrix coefficients $\End\,V_\l\,
 \subset\,C^\infty(\,G_\R,\,\C\,)$ of an irreducible representation are
 actually analytic functions with respect to the natural analytic structure
 on $G_\R$. Consequently all functions in the finite direct sum $\bigoplus
 \End\,V_\l$ are analytic as well and can thus be separated by their partial
 derivatives at the identity. In other words the second interpretation of
 the representation homomorphism is surjective simply because the matrix
 coefficients of an irreducible representation are analytic.
 \qed

 \begin{Corollary}
 \label{zentrum}\hfill\break
  The representation homomorphism $\U\g\longrightarrow\End\,V$ for a
  finite--dimensional representation $V$ of a complex semisimple Lie
  algebra $\g$ restricts to a surjective homomorphism:
  $$
   \Zent\,\U\g\;\longrightarrow\;\Zent\,\End_\g V
  $$
 \end{Corollary}

 \noindent
 Without giving a formal proof of this corollary we remark that the very
 definition of the commutant reads $\Zent\End_\g V\,=\,\End_\g V\,\cap\,
 \mathcal{A}_\g V$ while $\mathcal{A}_\g V$ is the image of the universal
 enveloping algebra $\U\g$ under the representation homomorphism by Lemma
 \ref{double}. Consequently every element in $\Zent\End_\g V$ is the image
 of some element of $\U\g$. With the representation homomorphism being
 $\g$--equivariant and $\g$ semisimple we can actually choose a
 $\g$--invariant preimage $\Q\,\in\,\U\g$ for a $\g$--invariant
 endomorphism or equivalently a preimage in the center $\Zent\U\g$.
 Thinking a little bit more about the proof of Lemma \ref{double}
 we can formulate a slightly different useful result about the
 structure of $\Zent\End_\g V$:

 \begin{Corollary}
 \label{structure}\hfill\break
  The center $\Zent\,\End_\g V$ of the algebra $\End_\g V$ of $\g$--invariant
  endomorphisms on a finite--dimensional representation $V\,=\,\bigoplus
  \Hom_\g(\,V_\l,\,V\,)\otimes V_\l$ of a complex semisimple Lie algebra
  $\g$ is spanned by the $\g$--invariant projections onto the isotypical
  components:
  $$
   \Zent\End_\g V\;\;=\;\;
   \bigoplus_\l\;\C\;\id_{\Hom_\g(V_\l,V)}\otimes\id_{V_\l}
  $$
 \end{Corollary}

 \noindent
 For applications in the Weitzenb\"ock machine the representation $V$
 will usually be a tensor product $T\otimes V_\l$ of an irreducible
 representation $V_\l$ of highest weight $\l$ defined over $\C$ with
 the complexified holonomy representation $T$ of the holonomy group
 $G_\R$ with Lie algebra $\g_\R$. In this case the algebra $\End_g(\,
 T\otimes V_\l\,)=\Zent\,\End_\g(\,T\otimes V_\l\,)$ is commutative.
\section{Module Generators and Higher Casimirs}
 \begin{Remark}[Module Generators for $\Zent\,\U\so_{2r+1}$]
 \hfill\label{ksoodd}\break
  The center of the universal enveloping algebra of $\so_{2r+1},\,r\geq 1,$
  is a free polynomial algebra $\Zent\U\so_{2r+1}\,=\,\C[P^{[2]},P^{[4]},
  \ldots,P^{[2r]}]$ in $r$ generators of degree $2,\,4,\,\ldots,\,2r$.
  Moreover:
  $$
   \begin{array}{lclcl}
    \Hom_{\so_{2r+1}}(&\!\!\!\C\!\!\!\!&,\;\U\so_{2r+1}\;)
    & \cong & \Zent \U\so_{2r+1}
    \\[1ex]
    \Hom_{\so_{2r+1}}(&\!\!\!\S^2_0T\!\!\!\!&,\;\U\so_{2r+1}\;)
    & \cong & \Zent \U\so_{2r+1}\,\<F_2,\,F_4,\,\ldots,\,F_{2r}>
    \\[1ex]
    \Hom_{\so_{2r+1}}(&\!\!\!\so_{2r+1}\!\!\!\!&,\;\U\so_{2r+1}\;)
    & \cong & \Zent\U\so_{2r+1}\, \<F_1,\,F_3,\,\ldots,\,F_{2r-1}>
  \end{array}
  $$
 \end{Remark}

 \begin{Remark}[Module Generators for $\Zent\,\U\so_{2r}$]
 \hfill\label{ksoeven}\break
  The center of the universal enveloping algebra $\U\so_{2r}$ of $\so_{2r},
  \,r\geq 2,$ is a free polynomial algebra $\Zent\U\so_{2r}\,=\,\C[P^{[2]},
  P^{[4]},\ldots,P^{[2r-2]},E^{[r]}]$ in $r-1$ generators of degree $2,\,4,
  \,\ldots,\,2r-2$ respectively and one additional generator in degree $r$.
  Moreover:
  $$
   \begin{array}{lclcl}
    \Hom_{\so_{2r}}(&\!\!\!\C\!\!\!\!&,\;\U\so_{2r}\;)
    & \cong & \Zent \U\so_{2r}
    \\[1ex]
    \Hom_{\so_{2r}}(&\!\!\!\S^2_0T\!\!\!\!&,\;\U\so_{2r}\;)
    & \cong & \Zent \U\so_{2r}\,\<F_2,\,F_4,\,\ldots,\,F_{2r-2}>
    \\[1ex]
    \Hom_{\so_{2r}}(&\!\!\!\so_{2r}\!\!\!\!&,\;\U\so_{2r}\;)
    & \cong & \Zent\U\so_{2r}\, \<F_1,\,F_3,\,\ldots,\,F_{2r-3},G_{r-1}>
  \end{array}
  $$
 \end{Remark}

 \begin{Remark}[Module Generators for $\Zent\,\U\g_2$]
 \hfill\label{kg2}\break
  The center of the universal enveloping algebra $\U\g_2$ of $\g_2$
  is a free polynomial algebra $\Zent\U\g_2\,=\,\C[\Cas^{[2]},\Cas^{[6]}]$
  in two generators of degree $2$ and $6$. Moreover:
  $$
   \begin{array}{lclcl}
    \Hom_{\g_2}(&\!\!\!\C\!\!\!\!&,\;\U\g_2\;)
    & \cong & \Zent\U\g_2
    \\[1ex]
    \Hom_{\g_2}(&\!\!\!\S^2_0T\!\!\!\!&,\;\U\g_2\;)
    & \cong & \Zent\U\g_2\,\<F_2,\,F_4,\,F_6>
    \\[1ex]
    \Hom_{\g_2}(&\!\!\!\g_2\!\!\!\!&,\;\U\g_2\;)
    & \cong & \Zent\U\g_2\,\<F_1,F_5>
    \\[1ex]
    \Hom_{\g_2}(&\!\!\!\g_2^\perp\!\!\!\!&,\;\U\g_2\;)
    & \cong & \Zent\U\g_2\,\langle G_3\rangle
  \end{array}
  $$
 \end{Remark}

 \begin{Remark}[Module Generators for $\Zent\,\U\spin_7$]
 \hfill\label{kspin_7}\break
  The center of the universal enveloping algebra $\U\spin_7$ of $\spin_7$
  is a free polynomial algebra $\Zent\U\spin_7\,=\,\C[\Cas^{[2]},\Cas^{[4]},
  \Cas^{[6]}]$ in three generators of degree $2,\,4$ and $6$. Moreover:
  $$
   \begin{array}{lclcl}
    \Hom_{\spin_7}(&\!\!\!\C\!\!\!\!&,\;\U\spin_7\;)
    & \cong & \Zent\U\spin_7
    \\[1ex]
    \Hom_{\spin_7}(&\!\!\!\S^2_0T\!\!\!\!&,\;\U\spin_7\;)
    & \cong & \Zent\U\spin_7\, \<F_2,\,F_4,\,F_6>
    \\[1ex]
    \Hom_{\spin_7}(&\!\!\!\spin_7\!\!\!\!&,\;\U\spin_7\;)
    & \cong & \Zent \U\spin_7\,\<F_1,\,F_3,\,F_5>
    \\[1ex]
    \Hom_{\spin_7}(&\!\!\!\spin_7^\perp\!\!\!\!&,\;\U\spin_7\;)
    & \cong & \Zent\U\spin_7\,\langle G_3 \rangle
   \end{array}
  $$
 \end{Remark}

 \begin{Remark}[Higher Casimirs for $\G_2$]
 \hfill\label{g2cas}\break
  The eigenvalues of the generators $\Cas^{[2]},\,\Cas^{[6]}$ of $\Zent\,
  \U\g_2$ of degrees 2 and 6 respectively on the irreducible representation
  $V_\l$ of highest weight $\l\,=\,a\o_1\,+\,b\o_2$ are given by:
  \begin{eqnarray*}
   \frac34\,\Cas^{[2]}_{V_\l}
   &=& a^2\;+\;3\,ab\;+\;3\,b^2\;+\;5\,a\;+\;9\,b\\[1.5ex]
   \frac{243}{11}\,\Cas^{[6]}_{V_\l}
   &=&
   4\,a^6\;+\;36\,a^5b\;+\;117\,a^4b^2
   \;+\;162\,a^3b^3\;+\;81\,a^2b^4\\
   &&
   +\;60\,a^5\;+\;414\,a^4b\;+\;954\,a^3b^2
   \;+\;810\,a^2b^3\;+\;162\,ab^4\\
   &&
   -\;408\,a^4\;-\;2808\,a^3b\;-\;8829\,a^2b^2
   \;-\;12636\,ab^3\;-\;6804\,b^4\\
   &&
   -\;6580\,a^3\;-\;33174\,a^2b\;-\;61362\,ab^2
   \;-\;40824\,b^3\\
   &&
   -\;6396\,a^2\;-\;32508\,ab\;-\;27756b^2\\
   &&
   +\;56520\,a\;+\;100440\,b
  \end{eqnarray*}
 \end{Remark}

 \begin{Remark}[Higher Casimirs for $\Spin(7)$]
 \hfill\label{spin7cas}\break
  The eigenvalues of the generators $\Cas^{[2]}$, $\Cas^{[4]}$ and
  $\Cas^{[6]}$ of $\Zent\,\U\spin_7$ of degrees $2,\,4$ and $6$
  on the irreducible representation $V_\l$ of highest weight
  $\l\,=\,a\o_1\,+\,b\o_2\,+\,c\o_3$ are:
  \begin{eqnarray*}
   2\,\Cas^{[2]}_{V_\l}
   &=&
   4\,a^2\;+\;8\,b^2\;+\;3\,c^2\;+\;8\,ab\;+\;4\,ac\;+\;8\,bc
   \;+\;20\,a\;+\;32\,b\;+\;18\,c\\[1.5ex]
   32\,\Cas^{[4]}_{V_\l}
   &=&
   16\,a^4
   \;+\;128\,b^4
   \;+\;21\,c^4
   \;+\;192\,a^2b^2
   \;+\;72\,a^2c^2
   \;+\;240\,b^2c^2\\
   &&
   \;+\;32\,a^3c
   \;+\;64\,a^3b
   \;+\;256\,b^3c
   \;+\;256\,b^3a
   \;+\;56\,c^3a
   \;+\;112\,c^3b\\
   &&
   \;+\;192\,a^2bc
   \;+\;384\,b^2ac
   \;+\;240\,c^2ab\\
   &&
   \;+\;160\,a^3
   \;+\;1024\,b^3
   \;+\;252\,c^3
   \;+\;768\,a^2b
   \;+\;432\,a^2c
   \;+\;1536\,b^2a\\
   &&
   \;+\;1632\,b^2c
   \;+\;1056\,c^2b
   \;+\;552\,c^2a
   \;+\;1632\,abc\\
   &&
   \;+\;800\,a^2
   \;+\;1152\,c^2
   \;+\;3040\,b^2
   \;+\;3040\,ab
   \;+\;1760\,ac
   \;+\;3424\,bc\\
   &&
   \;+\;2000\,a
   \;+\;3968\,b
   \;+\;2376\,c\\[1.5ex]
   512\,\Cas^{[6]}_{V_\l}
   &=&
   64\,a^6
   \;+\;2048\,b^6
   \;+\;183\,c^6
   \;+\;384\,a^5b
   \;+\;192\,a^5c
   \;+\;6144\,b^5c\\
   &&
   \,+\,6144\,b^5a
   \,+\,732\,c^5a
   \,+\,1464\,c^5b
   \,+\,1920\,a^4b^2
   \,+\,720\,a^4c^2
   \,+\,7680\,b^4a^2\\
   &&
   \;+\;9600\,b^4c^2
   \;+\;1260\,c^4a^2
   \;+\;4920\,c^4b^2
   \;+\;1920\,a^4bc
   \;+\;15360\,b^4ac\\
   &&
   \;+\;4920\,c^4ab
   \;+\;5120\,a^3b^3
   \;+\;1120\,a^3c^3
   \;+\;8960\,b^3c^3
   \;+\;7680\,a^3b^2c\\
   &&
   \,+\,4800\,a^3c^2b
   \,+\,15360\,b^3a^2c
   \,+\,19200\,b^3c^2a
   \,+\,6720\,c^3a^2b
   \,+\,13440\,c^3b^2a\\
   &&
   \;+\;14400\,a^2b^2c^2\\
   &&
   \;+\;960\,a^5
   \;+\;24576\,b^5
   \;+\;3294\,c^5
   \;+\;7680\,a^4b
   \;+\;4320\,a^4c
   \;+\;61440\,b^4a\\
   &&
   \;+\;65280\,b^4c
   \;+\;11100\,c^4a
   \;+\;22080\,c^4b
   \;+\;30720\,a^3b^2
   \;+\;11040\,a^3c^2\\
   &&
   \;+\;61440\,b^3a^2
   \;+\;84480\,b^3c^2
   \;+\;15120\,c^3a^2
   \;+\;60000\,c^3b^2
   \;+\;32640\,a^3bc\\
   &&
   \;+\;130560\,b^3ac
   \;+\;60000\,c^3ab
   \;+\;97920\,a^2b^2c
   \;+\;63360\,a^2c^2b\\
   &&
   \;+\;126720\,b^2c^2a\\
   &&
   \;+\;9600\,a^4
   \;+\;167424\,b^4
   \;+\;32592\,c^4
   \;+\;67200\,a^3b
   \;+\;38400\,a^3c\\
   &&
   \;+\;334848\,b^3a
   \;+\;365568\,b^3c
   \;+\;88032\,c^3a
   \;+\;175584\,c^3b
   \;+\;234624\,a^2b^2\\
   &&
   \;+\;92832\,a^2c^2
   \;+\;364128\,b^2c^2
   \;+\;257664\,a^2bc
   \;+\;548352\,b^2ac\\
   &&
   \;+\;364128\,c^2ab\\
   &&
   \;+\;56000\,a^3
   \;+\;684032\,b^3
   \;+\;193464\,c^3
   \;+\;413952\,a^2b
   \;+\;251808\,a^2c\\
   &&
   \,+\,993024\,b^2a
   \,+\,1158912\,b^2c
   \,+\,397968\,c^2a
   \,+\,790656\,c^2b
   \,+\,1125888\,abc\\
   &&
   \;+\;160000\,a^2
   \;+\;1321856\,b^2
   \;+\;562848\,c^2
   \;+\;1189760\,ab
   \;+\;759040\,ac\\
   &&
   \;+\;1607552\,bc
   \;+\;200000\,a
   \;+\;863744\,b
   \;+\;606240\,c
  \end{eqnarray*}
 \end{Remark}
\end{appendix}
\vskip1cm
\begin{center}
 \qquad
 \parbox{60mm}{Uwe Semmelmann\\
  Mathematisches Institut\\
  Universit{\"a}t zu K{\"o}ln\\
  Weyertal 86-90\\
  D-50931 K{\"o}ln, Germany\\[1mm]
  \texttt{semmelma@math.uni-koeln.de}}
 \qquad\qquad
 \parbox{60mm}{Gregor Weingart\\
  Institut f\"ur Mathematik\\
  Humboldt--Universit{\"a}t zu Berlin\\
  Unter den Linden\\
  D-10099 Berlin, Germany\\[1mm]
  \texttt{weingart@math.hu-berlin.de}}
 \end{center}

\begin{thebibliography}{9999}
%
 \bibitem[BH02]{branson}
  \textsc{Branson, T.~, Hijazi, O.:}\quad
  \textit{Bochner--Weitzenb\"ock formulas associated with the
   Rarita--Schwinger operator,\ }
  \textrm{Internat. J. Math. 13 (2002), no. 2, 137--182.}
%
 \bibitem[CGH00]{cgh}
  \textsc{Calderbank, D.~, Gauduchon, P.~\& Herzlich, M.:}\quad
  \textit{Refined Kato inequalities and conformal weights in Riemannian
   geometry,\ }
  \textrm{J. Funct. Anal. 173 (2000), no. 1, 214--255.}
%
 \bibitem[DW]{dw}
  \textsc{Diemer, T.~\& Weingart, G.:}\quad
  \textit{private communication.}
%
 \bibitem[F76]{fegan}
  \textsc{Fegan, H.~D.:}\quad
  \textit{Conformally invariant first order differential operators,\ }
  \textrm{Quart. J. Math. Oxford (2) 27 (1976), no. 107, 371--378.}
%
 \bibitem[G91]{pg1}
 \textsc{Gauduchon, P.:}\quad
  \textit{Structures de Weyl et theoremes d'annulation sur une variete
   conforme autoduale, }
  \textrm{Ann. Scuola Norm. Sup. Pisa Cl. Sci. (4) {\bf 18}  (1991),
   no. 4, 563--629.}
%
 \bibitem[H04]{homma}
  \textsc{Homma, Y.:}\quad
  \textit{Casimir elements and Bochner identities on Riemannian manifolds,\ }
  \textrm{Prog. Math. Phys., 34, Birkh\"auser Boston, Boston, MA, 2004.}
%
 \bibitem[H05]{homma1}
  \textsc{Homma, Y.:}\quad
  \textit{Bochner identities for K\"ahlerian gradients,\ }
  \textrm{Math. Ann. 333 (2005), no. 1, 181--211.}
%
 \bibitem[H06]{homma2}
  \textsc{Homma, Y.:}\quad
  \textit{Bochner--Weitzenb\"ock formulas and curvature actions on
      Riemannian manifolds,\ }
  \textrm{Trans. Amer. Math. Soc. 358 (2006), no. 1, 87--114.}
%
 \bibitem[PP67]{perel}
  \textsc{Perelomov, A.~\& Popov, S.:}\quad
  \textit{Casimir Operators for Classical Groups,\ }
  \textrm{Soviet Math.~Dokl. 8 no. 3 (1967)  631---634.}
%
 \bibitem[S06]{uwe}
  \textsc{Semmelmann, U.:}\quad
  \textit{Killing forms on $\G_2$-- and $\Spin(7)$-- manifolds,\ }
  \textrm{ J. Geom. Phys. 56 (2006), no. 9,1752--1766.}
%
 \bibitem[SW]{mp}
  \textsc{Semmelmann, U.~\& Weingart, G.:}\quad
  \textit{Maple Program for Explicit Calculations for $\G_2$ and $\Spin(7)$,\ }
  \texttt{http://www.math.uni-bonn.de/~gw/g2spin7.mws.}
%
\bibitem[SW]{curvature}
  \textsc{Semmelmann, U.~\& Weingart, G.:}\quad
  \textit{On the Curvature Terms associated to Weitzenb\"ock Formulas,\ }
  \texttt{in preparation}
%
\end{thebibliography}
\end{document}